\documentclass[11pt]{article}
\usepackage[utf8]{inputenc}
\usepackage{amsfonts}
\usepackage{amsmath}
\usepackage{amsthm}
\usepackage{amssymb}
\usepackage{thmtools}
\usepackage{xcolor}
\usepackage{soul}
\usepackage[hypertexnames=false]{hyperref}
\usepackage[capitalize,noabbrev]{cleveref}
\usepackage{autonum}
\usepackage{comment}
\usepackage{geometry}     
\usepackage{microtype}

\usepackage{tikzsymbols}

\title{Density of shapes of periodic tori in arbitrary rank}
\author{
  Nguyen-Thi Dang\\
  {\normalsize Institut de Mathématique d'Orsay, Université Paris-Saclay}\\
  {\normalsize \href{mailto:nguyen-thi.dang@universite-paris-saclay.fr}{nguyen-thi.dang@universite-paris-saclay.fr}}
  \and
  Nihar Gargava\\
  {\normalsize Institut de Mathématique d'Orsay, Université Paris-Saclay}\\
  {\normalsize \href{mailto:nihar.gargava@universite-paris-saclay.fr}{nihar.gargava@universite-paris-saclay.fr}}
  \and
  Jialun Li\\
  {\normalsize CNRS and CMLS, École Polytechnique}\\
  {\normalsize \href{mailto:jialun.li@polytechnique.edu}{jialun.li@polytechnique.edu}}
}
\date{March 2025}

\newtheorem{theorem}{Theorem}[section]

\newtheorem{ques}[theorem]{Question}

\newtheorem{lemma}[theorem]{Lemma}
\newtheorem{proposition}[theorem]{Proposition}

\newtheorem{corol}[theorem]{Corollary}

\theoremstyle{definition}
\newtheorem{defin}[theorem]{Definition}
\crefname{defin}{Definition}{Definitions}
\Crefname{defin}{Definition}{Definitions}
\Crefname{corol}{Corollary}{Corollary}

\theoremstyle{remark}

\newtheorem{remark}[theorem]{Remark}

\newcommand{\sel}{\textup{\bf \scalebox{0.5}{(Sel)}}}
\DeclareMathOperator{\PGL}{PGL}
\DeclareMathOperator{\Hi}{H}
\DeclareMathOperator{\Ni}{N}
\DeclareMathOperator{\adj}{adj}
\DeclareMathOperator{\Id}{Id}

\DeclareMathOperator{\sgn}{sgn}
\DeclareMathOperator{\GL}{GL}
\DeclareMathOperator{\M}{M}
\DeclareMathOperator{\SL}{SL}
\DeclareMathOperator{\G}{G}

\DeclareMathOperator{\disc}{disc}
\DeclareMathOperator{\reg}{reg}
\DeclareMathOperator{\Gal}{Gal}

\DeclareMathOperator{\vol}{vol}
\DeclareMathOperator{\sys}{sys}

\DeclareMathOperator{\len}{len}

\newcommand{\calO}{\mathcal{O}}
\newcommand{\Z}{\mathbb{Z}}
\newcommand{\ZZ}{\mathbb{Z}}

\newcommand{\N}{\mathbb{N}}
\newcommand{\Q}{\mathbb{Q}}
\newcommand{\diag}{\mathrm{diag}}

\newcommand{\sh}{\mathrm{sh}}
\newcommand{\card}{\#}
\newcommand{\val}{\mathrm{v}}

\newcommand{\plus}{+}
\newcommand{\R}{\mathbb{R}}

\newcommand{\submat}[1]{\raisebox{0.2ex}{$\ulcorner$}\mkern-5mu#1}

\newcommand{\F}{\mathbb{F}}

\begin{document}

\maketitle
\begin{abstract}
Under the action of diagonal matrices $A \curvearrowright \PGL_{n+1}(\mathbb{R})/\PGL_{n+1}(\mathbb{Z})$, we consider the compact orbits. 
To each such orbit $F$, one can associate the shape point $\sh(F) \in \PGL_{n}(\mathbb{R})/\PGL_{n}(\mathbb{Z})$ corresponding to the lattice of periods of $F$ in the Lie algebra $\mathfrak{a} \simeq \mathbb{R}^{n}$ of $A$.

We show that the set of shapes of compact orbits is dense in $\PGL_{n}(\mathbb{R})/\PGL_{n}(\mathbb{Z})$ for all $n \geq 3$.
This extends our previous work for $n=2$ \cite{dangDensityShapesPeriodic2025}.
\end{abstract}

\section{Introduction}

Let $n \geq 2$ and let $A$ be the group of diagonal $(n+1) \times (n+1)$ matrices in $\PGL_{n+1}(\R)$.
In this paper, we are motivated by the action of $A$ on $\PGL_{n+1}(\mathbb{R})/\PGL_{n+1}(\mathbb{Z})$ and the compact orbits of this action.
We call such a compact orbit $F= A g \PGL_{n+1}(\mathbb{Z}) \subseteq \PGL_{n+1}(\mathbb{R})/\PGL_{n+1}(\mathbb{Z})$ a periodic torus (cf.~\cref{de:pt}).
It is classically known that there are countably many periodic tori. 

The prime geodesic theorem, first proven by Huber \cite{huberZurAnalytischenTheorie1959}, considers the problem for $n=1$ and says that for $T \xrightarrow[]{} \infty$, the counting
\begin{equation}
    \card \{F \subseteq \SL_{2}(\mathbb{R})/\SL_{2}(\mathbb{Z}) \text{ periodic torus } \mid \len(F) \leq T     \} \sim \frac{e^{T}}{T} .
\end{equation}

For a large parameter $T \xrightarrow[]{} \infty$, the following problems are the natural extensions of the prime geodesic theorem in the situation where $n \geq 2$.
\begin{enumerate}
\item Counting by listing the periodic tori $F \subseteq \PGL_{n+1}(\R)/\PGL_{n+1}(\mathbb{Z})$ in increasing order of the induced volume $\vol(F)$. This would involve estimating
\begin{equation}
\label{eq:vol_or}
N_{1}(T) = \card \{ F \subseteq \PGL_{n+1}(\mathbb{R})/\PGL_{n+1}(\mathbb{Z}) \text{ a periodic torus }\mid  \vol(F) \leq T \} 
\end{equation}
\item Counting by listing the periodic torus in increasing order of the systole length $\sys(F)$, i.e. the length of the smallest closed geodesic within $F$. That is, we want to estimate
\begin{equation}
\label{eq:sys_order}
N_{2}(T) = \card \{ F \subseteq \PGL_{n+1}(\mathbb{R})/\PGL_{n+1}(\mathbb{Z}) \text{ a periodic torus }\mid  \sys(F) \leq T \} 
\end{equation}
\end{enumerate}

It is clear that if the volume is small, the systole is small. This implies that $N_{1}(T) \leq N_{2}(c_{n} T^{\tfrac{1}{n}})$ for some $c_{n}>0$ as $T$ grows. However, the converse may or may not be true depending on the geometry of $F$. Indeed, when $A$ is connected, as a manifold $F$ is topologically isomorphic to the high-dimensional flat torus $\mathbb{R}^{n}/\mathbb{Z}^{n}$. Hence, it is a natural quest to consider the shapes of periodic tori in order to compare various generalizations of the prime geodesic theorem in higher rank.

For this, we define a shape of a periodic torus $F$ to be a point $\sh(F) \in \PGL_{n}(\mathbb{R})/\PGL_{n}(\mathbb{Z})$ as the following. Observe that the Lie algebra $\mathfrak{a}$ of the Lie group $A$ is isomorphic to $\mathbb{R}^{n}$. 
With one such isomorphism fixed, one can consider for the periodic torus $F = A g \PGL_{n+1}(\mathbb{Z})$, the subgroup of periods in $A$ (cf. ~\cref{eq:periods_defi}). 
Under the logarithm map in~\cref{eq:log_defi}, one associates this subgroup of periods to a lattice in $\mathfrak{a}$, and hence a point $\sh(F) \in \PGL_{n}(\mathbb{R})/\PGL_{n}(\mathbb{Z})$.

The following is our main theorem.
\begin{theorem}
    \label{th:main}
    
    The following set is dense in $\PGL_{n}(\mathbb{R})/ \PGL_{n}(\mathbb{Z})$ for $n \geq 2$.
     $$\{\sh(F)  \mid F \text{ periodic torus in }\PGL_{n+1}(\mathbb{R})/\PGL_{n+1}(\mathbb{Z}) \}.$$
     
     Here $\sh(\cdot)$ and periodic tori are as defined in~\cref{se:defi_of_shape}
\end{theorem}
Hence, this shows that bounding $N_{2}$ in terms of $N_{1}$ is not straightforward.
The notion of a shape of a periodic torus is not new. Our main theorem was previously a conjecture of David-Shapira for the case of $n=2$
in 
~\cite{davidShapesUnitLattices2017}. 
Their motivation to study the set of shapes was driven by a mathoverflow question \cite{56153}. The answer by user Andre Reznikov on this discussion mentions that the distribution of shapes is a question going back to Gromov and Margulis. Recently, there has been a renewed interest in understanding shape points for multiplicative unit lattices ~\cite{harron2026shapes,chari2025unit}. 

\subsection{A number theoretic version}

Let us now state~\cref{th:main} in a more number theoretic way. The main point is that considering the shapes of periodic tori is the same thing as considering shapes of multiplicative units of orders. Let us make this more precise.

Let $K$ be a totally real number field of degree $n+1$ and let $\mathcal{O} \subseteq K$ be an order (see~\cref{de:numberfield,de:order}).
By Dirichlet's unit theorem, the image of $\mathcal{O}^{\times}$ under
$$
    u \longmapsto (\log|\sigma_0(u)|,\dots,\log|\sigma_n(u)|)
$$
is a lattice in $\mathfrak{a}=\{(x_0,\dots,x_n)\in\mathbb{R}^{n+1}:\sum_{i=0}^{n}x_i=0\}$,
where  $(\sigma_0,\dots,\sigma_n)$ are an ordered set of all the embeddings $K \hookrightarrow  \R$.
Using our fixed identification $\mathfrak{a}\simeq\mathbb{R}^n$, denote by $\sh(\mathcal{O}) \in \PGL_{n}(\R)/\PGL_{n}(\mathbb{Z})$ the point corresponding to the image of $\mathcal{O}^{\times}$ in $\mathfrak{a}$ via the above map.
The equivalence that we will encounter in Proposition~\ref{pr:equivalence_of_shapes} can be used to obtain the following restatement of~\cref{th:main}.
\begin{theorem}
    \label{th:main_orders}
    For every $n\geq 2$, the shapes $\sh(\mathcal{O})$, as $K$ ranges over totally real number fields of degree $n+1$ with an ordered set of real embeddings and $\mathcal{O}\subseteq K$ ranges over all orders, are dense in $\PGL_n(\mathbb{R})/\PGL_n(\mathbb{Z})$.
\end{theorem}

This is of interest in number theory because it tells us that unit lattices can be used to approximate any lattice in number theory.
As mentioned in~\cite{chari2025unit}, this is of relevance in the theory of lattice-based cryptography.

See our paper \cite{dangDensityShapesPeriodic2025} for further information and motivation about the problem.

\subsection{About \texorpdfstring{$S$}{S}-arithmetic diagonal orbits}
The main goal of the current paper is two-fold.
\begin{enumerate}
    \item  To extend our main theorem to arbitrary rank with $n \geq 2$,
    \item  To make the previous version with $n=2$ a bit more clear conceptually by relying on the language of $S$-arithmetic diagonal flows.
\end{enumerate}

We emphasize that we see the second as a more significant contribution in our work than the first. 
Indeed, firstly we have since found alternative ways to prove~\cref{th:main} and we will talk about it shortly. 
Secondly, we think that the second part should be of independent interest in the community of $S$-arithmetic homogeneous dynamics.

In~\cref{se:lal}, we prove the following theorem. See~\cref{ss:sarith} for preliminaries on $S$-arithmetics.
\begin{theorem}[Informal]
    \label{th:main_sarith}
Let $n\geq 2$, let $S$ be a finite set of primes with $\card S\geq 2n+1$. 
Let $x\in(\mathbb{Q}_S)^{n-1}$ be a sufficiently generic point in the sense of~\cref{de:rfr}.
Write $S^{*}=S\cup\{\infty\}$. 
Denote $\G=\PGL_n$, and
$$
    \eta(x)=\begin{pmatrix}1&x\\0&I_{n-1}\end{pmatrix},
    \qquad \Hi=\{\diag(t,1,\dots,1):t\neq 0\}\subseteq \G.
$$
Then the diagonal orbit
$$
    \Hi(\mathbb{Q}_S)\eta(x)\G(\mathcal{O}_S)
    \subseteq \G(\mathbb{Q}_{S^{*}})/\G(\mathcal{O}_S)
$$
has dense image in $\G(\mathbb{R})/\G(\mathbb{Z})$ under the natural projection
$$
    \pi:\G(\mathbb{Q}_{S^{*}})/\G(\mathcal{O}_S)
    \longrightarrow \G(\mathbb{R})/\G(\mathbb{Z}).
$$
Here the finite-place factors are embedded with identity real component (See~\cref{eq:implicit_inclusion}).
\end{theorem}

In brief, the theorem tells us that a diagonal group action in $S$-arithmetic orbits produces dense sets generically. The interesting thing here is that this is for a predetermined starting point $x \G(\mathcal{O}_{S}) \in \G(\mathbb{Q}_{S^{*}})/\G(\mathcal{O}_{S})$ and not for a randomly chosen starting point. 
This theorem is the main dynamical input for~\cref{th:main} and occupies the entire second half of our paper.

The main inputs for~\cref{th:main_sarith} are the following:
\begin{enumerate}
\item  ``Banana lemma'' (\cref{le:banana}),
\item  The $\mathbb{Q}$-linear independence of logarithms of distinct primes (used in~\cref{le:diagonal_primes}),
\item Bourgain's exponential sum estimates on subgroups of multiplicative totient groups (\cref{th:bourgain}).
\end{enumerate}

 We don't think that the condition $\card S \geq 2n+1$ is optimal. 

\subsubsection{Can we go beyond density for \texorpdfstring{$S$}{S}-arithmetic diagonal orbits?}

Can one upgrade the density of $S$-arithmetic diagonal orbits in the sense of ~\cref{th:main_sarith} to some kind of equidistribution theorem?
It depends on how we arrange the points of course. 
One can come up with some box-like averages and try to formulate a question in~\cref{ss:problems} (see Question~\ref{qu:sarithmetic-box-equidistribution}).

We don't think this question is easy and we can try to explain this with a toy example.
\begin{ques}
    \label{qu:crt}

    Consider $a,b \in \mathbb{Z}$ such that $a \neq b$. Denote $f(m,n) \in \mathbb{Z}$ to be a solution of 
    \begin{align}
	f(m,n) & \equiv a\pmod{2^{m}}, \\
	f(m,n) & \equiv b\pmod{3^{n}},
    \end{align}
    Consider as $N \to \infty$ the set 
    \begin{equation}
    A_{N} = \left\{ \frac{f(m,n)}{2^{m} 3^{n}}\pmod{1}\right\}_{\substack{1 \leq m \leq N \\  1 \leq n \leq N}} \subseteq  \R/\mathbb{Z}.
    \end{equation}
   Do the normalized counting measures on $A_{N}$ weakly converge to
    the Lebesgue measure on $\R/\mathbb{Z}$?
\end{ques}
The question of density is much easier. 
The only thing to check is that $\overline{\bigcup_{N \geq 1} A_{N}}$ is invariant under $\times 2$ as well as $\times 3$ and is an infinite set.
Then Furstenberg's $\times 2,\times 3$ theorem~\cite{furstenberg1967disjointness} says that every infinite closed subset of $\R/\mathbb{Z}$ invariant under both maps is the whole circle. 

This problem of equidistribution in~\cref{th:main_sarith} is similar to Question~\ref{qu:crt}. Indeed, one can write the set $\bigcup_{N \geq 1}A_{N}$ in the following way. Consider the action
\begin{equation}
\times \tfrac{1}{2}, \times \tfrac{1}{3} \curvearrowright 
\frac{\R \times \mathbb{Q}_{2} \times \mathbb{Q}_{3}}{\mathbb{Z}[\tfrac{1}{2},\tfrac{1}{3}]}
\end{equation}
Take the point $(0,-a,-b)$ in $\R \times \mathbb{Q}_{2} \times \mathbb{Q}_{3}/ \mathbb{Z}[\tfrac{1}{2},\tfrac{1}{3}]$ and consider the orbit obtained by applying $2^{-m} 3^{-n}$ for $m,n \geq 1$.
Consider the natural projection map (see Corollary~\ref{co:sarithmetic_torus})
\begin{equation}
    \pi: 
\frac{\R \times \mathbb{Q}_{2} \times \mathbb{Q}_{3}}{\mathbb{Z}[\tfrac{1}{2},\tfrac{1}{3}]}
\rightarrow  \R/\mathbb{Z}.
\end{equation}
The image is precisely $\bigcup_{N \geq 1} A_{N}$. 
With this analogy, \cref{th:main_sarith} can be seen as a specialized version of Furstenberg's theorem in the $\PGL_{n}$ setting. Indeed, the main ingredient in Furstenberg's proof is the $\mathbb{Q}$-linear independence of $\log 2$ and $\log 3$.

\subsection{Outline of the proof and relation to the previous work}
To explain the differences between the new work and the previous work, let us briefly recall how the proof for $n=2$ worked in
\cite{dangDensityShapesPeriodic2025}.
\subsubsection{Outline of the proof for the case  \texorpdfstring{$n=2$}{n=2}}

To prove density of shapes of periodic tori in the rank $2$ case, we start with the polynomial $f(X) = X(X-a_{1})(X-a_{2}) - 1$ and vary the parameters $a_{1},a_{2}$ while maintaining $a_{1},a_{2},a_{2}-a_{1} > O(1)$. 
\begin{enumerate}
\item[{\bf Step 1:}]
It is clear that $f(X)$ is monic, irreducible and totally real. We take $K = \mathbb{Q}[X]/\langle f(X)\rangle$ and consider the order $\mathcal{O} = \mathbb{Z}[X] / \langle f(X) \rangle$.
We invoke the following tight inequality of Cusick \cite{cusickLower1983}:
\begin{equation}
    \vol\bigg( \frac{\mathfrak{a}}{ \log \mathcal{O}^{\times}} \bigg) \geq \tfrac{1}{16} [ \log \tfrac{1}{4} \disc f ]^{2} ,
\end{equation}
where $\log : \mathcal{O}^{\times} \rightarrow  \mathfrak{a}$ is the map given in \cref{eq:logmap}.
With this, one is able to conclude that $\mathcal{O}^{\times} = \pm \langle X, (X-a_{1})\rangle$.
Then we show that under the choice of these parameters, the shape points produced live in a compact set. Hence, for any sequence of parameters $(a_{1},a_{2})$, the shape points of $\log(\calO^\times)$ live in a compact set of shapes.
    
\item[{\bf Step 2:}]
We consider the shapes $\mathcal{O}_{(k)} = \mathbb{Z} + q \mathcal{O}$ for the shapes $\mathcal{O}$ considered above for $q = 2^{k_{1}} 3^{k_{2}} 5^{k_{3}}$ for $k = (k_{1},k_{2},k_{3})$ varying in a set $\{1,\dots,N\}^{3}$. 
We let $N \rightarrow \infty$.
Assuming that the parameters $a_{1},a_{2}$ have certain congruence conditions modulo $(2 \cdot 3 \cdot 5)^{N}$ 
for each $N$, we get that 
\begin{equation}
\mathcal{O}^{\times}_{(k)} = \{ \pm X^{m_{1}}(X-a_{1})^{m_{2}} \mid  (m_{1},m_{2}) \in \Lambda_{(k)} \subseteq \mathbb{Z}^{2} \},
\end{equation}
where $\Lambda_{(k)} \subseteq \mathbb{Z}^{2}$ is some lattice that can be explicitly written down in terms of 
congruence conditions modulo some powers of $2,3,5$. 

\item[{\bf Step 3:}]
The fact that the shapes of $\log \mathcal{O}^{\times}$ for the varying parameters $a_{1},a_{2}$ live in a compact set implies that for any sequence of parameters $(a_{1},a_{2})$, we can find a subsequence of parameters such that $\log \mathcal{O}^{\times}$ converges. Hence, using a diagonal extraction argument, we can assume that $\log \mathcal{O}^{\times}$ is a fixed point.

\item[{\bf Step 4:}]
It is therefore sufficient to show that the lattice $\{\Lambda_{(k)}\}_{k \leq \{1,\dots,N\}^{3}}$, when appropriately rescaled, are dense in $\PGL_{2}(\mathbb{R})/\PGL_{2}(\mathbb{Z})$.
We do this using the banana lemma. This involves showing that the 
closure of $\bigcup_{N \geq 1} \{\Lambda_{k}\}_{k \in \lbrace 1 , \ldots, N \rbrace^3}$ is invariant under the diagonal group action and contains a horocyclic piece.
\end{enumerate}

\subsubsection{New inputs for generalizing to \texorpdfstring{$n>2$}{n>2}}

When we try to generalize the previous result, we run into several difficulties.
Luckily, all of these obstacles can be resolved using some pleasant and well-known ideas from number theory.

\begin{enumerate}
\item To generalize {\bf Step 1},  we try to work with generalization of our previous choice of polynomial. We work with the ``split minus one'' polynomials 
\begin{equation}
    f(X) = X(X-a_{1}) \cdots (X-a_{n})-1
\end{equation}
and investigate the orders 
\begin{equation}
 \label{eq:defi_of_O}   
 \mathcal{O}^{\times} \text{ for } \mathcal{O} = \mathbb{Z}[X]/\langle f(X) \rangle.
\end{equation}
Clearly $X,X-a_{1},\dots,X-a_{n} \in \mathcal{O}^{\times}$ make a nice candidate for a multiplicative basis of units of $\mathcal{O}^{\times}$. But this idea runs into several problems which we list below.
    \begin{enumerate}
    \item  For the case of $n>2$, the generalization of the Cusick inequality for orders in high degree number fields did not seem to be written anywhere. 
	A theorem of Silverman \cite{silvermanInequalityRelatingRegulator1984}
    works for maximal orders in high degree. Hence, for our purpose, we rewrite in~\cref{ss:cusick} the proof of Cusick's inequality for orders in totally real number by adapting the version for cubic fields. 
    We thus obtain that there exists real positive numbers $c_n, \gamma_n$ such that for a totally real number field $K$ of degree $n+1$ {\em with no non-trivial subfields} 
    and any order $\mathcal{O} \subseteq K$
    \begin{equation}
	\label{eq:cusick_silv_intro}
 \vol \Big( \frac{\mathfrak{a}}{\log \mathcal{O}^{\times}}\Big) \geq c_n \log \big(\gamma_n \disc(\calO)\big)^{n}.
    \end{equation}
For a field $K$ that admits non-trivial subfields, the exponent on the $\log \disc \mathcal{O}$ term is not strong enough for our use.
    \item 
    We need to make the field $K$ associated with the ``split minus one'' polynomials not admit any non-trivial subfields.
    We are able to do this by choosing the parameters $a_{1},a_{2},\dots,a_{n} \in  \mathbb{Z}_{\geq 1}$ sufficiently separated. We show that for these generic parameters, the polynomial $f(X)$ generates a number field with no non-trivial subfield.
    
    This is consistent with the results of Widmer \cite{widmer2011number} that shows a similar result for a more natural model of random number fields. Widmer's work was a key step in Bhargava's proof of the Van der Warden conjecture 
    \cite{bhargava2025galois}.
    \item 
	For the $n=2$ case, the inclusion 
	$\pm \langle X, X-a_{1},X-a_{2}\rangle \subseteq \mathcal{O}^{\times}$ 
	is an equality and we are able to prove the shape of the sublattice is in a compact set of $\PGL_2(\R)/\PGL_2(\Z)$.
	In the general case, Cusick-Silverman inequality given in~\cref{eq:cusick_silv_intro}
	is not tight enough to show the same equality and we do not expect it. 
	All we can conclude is that 
	$\pm \langle X, X-a_{1},\dots,X-a_{n}\rangle \subseteq \mathcal{O}^{\times}$ is a finite index sublattice of bounded shape and whose index is bounded by some constant $C_{n} > 0$ that depends only on $n$ and not on $K$ or $\mathcal{O}$.
	Hence, there exists a matrix $[m_{ij}] = M \in \M_{n}(\mathbb{Z})$ with $|\det M| \leq C_{n}$ and a basis $\beta_{1},\dots,\beta_{n} \in \mathcal{O}^{\times}$ such that 
	\begin{equation}
	\prod_{i=1}^{n} \beta_{i}^{m_{ij}} = X - a_{j} \text{ for } j \in\{1,\dots,n\}.
	\end{equation}
	Up to changing the choice of $\beta_{i}$, we can assume that the matrix $M$ is among a finite set $\tilde{\mathcal{M}}$ of matrices in $\M_{n}(\mathbb{Z})$. 
    \end{enumerate}
    
    \item  Let $S$ be a finite set of prime numbers.
    We study the orders 
    \begin{equation}
    \label{eq:orders}
        \mathcal{O}_{(k)} = \mathbb{Z} + \prod_{p\in S} p^{k_p} \mathcal{O}\text{ for }1 \leq k_p \leq N.
    \end{equation}

    To generalize {\bf Step 2}, we caracterize the lattice $\Lambda_{(k)} \subseteq \mathbb{Z}^{n}$ such that 
    \begin{equation}
    \mathcal{O}_{(k)}^{\times} = \bigg\{  \pm \prod_{i=1}^{n} \beta_{i}^{m_{i}} \mid (m_{1},\dots,m_{n}) \in \Lambda_{(k)} \bigg\},
    \end{equation}
    using the basis of units $\beta_{1},\dots,\beta_{n}$ of $\mathcal{O}^{\times}$ chosen above.
    Understanding the congruence conditions defining the lattice $\Lambda_{(k)} \subseteq \mathbb{Z}^{n}$ as a function of $k$
    is somewhat more involved because of the introduction of the matrix $M$ above.
    We are able to somewhat simplify these computations by resorting to the following ideas from number theory.
    \begin{enumerate}
\item 
Consider $f^{\sel}(X)=(X+1)X^{n} - 1$. 
This polynomial is closely related to the Selmer polynomial\cite{nartEquationsTypeXn+aX+b1979, osadaGaloisGroupsPolynomials1987}.
and is a very well-studied polynomial in number theory.
	    
Denote $\mathcal{O}^{\sel} = \mathbb{Z}[X]/\langle f^{\sel}(X)\rangle$. Then the choice of primes is such that 
$p{\mathcal{O}^{\sel}} \subseteq \mathcal{O}^{\sel}$ is a prime ideal in $\mathcal{O}^{\sel}$ for all $p\in S$. 
In particular, this implies that the ring $\mathcal{O}^{\sel} \otimes \mathbb{Z}_{p}$ is the integral closure of $\mathbb{Z}_{p}$ in $K_{p}$. We use this simplification in~\cref{le:roots_of_unity}, for example.
        
\item 
Our choice of congruence conditions on $a_{1},\dots,a_{n}$ is such that for all prime $p\in S$,
\begin{equation}
f(X) \equiv f^{\sel}(X) \text{ or } f^{{\sel}}(X+1) \pmod{p^{N}},
\end{equation}
so that for every $1 \leq k_p \leq N$,
\begin{equation}
\frac{\mathcal{O}}{p^{k_{p}} \mathcal{O}} \simeq \frac{\mathcal{O}^{\sel}}{p^{k_{p}} \mathcal{O}^{\sel}}.
\end{equation}

\item We obtain a simpler description of $\mathcal{O}_{(k)}^{\times}$ in Proposition~\ref{prop:structure-suborders}.
The main technical difficulty in this part is that from the sublattice of known shape generated by $X-a_1, \ldots, X-a_n$, we can only recover $\beta_1^{\det M}, \ldots, \beta_n^{\det M}$.
Now to extract the suborder $(\Z + p^k \calO )^\times$, we can only do it up to some root of unity of $\frac{\mathcal{O}}{p^{k_p} \mathcal{O}}$.
Consequently, we require as parameter a generator of the roots of unity of $\calO^\sel \otimes \Z_p$ and deduce that $\Lambda_{(k)}$ is caracterized by congruence equations of two types : the first one deals with the non-torsional part (in a very similar fashion to the case $n=2$); the second one deals with the torsional part.  
The key points are that once the shape change $M$ and $p$ are fixed, the non-torsional equations are deterministic functions of the residue modulo $p^N$ of $a_1, \ldots, a_n$; there can be only $\varphi(p^{n+1}-1)^{n+1}$ possibilities for torsional part.
\end{enumerate}

\item Since there can only be finitely many shape change matrices and torsional equations, up to going into a subsequence for the parameters $(a_1,\ldots,a_n)$, we can assume the shape of $\calO^\times$ converges, the shape change matrix $M$ and the torsional equations are stationary.
This generalizes {\bf Step 3}
    
\item 
\cref{th:main_sarith} works to generalize the {\bf Step 4}.
Because of the reduction of {\bf Step 3}, we can now show that the shape points arrange themselves in an $S$-arithmetic diagonal orbits. To prove \cref{th:main_orders}, we show that the image of this orbit 
\begin{enumerate}
    \item is invariant by a sequence of diagonal elements of the form $\diag(a_1(i), a_2(i), \ldots, a_n(i)) \in \PGL_n(\R)$, with $a_1(i) \rightarrow + \infty$ and the $a_2(i) \geq \ldots \geq a_{n}(i)$ are bounded in absolute value; 
    {We prove invariance by the full diagonal group of $\PGL_n(\R)$ in \cref{prop:diag-inv}. }
    \item contains a positive measure piece of a horospherical orbit for the following horospherical group
    $$ U^+ = \lbrace g \in \PGL_n(\R) \; : \; a(i)^{-1}ga(i)\rightarrow I_n \text{ as } i \rightarrow \infty \rbrace .$$ 
    {For the subgoup 
    $$\Hi (\R) = \lbrace \diag(t, 1, \ldots, 1) \; : \; t \in \R_{> 0} \rbrace,$$ the relevant horospherical group is 
    $$\Ni(\R) = \bigg\lbrace \begin{pmatrix} 1 & x^{(1)} \cdots x^{(n-1)} \\ & I_{n-1}\end{pmatrix} \; : \; x \in \R^{n-1} \bigg\rbrace .$$
    We prove in \cref{prop-horos-piece} that the set of shapes contains an $\Ni(\R)$-orbit of starting point in $\PGL_n(\R)$. }
\end{enumerate}
 For $n=2$ we had used $3$ primes in \cite{dangDensityShapesPeriodic2025}. 
 For $n \geq 3$, our method requires $2n+1$ primes: the key restriction being given by \cref{le:diagonal_primes} for the diagonal invariance part.

A notable feature of our proof of density of the $S$-arithmetic diagonal orbit is our use of Bourgain's deep results about equidistribution of large multiplicative cosets in the rings $\mathbb{Z}/N\mathbb{Z}$ (cf.~\cref{ss:bourgain}). 
This equidistribution is a key step for showing the existence of a positive measure piece of horospherical orbit inside the closure of our diagonal orbit when projected to $\mathrm{PGL}_{n}(\mathbb{R})/\mathrm{PGL}_{n}(\mathbb{Z})$.

\end{enumerate}

\subsection{Some related problems}
\label{ss:problems}

We mention some open problems that we think would interest our readers. 
We think some of the techniques we have build in solving~\cref{th:main} could be useful in the following problems.

\begin{enumerate}
\item Here are some problems inspired from number theory. In our work, we exploit the fact that the unit group of an order $\mathcal{O}$ in a totally real number field $K$ of degree $n+1$ gives us some periodic torus $F$ and the 
     shape $\sh(F)$ of this periodic torus $F$ is, up to rescaling the shape $\sh(\mathcal{O})$ associated to the order $\mathcal{O}$, by which we mean the log embedding $\log \mathcal{O}^{\times} \subseteq \mathbb{R}^{n}$ of the unit lattice (cf.~\cref{ss:equivalence_of_shapes}).
      \begin{enumerate}
     \item In order to prove the density of shapes, we consider unit lattices $\log \mathcal{O}^{\times} \subseteq \mathbb{R}^{n}$ for a varying choice of orders $\mathcal{O} \subseteq K$ and a varying choice of $K$ among totally real subfields of degree $n+1$. 
     It is natural to ask if one can fix $K$ and prove density with $\log \mathcal{O}^{\times}$ as $\mathcal{O} \subseteq \mathcal{O}_{K}$ varies among suborders in the maximal order $\mathcal{O}_{K}$ of $K$.
     \item 
     It would be interesting if one could prove density of shapes of unit lattices $ \log \mathcal{O}^{\times}$ for $\mathcal{O}$ varying among {\em maximal} orders of totally real number fields. In fact, in the maximal order case, even showing that the shapes are unbounded or contain a non-empty interior in the closure would be interesting.
	 
     For $n=2$, the set of shapes of unit lattices of maximal orders was conjectured to be dense in \cite{chari2025unit} based on some numerical evidence.
         
     \item 
     The work of~\cite{harron2026shapes,chari2025unit} makes interesting observations about how the Galois group action of a number field $K$ can potentially induce symmetries on the shape points $\sh \mathcal{O}$. 
     It would be interesting if one could prove density of shapes in proper subvarieties of $\mathrm{PGL}_{n}(\mathbb{R})/\mathrm{PGL}_{n}(\mathbb{Z})$ using these Galois symmetries.
     \end{enumerate} 
 
\item Changing the space gives related questions.
 \begin{enumerate}
 \item 
 One could ask the same question by replacing the action of the diagonal group $A \curvearrowright \mathrm{PGL}_{n+1}(\mathbb{R})/\mathrm{PGL}_{n+1}(\mathbb{Z})$ by the action on some other homogeneous space. In each of the following examples, one could build a pleasant equivalence of shapes as in \cref{ss:equivalence_of_shapes}.
 \begin{enumerate}
 \item One could replace the group $\PGL_{n+1}(\Z)$ with some congruence subgroup. This imposes congruence conditions on the units in $\calO^{\times}$.
 \item 
One could consider shapes in compact homogeneous spaces. The dictionary for toral orbits in this case is also known \cite[Prop. 4.1]{ELMV}. 

 \item Another interesting case if that of Hilbert modular varieties. Take for example $\SL(2,\R) \times \SL(2,\R) / \SL_2(\mathbb{Z}[\sqrt{2}])$ where $\SL_2(\Z[\sqrt{2}])$ acts diagonally on each copy of $\SL_2(\R)$ via a different Galois embedding. In this case, one has a 2-dimensional torus acting on the homogeneous space and each periodic torus is related to an order in a quadratic extension $L$ of $\mathbb{Q}[\sqrt{2}]$ that has a $\Z[\sqrt{2}]$-module structure \footnote{It was reported to us by Menny Aka and Andreas Wieser that they have a partial result towards density in this case.}.
 \end{enumerate}
 \item Replacing $A$ by some other tori gives corresponds to orders in number fields with different signatures. For example the set of shapes of unit lattices of orders in cubic number fields with exactly one real embedding correspond to compact orbits of the action of the following group on $\SL(3,\mathbb{R})/\SL(3,\mathbb{Z})$:
 \begin{equation}
   \Bigg\{  \begin{pmatrix}
      a^{-2} & & \\
      & a \cos( \theta) & a \sin (\theta) \\ 
      & - a \sin(\theta)  & a \cos( \theta) \\
     \end{pmatrix} , a\in \R_{>0} , \theta \in \R \Bigg\}.
 \end{equation}

     \end{enumerate}
    
     \item In \cite{56153}, different types of equidistributions were conjectured based on numerical tests. Other than the volume or systole ordering mentioned in \cref{eq:vol_or} and (\ref{eq:sys_order}), one can also consider ordering the periodic tori with respect to the discriminant of the corresponding order and expect equidistribution.
     
     In particular, it is stated in \cite{56153} that Margulis expects that under the systole ordering used in \cref{eq:sys_order}, one does not expect equidistribution in the $n=2$ case, and instead expects an escape of mass to the cusp. We do not have any intuition for why such a non-equidistribution is expected and we have not seen any numerical evidence to confirm or deny such a behaviour in the systole ordering \footnote{We were recently informed by Radu Toma that he has generated some numerical examples about this escape of mass phenomenon with the systole ordering.}. 

     \item 
	 Can one upgrade~\cref{th:main_sarith} to equidistribution along exponent boxes?
     Fix $n\geq 2$, a finite set of primes $S$ with $\card S\geq 2n+1$. Take a rational, strongly full rank point $x\in\mathbb{Q}_S^{n-1}$ as in~\cref{de:rfr}. 
     Write $G=\PGL_n$.
     \begin{ques}\label{qu:sarithmetic-box-equidistribution}
     For every such fixed $S$ and $x$, do the probability measures
     $$
         \mu_N=\frac{1}{(N+1)^{\card S}}
         \sum_{k\in\{0,\ldots,N\}^{S}}
         \delta_{\pi\big(a_S(-k)\eta(x)G(\mathcal{O}_S)\big)}
     $$
     converge to a nice limit measure as $N\to\infty$, where $a_S(k)$ is defined in~\cref{eq:defi_of_ask}?
     \end{ques}
 The starting point is deterministic and fixed. We only require it to satisfy the stated arithmetic hypotheses. 
 With these hypotheses, we know by~\cref{th:main_sarith} that the image forms a dense set so one can ask for a continuous limit distribution. A partial escape of mass could be plausible.

\end{enumerate}

\subsection{Acknowledgements}

We have also had fruitful discussions about the project with Anurag Rao, Menny Aka, Andreas Wieser, Radu Toma, Matthew de-Courcy Ireland and Noy Soffer Aranov.

NG acknowledges funding from the SNSF Postdoc.Mobility grant 225437, ``Random Geometries with Arithmetic Constraints'', and from ERC Advanced Grant 101096550, ``Integrating Spectral and Geometric data on Moduli Space''.

We thank Fudan University for their kind hospitality. Additionally, 
we would like to thank the fruit shop on Zhengmin road in Shanghai for their mangoes, durians, pipas and kumquats.

\paragraph{Role of LLMs/AI}
There was no role of AI in the process of proof generation for this paper. The ideas are based on our original pre-AI work for the cubic case \cite{dangDensityShapesPeriodic2025} and most of the work done was during the summer of 2025 in Paris and spring 2026 in Shanghai. 
The preprint was delayed because we found some gaps in a previous version and used our human brains to fix the gap\footnote{The gap was we were trying to show diagonal invariance in \cref{ss:diagonal} which had some problems when we worked only $n+1$ primes. This fix was to use $2n+1$ primes instead.}.

That being said, we have used AI for language, spelling, grammar and proof audits. 
We also ran various LLMs to try to answer as many questions as possible in \cref{ss:problems}.
They made considerable progress on several of those questions but we chose not to report it in this article, to keep AI-assisted work separate from AI-free work.

\section{Preliminaries}

\subsection{Number fields and orders}

\begin{defin}
\label{de:numberfield}
A number field $K$ is defined as the field 
$\mathbb{Q}[X]/\langle f(X) \rangle$
for some irreducible polynomial $f(X) \in \mathbb{Q}[X]$. 
The degree of $K$ is defined to be the degree of $f$. When $f$ is totally split over $\mathbb{R}$, we say $K$ is totally real.
\end{defin}

Throughout this text, we will work with number field such that $\dim_{\mathbb{Q}} K =  n+1$. For a given irreducible polynomial $f(X) \in \mathbb{Q}[X]$, we say that $K = \mathbb{Q}[X]/\langle f(X) \rangle$ is the number field associated to $f$.
 
 \begin{remark}
  Each totally split irreducible polynomial $f(X) \in \Q[X]$ of degree $n+1$ comes with a natural ordered set of embeddings $\{\sigma_0, \cdots , \sigma_{n}\} $ where $\sigma_i:K \to \R$ and $K = \Q[X]/\langle f(X) \rangle$. This is because we can assign $\sigma_i: X \mapsto \alpha_i$ where $\alpha_0 < \cdots < \alpha_n$ are the roots of the polynomial $f(X)$ (there are no repetitions among roots). 
  
  This way, every totally real polynomial comes with an ordered set of embeddings.

 \end{remark}

\begin{defin}
\label{de:order}
Given a number field $K$, we define an order $\mathcal{O}$ to be a subring in $K$ that contains $1$, is a finitely generated $\mathbb{Z}$-module and satisfies $\mathcal{O} \cdot \mathbb{Q} = K$.
\end{defin}

For an irreducible monic polynomial in $\mathbb{Z}[X]$ that is totally split polynomial over $\mathbb{R}$, we know that $\mathbb{Z}[X]/ \langle f(X) \rangle \subseteq \mathbb{Q}[X]/ \langle f(X) \rangle$ is an order. We will call this construction the order associated to $f(X)$.

\subsection{Shapes}
\label{se:defi_of_shape}
\subsubsection{Shapes in homogeneous dynamics}

Let ${A}$ be the group of diagonal matrices in $\PGL_{n+1}(\mathbb{R})$ for $n \geq 2$.
More precisely,
\begin{equation}
A=\left\{\operatorname{diag}(a_0,\ldots,a_n)\cdot \mathbb R^\times
\;\middle|\;a_0,\ldots,a_n\in\mathbb R^\times\right\} \subseteq \PGL_{n+1}(\R).
\end{equation}
Consider 
the action of $A$ on $\PGL_{n+1}(\mathbb{R})/\PGL_{n+1}(\mathbb{Z})$. 
Observe that the Lie algebra $\mathfrak{a}$ of {the identity component of} ${A}$ is 
\begin{equation}
\mathfrak{a}= \{ (x_{0},\dots,x_{n})\in \mathbb{R}^{n+1} \mid  \sum_{i=0}^{n} x_{i} = 0\}.
\end{equation}

\begin{defin}
\label{de:pt}
For $n \geq 2$, we define a periodic torus (under ${A}$-action) $F$
to be an orbit
$$
F={A}g\PGL_{n+1}(\mathbb{Z}),\qquad
 g\PGL_{n+1}(\mathbb{Z})\in\PGL_{n+1}(\mathbb{R})/\PGL_{n+1}(\mathbb{Z}),
$$
such that $F\subseteq\PGL_{n+1}(\mathbb{R})/\PGL_{n+1}(\mathbb{Z})$ is compact.

We denote 
\begin{equation}
  \mathcal{P}(F) = \{a \in {A} \mid a g \PGL_{n+1}(\mathbb{Z}) = g \PGL_{n+1}(\mathbb{Z})\}
  \label{eq:periods_defi}
\end{equation}
to be the periods of $F$.
\end{defin}

We fix a linear bijection 
\begin{equation}
   \mathfrak{a} \simeq \mathbb{R}^{n} \label{eq:lin_bijection} 
\end{equation}
for the rest of this paper.
Consider the logarithm map $\log : {A} \rightarrow  \mathfrak{a}$ given by 
\begin{equation}
\label{eq:log_defi}
\diag(a_{0},\dots,a_{n}) \R^{\times} \mapsto   ( \log |a_{i}|- \tfrac{1}{n+1} \textstyle\sum_{j=0}^{{n}}\log| a_{j}| )_{i=0}^{n} .
\end{equation}
Because of the bijection in~\cref{eq:lin_bijection}, one can get a map
$
    \log: {A} \rightarrow  \mathbb{R}^{n}
    $.

Given a periodic torus $F$, we define the shape point associated to $F$ to be the following.
\begin{defin}
	
Given a periodic torus $F$, one considers the subgroup $\mathcal{P}(F) \subseteq {A}$. We know that $\mathcal{P}(F) \subseteq {A}$ is discrete and cocompact.
Therefore, under the logarithm map $\log \mathcal{P}(F)  \subseteq \mathbb{R}^{n}$ is a lattice. Hence, for some $h \in \GL_{n}(\mathbb{R})$, we know that $\log \mathcal{P}(F) =  h \cdot \mathbb{Z}^{n}$. We denote $\sh(F) = h \PGL_{n}(\mathbb{Z}) \in \PGL_{n}(\mathbb{R})/\PGL_{n}(\mathbb{Z})$.
\end{defin}
This $\sh(F)$ is the homothety class of $\log \mathcal{P}(F)$.

\subsubsection{Shapes in number theory}

Given a totally real number field $K$ of degree $n+1$ and an ordered set of embeddings $\sigma_{0},\dots, \sigma_{n} : K \rightarrow  \mathbb{R}$, one can construct the logarithm map 
\begin{align}
\log : K^{\times}  & \rightarrow  \mathbb{R}^{n+1}, \\ 
x & \mapsto  ( \log |\sigma_{i}(x)  |)_{i=0}^{n}
\label{eq:logmap}
\end{align}

Given $u \in \mathcal{O}^{\times}$, observe that $\prod_{i=0}^{n}|\sigma_{i}(u)|=1$. 
This is because the multiplicative norm of $u$ must be $\pm 1$.
Hence, for any order $\mathcal{O} \subseteq K$, we can map the unit group $\mathcal{O}^{\times}$ to $\mathfrak{a}$. The image is a lattice because of the following theorem of Dirichlet.
\begin{theorem}
    
	{\bf (Dirichlet's unit theorem)}
	
Let $K$ be a totally real number field of degree $n+1$ and let $\mathcal{O} \subseteq K$ be an order. Then, there exists multiplicatively independent units $u_{1},\dots,u_{n} \in \mathcal{O}^{\times}$ such that 
\begin{equation}
    \mathcal{O}^{\times} = \{ \pm u_{1}^{m_{1}}\dots u_{n}^{m_{n}}  \mid (m_{1},\dots,m_{n}) \in \mathbb{Z}^{n}  \}.
\end{equation}
That is, $\mathcal{O}^{\times}/ \{\pm 1\}$ is a free Abelian group of rank $n$.
\end{theorem}

We define shape of an order as the following
\begin{defin}
	
Given a totally real number field $K$ of degree $n+1$, an ordered set of embeddings $\sigma_{0},\dots,\sigma_{n}$ and an order $\mathcal{O} \subseteq K$, one denotes $\sh(\mathcal{O})$ to be the point $h \PGL_{n}(\mathbb{Z}) \in  \PGL_{n}(\mathbb{R})/\PGL_{n}(\mathbb{Z})$  such that the lattice $\log \mathcal{O}^{\times} =  h \cdot \mathbb{Z}^{n}$, where $h \in \GL_n({R})$ is some matrix and $\log : \mathcal{O}^{\times}\rightarrow \mathbb{R}^{n}$ is the map in~\cref{eq:logmap} after identifying $\mathfrak{a} \simeq \mathbb{R}^{n}$ with respect to the bijection in~\cref{eq:lin_bijection}.
\end{defin}

\subsubsection{Equivalence of shapes}
\label{ss:equivalence_of_shapes}

\begin{proposition}
\label{pr:equivalence_of_shapes}

  The following two subsets of $\PGL_n(\mathbb{R})/\PGL_n(\mathbb{Z})$ are equal.
  \begin{enumerate}
      \item $\{\sh(F)  \mid F \text{ periodic torus in }\PGL_{n+1}(\mathbb{R})/\PGL_{n+1}(\mathbb{Z}) \}$,
      \item $\{\sh(\mathcal{O})  \mid K \text{ totally real field with an ordered set of embeddings}, \mathcal{O} \subseteq K \text{ some order}, \deg K = n+1\}$.
  \end{enumerate}
\end{proposition}
\begin{proof}
Let us show that the second set is contained in the first set. Let $K$ be a totally real number field of degree $n+1$. Let $\sigma_{0},\dots,\sigma_{n}:K \rightarrow  \mathbb{R}$ be $n+1$ real embeddings and let $\mathcal{O} \subseteq K$ be an order.
Under the identification
\begin{align}
    K &  \rightarrow  \mathbb{R}^{n+1} \\ 
    x & \mapsto  (\sigma_{i}(x))_{i=0}^{n}
\end{align}
Under this map, we know that $\mathcal{O}$ maps to some lattice $c \cdot g \cdot \mathbb{Z}^{n+1}$ for some $g \in \GL_{n+1}(\mathbb{R})$ and $c \in \mathbb{R} \setminus \{0\}$. We claim that the point $F = {A} g \PGL_{n+1}(\mathbb{Z})$ is a periodic tori whose shape $\sh(F) = \sh(\mathcal{O})$. 
Indeed, note that the following is an isomorphism of Lie groups.
\begin{align}
    (K \otimes \mathbb{R} )^{\times}/ \R^{\times}  &  \xrightarrow[]{} {A} \\
    x & \mapsto   \diag( \sigma_{0}(x), \sigma_{1}(x),\dots,\sigma_{n}(x))
\end{align}
Then, under this isomorphism, we see that the periods $\mathcal{P}(F)$ defined in~\cref{eq:periods_defi},
are $\mathcal{O}^{\times}/\{ \pm 1\}$ since
\begin{equation}
\{x \in K\otimes \mathbb{R} \mid x \mathcal{O} = \mathcal{O}\} = \mathcal{O}^{\times}.
\end{equation}
Therefore, we know that $\sh(\mathcal{O}) = \sh(F)$.

Now let us show that the first set is contained in the second set.
Let $F = {A} g \PGL_{n+1}(\mathbb{Z})$ be a periodic tori. Consider the set of periods $\mathcal{P}(F)$ as in~\cref{eq:periods_defi}. We know by the compactness of $F$ that there must exist some multiplicatively independent elements $u_{1},\dots,u_{n} \in \mathcal{P}(F)$ such that 
\begin{equation}
    \mathcal{P}(F) = \{  u_{1}^{m_{1}} \dots u_{n}^{m_{n}}  \mid  (m_{1},\dots,m_{n}) \in \mathbb{Z}^{n} \}.
\end{equation}
Observe that $g^{-1} u_{i} g \in \PGL_{n+1}(\mathbb{Z})$. 
Hence we can chose a representative of $u_{i} \in A$ in $\GL_{n+1}(\R)$
that has a characteristic polynomial that is monic and lies in $\mathbb{Z}[X]$. Let us also denote those representatives by $u_{i}$.
Hence, each diagonal entry of this $u_{i}$ is an algebraic integer.
Consider the $\mathbb{Q}$-algebra
\begin{equation}
K = \mathbb{Q}[u_{1},\dots,u_{n}] \subseteq \M_{n+1}(\mathbb{R}).
\end{equation}
We claim that $K$ is a totally real number field of degree $n+1$.
Since $K \otimes \mathbb{R} = \mathbb{R}[u_{1},\dots,u_{n}] \subseteq \text{ diagonal matrices in }\M_{n+1}(\mathbb{R})$, 
one knows that $K \simeq K_{1} \oplus \dots  \oplus K_{r}$ for some totall real number feilds $K_{1},\dots,K_{r}$.
We know that $\deg K = \sum_{} \deg K_{i} \leq n+1$. 
The $u_{i}$ are multiplicatively independent units in $K$ so they generate an $\mathbb{Z}$-module of rank $n$.
This means that $\sum_{i=1}^{r}(\deg K_{i} - 1) \geq n$, because each $K_{i}$ can contribute to a rank of at most $\deg K_{i}-1$. Hence $r=1$ and we show the claim.

One constructs the embeddings $\sigma_{0},\dots,\sigma_{n}: K \rightarrow  \mathbb{R}$ by mapping to the various diagonal entries of $K \subseteq \M_{n+1}(\mathbb{R})$. Finally, we can construct $\mathcal{O}$ via
\begin{equation}
    \mathcal{O} = \{ x \in K \mid x  g \mathbb{Z}^{n+1} \subseteq g \mathbb{Z}^{n+1}\}.
\end{equation}
With this order $\calO$, we can check that 
\begin{equation}
\mathcal{P}(F)  = \mathcal{O}^{\times}/ \pm 1.
\end{equation}
\end{proof}

\subsection{\texorpdfstring{$S$}{S}-arithmetics}
\label{ss:sarith}
\subsubsection{Local completions}

For each prime number $p \in \mathbb{Z}_{\geq 2}$, we induce an absolute value $|\cdot|_{p}$ on $\mathbb{Q}$ according to the following principle.

\begin{equation}
    | x|_p =  p^{-\val_{p}(x)},
\end{equation}
where for $x \in \mathbb{Z} \setminus \{0\}$, we define
\begin{equation}
    \val_{p}(x)  = \max\{ n \in \mathbb{Z} \mid p^{-n}x \in \mathbb{Z} \},
\end{equation}
and for $x  = a/b \in \mathbb{Q}$, $\val_{p}(x) = \val_{p}(a) - \val_{p}(b)$. We also define $\val_{p}(0)=\infty$

One can then check that $(x,y) \mapsto  |x-y|_{p}$ is a metric on $\mathbb{Q}$. We denote $\mathbb{Q}_{p}$ to be the metric completion of $\mathbb{Q}$ with respect to this metric and $\mathbb{Z}_{p}$ to be the metric completion of $\mathbb{Z}$. Note that 
\begin{equation}
    \mathbb{Z}_{p} =\{x  \in \mathbb{Q}_{p} \mid |x|_{p} \leq 1 \}.
\end{equation}
Under the metric space topology, we get that $\mathbb{Z}_{p}$ is compact.

For each prime $p$, the above completion $\mathbb{Q}_p$ is called a non-Archimedean completion.
The usual absolute value $x \mapsto x\sgn(x)$ on $\mathbb{Q}$ is called the Archimedean absolute value, which we shall denote as 
$|\cdot |_{\infty}$. We denote $\mathbb{Q}_{\infty}  = \mathbb{R}$. 

The natural embedding of $\mathbb{Q} \rightarrow \mathbb{Q}_{v}$ is dense for any place $v$.

\begin{remark}
A classical theorem of Ostrowski says that any absolute on $\mathbb{Q}$ is one of $|\cdot|_{p}$ for $p = \infty$ or $2,3,5,\dots$ etc.
\end{remark}

\subsubsection{\texorpdfstring{$S$}{S}-integers}
We let $S$ be a finite set of places not including $\infty$. Denote $S^{*} = S \cup \{\infty \}$. We define 
\begin{align}
    \mathbb{Z}_{S} & =  \prod_{p \in S} \mathbb{Z}_{p} ,\\
    \mathbb{Q}_{S^{*}} & =  \prod_{p \in S^{*}} \mathbb{Q}_{p}.
\end{align}
We finally define the $S$-integers to be
$$ \calO_S = \Z [\{\tfrac{1}{p} \}_{p \in S}].$$
\begin{lemma}
$$  \mathcal{O}_{S} =\{x  \in  \mathbb{Q} \mid  x \in \mathbb{Z}_{q} \text{ for all primes } q \notin S ,  q\neq \infty\}$$
\end{lemma}
\begin{proof}
A generating set of $\calO_S$ as a $\mathbb{Z}$-module is $\prod_{p\in S}p^{-a_p}$ where $(a_p)_{p\in S} \in \Z_{\geq 0}^S$.

$\subset$ For every other prime $q \notin S \cup \lbrace \infty \rbrace$, each prime in $S$ is invertible in $\Z/q\Z$. 
So in particular, these generators of $\calO_S$ is in $\Z_q$ for every other prime $q\notin S$.

$\supset$ Let $x\in \mathbb{Q}$ such that $x\in \Z_q$ for every other prime $q\notin S$.
We write $x=\frac{a}{b}$ with $a$ and $b$ coprime. 
The prime factorisation of $b$ can only contain primes in $S$ otherwise, there would be a prime $q\notin S$ for which $b$ is not invertible in $\Z_q$ which is a contradiction. 
\end{proof}

One observes that there is a natural diagonal inclusion 
\begin{equation}
    \mathcal{O}_{S} \rightarrow  \mathbb{Q}_{S^{*}},
\end{equation}
which uses the inclusion $\mathbb{Q} \rightarrow \mathbb{Q}_{p}$ for each $p \in S^{*}$. 
One also has the inclusion
\begin{equation}
    \mathbb{Z}_{S} \times \{0\}  \subseteq  \mathbb{Q}_{S^{*}}.
\end{equation}
When the context is clear, we will abuse our notations slightly and write $\mathbb{Z}_{S} \subseteq \mathbb{Q}_{S^{*}}$ without mentioning the $0$ and the infinite place.
\begin{lemma}
  \label{le:cocompact_inters}
  The additive quotient $\mathbb{Q}_{S^{*}}/\mathcal{O}_{S}$ is compact.
\end{lemma}
\begin{proof}

    Any element of $\mathbb{Q}_{S^{*}}$ can be denoted by 
    \begin{equation}
        ( x_{p}+z_{p})_{p \in S} \times (y),
    \end{equation}
    where $x_{p} \in \mathbb{Z}[\tfrac{1}{p}], z_{p} \in \mathbb{Z}_{p}$ and $y \in \mathbb{R}$. 
    We act additively with $-\sum_{q\in S} x_q$, which is a rational number in $\calO_S$.
    Then for every $p\in S$ we have
    $$x_p+z_p - \sum_{q\in S} x_q = z_p - \sum_{q\in S\setminus \lbrace p\rbrace} x_q.$$
    Now every $q\in S \setminus \lbrace p\rbrace $ is prime and coprime with $p$, so $\Z[\tfrac{1}{q}]\subset \Z_p$ when $q\neq p$.
    We denote by $z'_p = z_p - \sum_{q\in S\setminus \lbrace p\rbrace} x_q$ and deduce that
    $$ ( x_{p}+z_{p})_{p \in S} \times (y) \equiv ( z'_p )_{p\in S} \times (y-\sum_{p\in S} x_p)  \pmod{\calO_S}, $$
    where $z'_p\in \Z_p$ for every $p\in S$.
    Finally, by acting additively with elements of $\mathbb{Z}$, we can assume that $y \in (0,1]$

Hence up to $\mathcal{O}_{S}$-action, every $x \in \mathbb{Q}_{S^{*}}$ has a representative in $\mathbb{Z}_{S} \times (0,1]$ which is (pre)compact.
\end{proof}

\begin{corol}
    \label{co:sarithmetic_torus}
    The above lemma proves that 
    \begin{equation}
        \mathbb{Z}_{S} \backslash \mathbb{Q}_{S^{*}} / \mathcal{O}_{S} \simeq \mathbb{R}/\mathbb{Z}.
    \end{equation}
\end{corol}

\subsection{Simultaneous local approximation of rationals}

\label{ss:defi_of_approx}

We defined three different levels of approximations of elements of $\mathbb{Q}_{S}$ by elements of $\mathbb{Q}$.

\begin{defin}
    \begin{enumerate}

        \item 
	    For a prime $p$, $x \in \mathbb{Q}_{p}$ and $k \in \mathbb{Z}$ some exponent, we denote $[x]_{k} \in \mathbb{Z}[\tfrac{1}{p}]$ some rational 
    such that $$x-[x]_{k} \in p^{k} \mathbb{Z}_{p}.$$
    \item 
	For a finite set of primes $S \subseteq \mathbb{P}$, given $(x_{p})_{p \in S} \in \mathbb{Q}_{S}$ and a tuple of integers $ (k_{p})_{p \in S} \in \mathbb{Z}^{S}$, we write $[x]_{k} \in \mathcal{O}_{S}$ some rational number such that 
	\begin{equation}
	    x_{p} - [x]_{k} \in p^{k_{p}} \mathbb{Z}_{p} \text{ for every } p \in S.
	\end{equation}

    \item For a tuple $x=(x^{(i)})_{1 \leq i \leq n} \in \mathbb{Q}_{S}^{n}$ and a tuple of exponents $(k_{p})_{p \in S} \in \mathbb{Z}^{S}$, 
    one denotes by 
$[x]_{k}  = ([x]_{k}^{(i)})_{1 \leq i \leq n} \in \mathcal{O}_{S}^{n} $ to be a solution of 
	    \begin{equation}
		x^{(i)}_{p} - [x]_{k}^{(i)} \in p^{k_{p}} \mathbb{Z}_{p}\text{ for every } i \in \{1,\dots,n\}\text{ and } p \in S.
	    \end{equation}
    \end{enumerate}
\end{defin}

While using this notation, one should always note that $[x]_{k}$ is not uniquely defined. However, different choice only differs by $\prod_{p\in S}p^{k_p} \Z $.

Here's a lemma about where to expect such approximants.
\begin{lemma}
  \label{le:where_to_find}
  Suppose $x \in \mathbb{Z}_{S}$ and $k \in \mathbb{Z}_{\geq 0}^{S}$. 
  Then, one can choose $[x]_{k} \in \{0, 1 ,\dots ,N-1\}$ where $N= \prod_{p \in S} p^{k_{p}}$.
\end{lemma}

\subsection{\texorpdfstring{$S$}{S}-arithmetic homogeneous space}

Now we will establish homogeneous spaces. For every $\mathbb{Z}$-algebra $R$, we denote by
\begin{equation}
    \PGL_{n}(R) = \GL_{n}(R) / R^{\times}.
\end{equation}
For notational convenience, we will simply write $\G = \PGL_{n}$. 
In our writing, we will often write elements of $\PGL_{n}(\mathbb{R})$ and $\PGL_{n}(\mathbb{Q})$ as matrices in $\GL_{n}(\mathbb{R})$ and $\GL_{n}(\mathbb{Q})$ respectively 
with the quotient by $\mathbb{R}^{\times}, \mathbb{Q}^{\times}$ implicitly understood.

There is a diagonal inclusion as before, given by 
\begin{equation}
    \G(\mathcal{O}_{S}) \rightarrow \G(\mathbb{Q}_{S}^{*}).
\end{equation}
We also assume the inclusion $\G(\mathbb{Z}_{S}) \rightarrow  \G(\mathbb{Q}_{S^{*}})$ to be 
\begin{equation}
    \G(\mathbb{Z}_{S}) \times \{ 1_{\G(\mathbb{R})}\} \rightarrow  \G(\mathbb{Q}_{S^{*}}).
\end{equation}

Then, by an $S$-arithmetic homogeneous space, we mean the following quotient: 
\begin{equation}
    \G(\mathbb{Q}_{S^{*}})/ \G( \mathcal{O}_{S}) .
\end{equation}

The following proposition gives us a surjection $ \G(\mathbb{Q}_{S^{*}})/ \G( \mathcal{O}_{S}) \to  \G(\mathbb{R})/\G(\mathbb{Z}).$
It is a classical fact and we state it without a proof. This is a standard consequence of strong approximation; see for example \cite[Section 7]{PlatonovRapinchuk1994}.
\begin{proposition}
  \label{le:surjective_S}
  There is a natural isomorphism
  \begin{equation}
      \G( \mathbb{Z}_{S})  \backslash \G(\mathbb{Q}_{S^{*}}) / \G(\mathcal{O}_{S}) \simeq  \G(\mathbb{R})/\G(\mathbb{Z})
  \end{equation}
\end{proposition}

\begin{remark}
\cref{le:surjective_S} is a generalization of \cref{co:sarithmetic_torus} when $\G(R)$ is defined to be the additive group of a $\mathbb{Z}$-algebra $R$.
\end{remark}

Henceforth, we will denote $\pi$ to be the natural projection
\begin{equation}
       \pi : \G(\mathbb{Q}_{S^{*}}) / \G(\mathcal{O}_{S}) \rightarrow   \G(\mathbb{R})/\G(\mathbb{Z}).
\end{equation}
We will also use the following inclusion implicitly:
\begin{align}
    \label{eq:implicit_inclusion}
    \iota: \G(\mathbb{Q}_{S}) & \rightarrow  \G(\mathbb{Q}_{S^{*}})\\
    x & \mapsto  (\Id,x).
\end{align}
\subsubsection{Unipotent and diagonal element}
\label{ss:diagonal}

\begin{enumerate}
    \item[(a)]{\bf Unipotent elements}
    
    We are interested in the unipotent group $\Ni \subseteq \G$ given by matrices of the form.
\begin{equation}
    \begin{pmatrix}
	    1 & * \\
	       & I_{n-1}  
    \end{pmatrix}    
\end{equation}
Let us define for any ring $R$ and 
$x = (x^{(i)})_{1 \leq i \leq {n-1}} \in R^{n-1}$, we denote by $\eta(x)$ the matrix 
\begin{equation}
      \eta (x) = 
      \begin{pmatrix}
	  1 & x^{(1)} \cdots x^{(n-1)} \\
	   &  I_{n-1}  
      \end{pmatrix}
      \in 
\Ni(R)
 \subseteq
      \M_{n}(R).
      \label{eq:defi_of_eta}
\end{equation}
We will use this notation for $R= \mathcal{O}_{S}, \mathbb{Z}_{S}, \mathbb{Q}_{S}, \mathbb{R},\mathbb{Q}, \mathbb{Z}$.

\item[(b)]{\bf Diagonal elements}

We denote by $\Hi \subseteq \G$ the subgroup of matrices in $\G$ given as
\begin{equation}
    \begin{pmatrix}
	    * &         \\
	       &  I_{n-1}  
    \end{pmatrix}    
\end{equation}
For any ring $R$, for every $t \in R$ we consider the matrix
\begin{equation}
    a(t) = \diag(t,1,\dots,1) \in \Hi(R).
\end{equation}
We will use this notation for $R= \mathcal{O}_{S}, \mathbb{Z}_{S}, \mathbb{Q}_{S}, \mathbb{R},\mathbb{Q}, \mathbb{Z}$, etc.

Notationally, we will deal with the following types of diagonal elements in $\Hi(\mathbb{Q}_{S})$ for which we designate a special notation.
For every $k = (k_{p})_{p \in S} \in \mathbb{Z}^{S}$, we consider the following diagonal matrix
\begin{equation}
    a_{S}(k) = 
    \big(\; a(p^{k_{p}})\; \big)_{p \in S} \in  \Hi(\mathbb{Q}_{S}).
    \label{eq:defi_of_ask}
\end{equation}

\end{enumerate}

\subsection{Hecke congruence groups and equality of orbits}

We first establish the following lemma, which we will strengthen subsequently.
\begin{lemma}
  \label{le:easy_hecke}
  
Let $A, B \in \M_{n}(\mathbb{Z})$ be such that $\det A = \det B = \gamma \in \mathbb{Z}$ and $ \gamma \neq 0$. 
Then, $A \equiv B \pmod{{\gamma}} \implies  A \cdot \G(\mathbb{Z}) = B \cdot \G(\mathbb{Z})$.
\end{lemma}
\begin{proof}
    We need to check that $A^{-1} B \in \G(\mathbb{Z})$. This is implied if each entry of $ A^{\adj} B $ is divisible by $\gamma$ since $A^{-1} = \gamma^{-1} A^{\adj}$.
    But $A^{\adj} B \equiv A^{\adj} A \equiv \det A \cdot I_n \equiv 0 \pmod{\gamma}$ so we are done.
\end{proof}

With the same assumption on $A,B$, if $A,B \in \M_{n}(\mathbb{Z})$ are such that $A \cdot \G(\mathbb{Z}) = B \cdot \G(\mathbb{Z})$, then the elementary divisors of $A,B$ are the same up to possible reordering and multiplication by $\pm 1$. 
Hence, there is a choice of elementary divisors $\gamma_{1},\gamma_{2},\dots,\gamma_{n} \in \mathbb{Z}$ of $A$ and $B$ and there exists $g_{1},g_{2} \in \G(\mathbb{Z})$ such that
$$\left\lbrace
\begin{array}{cl}
  A \cdot \G(\mathbb{Z}) &= g_{1} \diag(\gamma_{1},\gamma_{2},\dots , \gamma_{n} ) \cdot \G(\mathbb{Z})     \\
  B \cdot \G(\mathbb{Z}) &= g_{2} \diag(\gamma_{1},\gamma_{2},\dots , \gamma_{n} ) \cdot \G(\mathbb{Z})    
\end{array}.
\right.$$

The following lemma is a strengthened version of \cref{le:easy_hecke} 
in a special case when $\gamma_{1}=\gamma$ and $\gamma_{2},\dots,\gamma_{n} = 1$.
\begin{lemma}
  \label{le:hecke_neighbour}
  Let $N \geq 1$.
  Then $g_{1},g_{2} \in M_{n}(\mathbb{Z})$ with $\det g_1=\det g_2 \neq 0$ satisfy
  \begin{equation}
      \label{eq:equal_cosets}
      g_{1} a(N) \G(\mathbb{Z}) = g_{2} a(N) \G(\mathbb{Z})
  \end{equation}
  if and only if
  \begin{equation}
    \frac{1}{\det g_2}g_{2}^{\adj} g_{1}  \in
    a(N) \G(\mathbb{Z})  a(N^{-1}).
  \end{equation}
\end{lemma}
\begin{proof}
    So the first line is equivalent to $g_2^{\adj}g_1 a(N) \G(\Z) = \det( g_2) a(N) \G(\Z)$.
    Hence since $\G(\Z)$ is a group, equivalent to the required conclusion.
\end{proof}

This inspires us to state the following definition.

\begin{defin}
    \label{de:hecke_subgroup}
Let $N\geq 1$. 
    We define the congruence subgroup $\Gamma(N)$ to be 
    \begin{equation}
        \Gamma(N) = 
    a(N) \G(\mathbb{Z})  a(N)^{-1}  \cap \G(\mathbb{Z}) \subseteq \G(\mathbb{Z}).
    \end{equation}

    Equivalently, 
    \begin{equation}
        \Gamma(N) =  \{ 
            [a_{ij}]_{1 \leq i,j \leq n}
            \in \G(\mathbb{Z})
            \mid  a_{12} \equiv \dots \equiv a_{1n} \equiv 0 \pmod{N}
        \}.
    \end{equation}
\end{defin}

The following Corollary is a variation of the previous lemmas where we want to understand equality of orbits while isolating the interaction of $g_{1},g_{2} \in \M_{n}(\mathbb{Z})$ with a prime $p$ coprime to the determinant $\gamma = \det g_{1} = \det g_{2}$. 
We will use this Corollary in~\cref{se:lal}.
Let us introduce notations for minors of matrices and submatrices.
Let $g\in \mathrm{M}_n(\R)$. For every $1 \leq i,j \leq n$, denote by $\Delta_{i,j}(g)$ the determinant of the submatrix of $g$ obtained by removing the $i-$th line and the $j-$th column.
Denote by $\submat g$ the submatrix of $g$ obtained by removing the first line and first column.
\begin{corol}
    Let $x = (x^{(i)})_{1 \leq i\leq n-1} \in \mathbb{Z}^{n-1}$.
    Let $g_{1},g_{0}\in \M_{n}(\mathbb{Z})$ and $\gamma \in \Z\setminus \lbrace 0\rbrace$ such that $\det g_{1} = \det g_{0} = \gamma$. 
    Let $\eta(x) \in \G(\mathbb{Z})$
    where the function $\eta$ is defined in \cref{eq:defi_of_eta}.
    Let $p \in \N$ be a prime 
    and $l \geq 1$. Assume also that
    $$
g_1^{\mathrm{adj}}\eta(x)g_0\equiv0\pmod{\gamma}.
$$

Then $g_{1} a({p^l}) \G(\mathbb{Z}) = \eta(x) g_{0} a(p^l) \G(\mathbb{Z})$
    if and only if the following conditions hold.
\begin{enumerate}
        \item[1.)] $\Delta_{1,1}(g_1) \big(x \; \submat g_0+ v(g_0) \big)+
    \big( (-1)^{2+1} \Delta_{2,1}(g_1) ,\cdots, (-1)^{n+1}  \Delta_{n,1}(g_1) \big) \; \submat g_0 \equiv 0 \pmod{p^l\gamma} $,
    \item[2.)] $ g_{1}^{\adj} \eta(x) g_{0}\equiv 0\pmod{\gamma}$,
\end{enumerate}
    where $v(g_0)=\big( (g_0)_{1,j} \big)_{2 \leq j \leq n} \in (\Z^{n-1})^t$.
    \label{co:difficult_hecke}
\end{corol}
\begin{proof}
    By \cref{le:hecke_neighbour}, having $g_{1} a({p^l}) \G(\mathbb{Z}) = \eta(x) g_{0} a(p^l)\G(\Z)$ is equivalent to
    $$ \frac{1}{\det(g_1)}g_{1}^{\adj} \eta(x) g_{0} \in \Gamma(p^l). $$
   By matrix computation, one checks that
    	$$ g_{1}^{\adj} \eta(x) g_{0}=\begin{pmatrix}
    	    * & \Delta_{1,1}(g_1) (x\; \submat g_0 + v(g_0) )+
    \big( (-1)^{2+1} \Delta_{2,1}(g_1) ,\cdots, (-1)^{n+1}  \Delta_{n,1}(g_1) \big)\; \submat g_0 \\
    	    * & *
    	\end{pmatrix} $$
    	Hence, it is equivalent to
	 	$$\Delta_{1,1}(g_1) \big(x \; \submat g_0+ v(g_0) \big)+
    \big( (-1)^{2+1} \Delta_{2,1}(g_1) ,\cdots, (-1)^{n+1}  \Delta_{n,1}(g_1) \big) \; \submat g_0\equiv 0 \pmod{p^l\gamma}. $$
      and $$ g_{1}^{\adj} \eta(x) g_{0}\equiv 0\pmod{\gamma} $$ 
\end{proof}

\begin{corol}\label{co:less_difficult_hecke}
With the same assumption as in \cref{co:difficult_hecke}, the following three conditions are sufficient to have $g_{1} a({p^l}) \G(\mathbb{Z}) = \eta(x) g_{0} a(p^l)\G(\Z).$
\begin{itemize}
    \item[1.)] $g_1-g_0 \equiv 0 \pmod{\gamma}$,
    \item[2.)] $x \equiv 0 \pmod{\gamma}$
	   \item[3.)] $\Delta_{1,1}(g_1) \big(x \; \submat g_0+ v(g_0) \big)+
    \big( (-1)^{2+1} \Delta_{2,1}(g_1) ,\cdots, (-1)^{n+1}  \Delta_{n,1}(g_1) \big) \; \submat g_0\equiv 0 \pmod{p^{l+r}}$.

\end{itemize}
Here $r=\val_{p}(\gamma)$.
\end{corol}
\begin{proof}
The second condition of \cref{co:difficult_hecke} that $g_1^{adj}\eta(x)g_0 \equiv 0 \pmod{\gamma}$ can be verified by assuming $ g_1\equiv g_0,\ x\equiv 0\pmod{\gamma}$, using \cref{le:easy_hecke}. 
For the first condition, the divisibility of the entries by $\gamma$ is also satisfied since $\gamma$ divides every entry of $g_1^{\adj}\eta(x)g_0$.
\end{proof}

\subsection{Equidistribution of multiplicative cosets of totient rings}
\label{ss:bourgain}

For the horospherical piece, we need some discussion about equidistribution of multiplicative subgroups of $\Z/N\Z$. 

 \begin{defin}\label{defin-equid}
 Let $(X_k)_{k\in \N}$ be a sequence of sets such that $X_k \subset \Z/N_k\Z$ for every $k\in \N$ where $(N_k)$ is an increasing sequence of integers. 
    We say $(X_k)$ equidistributes as $k\rightarrow\infty$ if the uniform measure of the rescaled set $\frac{X_k}{N_k}$ converges weakly to the Lebesgue measure on $\R/\Z$. 
    Equivalently, for every non-zero $\xi \in \Z$,
    $$ \frac{1}{\vert X_k \vert } \sum_{x \in X_k}  
    \exp\bigg( \frac{2i\pi \xi x  }{N_k} \bigg) \xrightarrow[k \rightarrow + \infty]{} 0. $$
 \end{defin}
The main tool is the following estimate of Bourgain.
\begin{theorem}[\cite{bourgain_exponential_2005}]\label{th:bourgain}
Let $N\in \Z_+$ and $H<(\Z/N\Z)^\times$ be a multiplicative subgroup such that $\vert H\vert>N^\delta$
 where $\delta>0$ is arbitrary.
 Then one has the exponential sum estimate
 $$  \underset{\xi \in (\Z/N\Z)^\times}{\mathrm{max}} \Bigg\vert \frac{1}{\vert H\vert}  \sum_{x \in H} \exp\bigg(\frac{2i\pi  \xi x}{N}\bigg) \Bigg\vert  < c N^{-\varepsilon}, \quad \text{where } \varepsilon=\varepsilon(\delta) >0.$$
\end{theorem}
\noindent The original article does not specify what the constant $c>0$ depends on. 
It can only depend on $\delta$ but not on $H$ nor on $N$ otherwise the statement would be empty.
Note that when $N$ is large enough, we can take a slightly smaller exponent $\varepsilon(\delta)$ which will simplify this constant to $c=1$.

The following Corollary of the previous theorem is the main equidistribution statement we use in the rest of the article.
\begin{corol}\label{corol-bourgain}
Let $(N_k)_{k\in \N}$ be an increasing sequence of integers and $(H_k)_{k\in \N}$ be a sequence of multiplicative subgroups $H_k<(\Z/N_k\Z)^\times$ such that $\vert H_k \vert > N_k^\delta$ where $\delta>0$ is arbitrary.

Then for every non-zero integer $\alpha\in \Z$ (which need not be coprime with the $N_k$), all sequences $(\lambda_k)_{k\in \N}$ and $(\mu_k)_{k\in \N}$ where $\lambda_k \in (\Z/N_k \Z)^\times$ and $\mu_k \in \Z/N_k\Z$, there exists $k_\alpha\in \N$ such that for all $k\geq k_\alpha$ large enough,
$$ \frac{1}{\vert  H_k\vert} \; \underset{\xi \in (\Z/N_k\Z)^\times}{\max} \bigg\vert \sum_{x \in  H_k} \exp\bigg( \frac{2i\pi \xi (\mu_k + \alpha \lambda_k x)}{N_k} \bigg) \bigg\vert < c\frac{|\alpha|^\varepsilon}{N_k^\varepsilon} \text{ where } \varepsilon=\varepsilon(\delta/2)>0 .$$
As a consequence, 
$ \mu_k + \alpha \lambda_k H_k $ equidistributes as $k \rightarrow \infty.$ 
In particular, $H_k$ also equidistributes as $k \rightarrow \infty.$
\end{corol}

\begin{proof}
Denote by $\alpha_k = \gcd(\alpha, N_k)$. Then for all $\xi \in (\Z/N_k\Z)^\times$,
$$ \sum_{x\in H_k} \exp\bigg( \frac{2 i \pi \xi(\mu_k + \alpha \lambda_k x)}{N_k} \bigg) =
\exp\bigg( \frac{2 i \pi \xi \mu_k}{N_k} \bigg) \sum_{x\in H_k} \exp\bigg( \frac{2 i \pi \xi (\alpha/\alpha_k) \lambda_k x }{N_k/\alpha_k} \bigg).$$
In absolute value, the term containing $\mu_k$ vanishes. 
Denote by $N_k'=N_k/\alpha_k$ and by $H_k'$ the projection of $H_k$ in $\Z/N_k'\Z$.
Note that $\xi (\alpha/\alpha_k) \lambda_k$ is coprime with $N_k$ hence also coprime with $N_k'$.
It remains to find the correct lower bound for $\vert H_k'\vert$.

The map $(\Z/N_k\Z)^\times \twoheadrightarrow (\Z/N_k'\Z)^\times$ is surjective and each fiber has cardinality $\varphi(N_k)/\varphi(N_k')$.
Consequently $$ \vert H_k'\vert \geq \vert H_k\vert \frac{\varphi (N_k')}{\varphi(N_k) } > N_k^\delta \frac{\varphi (N_k')}{\varphi(N_k) } = (N_k')^\delta \alpha_k^\delta \frac{\varphi (N_k')}{\varphi(N_k) }. $$
By properties of Euler totient function, 
$$\frac{\varphi (N_k')}{\varphi(N_k) } = \frac{\varphi (N_k') /N'_k}{\varphi(N_k)/N_k } \; \frac{N_k'}{N_k } = \frac{\prod_{p \mid N'_k} (1-1/p) }{ \prod_{p \mid N_k} (1-1/p) } \; \frac{1}{\alpha_k} \geq  \frac{1}{\alpha_k} .$$
Hence $\vert H_k'\vert > (N_k')^\delta (\alpha_k)^{\delta-1}. $
Now if $\delta \geq 1$, then we directly obtain the correct bound on $\vert H_k' \vert$ to apply \cref{th:bourgain}.
Otherwise, $\alpha_k = \gcd(\alpha,N_k) \leq |\alpha|$, so we bound below by $(N'_k)^\delta |\alpha|^{\delta-1}.$ 
Asymptotically as $k$ goes to infinity, $N_k'\asymp N_k \rightarrow + \infty$, hence the lower bound is $\vert H_k'\vert > (N_k')^{ \delta/2}$ for $k\geq k_\alpha$ large enough.

We apply \cref{th:bourgain}, there exists $\varepsilon(\delta/2)>0$ such that for all $\xi \in (\Z/N_k\Z)^\times,$ for $k\geq k_\alpha$,
\begin{equation}
    \begin{split}
   \frac{1}{\vert H_k\vert } \bigg\vert \sum_{x\in H_k} \exp\bigg( \frac{2 i \pi \xi(\mu_k + \alpha \lambda_k x)}{N_k}& \bigg) \bigg\vert \leq \frac{1}{\vert H_k \vert } \;  \underset{\xi' \in (\Z/N_k'\Z)^\times}{ \max } 
\bigg\vert  \frac{\vert H_k\vert}{\vert H_k'\vert}  \sum_{ x \in H_k'} \exp \bigg( \frac{2i\pi \xi' x}{N_k'} \bigg)  \bigg\vert \\
& \leq \frac{1}{\vert H_k' \vert} \; \underset{\xi' \in (\Z/N_k'\Z)^\times}{ \max } 
\bigg\vert \sum_{ x \in H_k'} \exp \bigg( \frac{2i\pi \xi' x}{N_k'} \bigg)  \bigg\vert < c(N_k')^{-\varepsilon} = c\frac{|\alpha|^\varepsilon}{N_k^\varepsilon}.
    \end{split}
\end{equation}
The image $\mu_k+\alpha\lambda_kH_k$ has $\vert H_k'\vert$ elements, and $\vert H_k\vert/\vert H_k'\vert\leq\varphi(N_k)/\varphi(N_k')\leq|\alpha|$.

Let us now deduce equidistribution of $\mu_k+\alpha \lambda_k H_k$.
Fix a non-zero $\zeta \in \Z \setminus \lbrace 0 \rbrace$.
Then we apply the main estimate for the sequence $(\zeta \mu_k + \zeta \alpha \lambda_k H_k)_{k\in \N}$ and for $k\geq k_{\alpha \zeta}$ large enough
$$\frac{1}{\vert H_k\vert}\bigg\vert\sum_{x\in H_k}\exp\bigg(\frac{2i\pi\zeta(\mu_k+\alpha\lambda_kx)}{N_k}\bigg)\bigg\vert\leq c\frac{|\alpha\zeta|^{\varepsilon}}{N_k^\varepsilon}\xrightarrow[k\to+\infty]{}0.$$
Hence the equidistribution towards the Lebesgue measure as $k$ goes to infinity.
\end{proof}

\newcommand{\calP}{\mathcal{P}}
\newcommand{\calM}{\mathcal{M}}

\section{Setup for the ``shapes machine''}

\subsection{``Split minus one'' polynomials}

We consider a class of polynomials that we choose to call ``split minus one'' polynomials. These polynomials have been studied before by \cite{casselsProductInhomogeneousLinear1952} and \cite{ankenyNoteClassNumbersAlgebraic1956}.

\begin{defin}
	\label{de:splitminusone}
	Let $n\geq 2$ and let $a_0,a_1, \cdots ,a_{n} \in \Z$ be distinct integers such that $a_{j+1}-a_j> 2$ for all $0 \leq j \leq n-1$. We denote the following polynomial to be ``split minus one'' with parameters $a_0,\cdots,a_{n}$.
\begin{equation}
f(X):= (X-a_0) \cdots (X- a_{n}) - 1.
\end{equation}
\end{defin}
Note that without loss of generality, we can assume that $a_0=0 \leq a_1 \leq \ldots \leq a_n$ by ordering the integer coefficients and doing a change of variable $X \mapsto X-a_0$.

The polynomial enjoys the following properties.
\begin{lemma}\label{lem:irred}
	Consider $f$ a "split minus one" polynomial
\begin{enumerate}
    \item All the roots of $f(X)$ are real and distinct.
    \item
    The polynomial $f(X)$ is irreducible.
    \item Let $\alpha_{0}<\dots<\alpha_{n} \in \mathbb{R}$ be the roots of $f(X)$. Then, one has for each $i \in \{0,\dots,n\}$ that
	    \begin{equation}
	        |\alpha_{i}- a_{i}| \leq \prod_{\substack{j=0 \\ j \neq i}}^{n} (|a_{j}-a_{i}|-1)^{-1} 
		\label{eq:inequality_of_roots}
	    \end{equation}
\end{enumerate}
\end{lemma}
\begin{proof}
\begin{enumerate}
    \item For each $a_i$, one sees that $f(a_i)=-1$ but at least one of $f(a_i \pm 1) >0$. Using that two consecutive roots are at a distance of at least 2, this implies that a real root exists by intermediate value theorem. 
    
    \item 
    
    Because $f(X)$ is monic, the only way it can factorize is
    $$
    f(X)=g(X)h(X),\qquad 0<\deg h(X),\deg g(X)<n+1,
    $$
    for some monic $g,h\in\mathbb{Z}[X]$.  
    
    By substituting $a_i$, we know that $g(a_i) h(a_i) = -1$ so $(g(a_i),h(a_i)) =  (\pm 1, \mp 1)$. This means that the polynomial $j(X) = g(X) + h(X)$ is zero at all the $a_i$. Hence $j(X)=0$, since the degree of $j$ is less than $n+1$.
    
    So $h(X) = - g(X)$ must hold and $f(X) = -h(X)^2$. But since all the roots of $f(X)$ are distinct, this is impossible.
    \item We know that $f(\alpha_{i})=0$. A priori, one also knows that $|\alpha_{i}-a_{i}| \leq 1$ since the roots are between $a_i$ and $a_{i} + 1$ or $a_{i}-1$ for each ${i} \in\{0,\dots,n\}$.
	    Since we have 
	        $f(X)=\prod_{i=0}^{n}(X-\alpha_{i})$, one has for each $i \in \{0,\dots,n\}$ that $|f(a_{i})|=1$. So
	    \begin{equation}
	        |\alpha_{i}- a_{i}| = \prod_{\substack{j=0 \\ j \neq i}}^{n} |\alpha_{j}-a_{i}|^{-1}.
	    \end{equation}
	    Then, replacing $\alpha_{j}$ with $a_{j}$ can at most make an error of $1$, we estimate this as
	    \begin{equation}
	        |\alpha_{i}- a_{i}| \leq \prod_{\substack{j=0 \\ j \neq i}}^{n} (|a_{j}-a_{i}|-1)^{-1}.
	    \end{equation}

\end{enumerate}
\end{proof}

\subsection{Higher dimensional Cusick-Silverman inequality}

\label{ss:cusick}

For any number field $K$ denote by $\reg(K)$ the regulator, the covolume of $\log$ of the units of maximal order of $K$. Let $\disc(K)$ be the square of the covolume of the maximal order of $K$ in $K \otimes \R \simeq \R^{n+1}$.
\begin{proposition}\label{prop:cusick-high-rank}
Let $K$ be a totally real number field of degree $n+1\geq 2$.
Assume $K$ has no non-trivial subfield.
Then 
$$\reg(K)\geq c_n \log(\gamma_n \disc(K))^{n},$$
for positive constants $c_n$ and $\gamma_n$ depending only on $n$.
If $\log(\gamma_n\disc K)<0$, replace the right-hand side by $0$.
\end{proposition}

This is the inequality from \cite{silvermanInequalityRelatingRegulator1984}. The proof is written for maximal orders and does not directly generalize to non-maximal orders. 
However, Cusick's proof \cite{cusickRegulatorSpectrumTotally1991,cusickLower1983} written for the case of cubic fields works for non-maximal orders.
Let us write a proof of the Cusick-Silverman {inequality} for orders here.  

\begin{proposition}\label{prop:cusick}
Let $K$ be a totally real number field of degree $n+1\geq 3$.
Assume $K$ has no non-trivial subfield.
Let $\calO \subseteq \calO_K$ be any order in $K$.

Then  
$$
\label{eq:ineq_cusik}
 \reg(\calO)
\geq c_n \log \big(\gamma_n \disc(\calO)\big)^{n}
,$$
for some constants $c_n, \gamma_n$ depending on $n$ but not on $K$ or $\calO$. If $\log(\gamma_{n} \disc \mathcal{O}) < 0$, we replace the right side of~\cref{eq:ineq_cusik} by 0.
\end{proposition}
\begin{remark}
\begin{enumerate}
\label{re:cusick-tight}
\noindent
    \item More precisely, one can take $\gamma_n,c_{n}$ as some explicit functions of $n$.
Note that these {constants} are definitely not optimal, as we compute worse constants than Cusick for the cubic case.
\item While dealing with the cubic case version of the same problem, the authors previously {leveraged} the fact that the constants $\gamma_n,c_n$ were tight enough to exactly pinpoint the shape of the unit lattice of the relevant orders \cite[Proposition 3.4]{dangDensityShapesPeriodic2025} . 
In the higher rank case, having a Cusick-Silverman inequality of the same quality is hopeless and this introduces some complications in the proof. The discussion about roots of unity in \cref{sec:roots-of-unity} is an artefact of this complication. 

\end{enumerate}

\end{remark}

\begin{proof}
    Let $\epsilon \in \calO^{\times}$ be some unit which we will choose more precisely later
    and let $\sigma_j : K \rightarrow \mathbb{R}$ be the real {embeddings} for $j\in\{0,1,\cdots, n\}$. Then, we observe that 
    $$
    \label{eq:vandermonde}
    { \disc(\epsilon)} 
    = \prod_{j > l} (\sigma_j(\epsilon) - \sigma_l(\epsilon) ) ^2.
    $$
    Here, the discriminant is the square of the covolume of $\mathbb{Z}[\epsilon] \subseteq \calO$ within $\calO \otimes \R = K \otimes \R$.
    For the non-torsion unit $\epsilon$ chosen below, $\mathbb{Z}[\epsilon]$ must be full-rank because by our assumption we get $K=\mathbb{Q}(\epsilon)$. The equality in ~\cref{eq:vandermonde} follows from the standard equality on Vandermonde matrices.
    
      Let $\tau$ be a permutation on $\{0,1\cdots,n\}$ such that $|\sigma_{\tau(0)}(\epsilon)|\leq \cdots\leq  |\sigma_{\tau(n)}(\epsilon)|$. We note that
    $$ \disc(\calO) 
    \leq { \disc(\epsilon)} = \prod_{ j > l} \left( 1- \frac{\sigma_{\tau(l)}(\epsilon)}{\sigma_{\tau(j)}(\epsilon)} \right)^2 \sigma_{\tau(j)}(\epsilon)^2 
    \leq 2^{n(n+1)} \prod_j \sigma_{\tau(j)}(\epsilon)^{2 j} .
    $$
    For the rightmost inequality, we use that $\left( 1- \frac{\sigma_{\tau(l)}(\epsilon)}{\sigma_{\tau(j)}(\epsilon)} \right)^2 \leq 2^2$ and that there are $\frac{n(n+1)}{2}$ such terms to bound above.
    So we get an inequality 
   $$
   \log \disc(\calO)
   \leq n(n+1)\log 2 \; + \; \sum_{j } 2j \log | \sigma_{\tau(j)}( \epsilon) |. 
   $$
   Then by Cauchy--Schwarz and $\tau$ being a permutation,
   $$ \log \disc(\calO)\leq n(n+1)\log 2 \; + \; (\sum_{j } (2j)^2)^{1/2}(\sum_j (\log | \sigma_{j}( \epsilon) |)^2)^{1/2}.$$

   Let $\epsilon_1,\cdots, \epsilon_n$ be the generators of $\calO^{\times}/\{\pm 1\}$ as a $\Z$-module. 
   Then $\epsilon = \pm \epsilon_1^{x_1} \cdots \epsilon_{n}^{x_{n}},$ and for every $1 \leq j \leq n$ we write 
   $$\log\vert \sigma_j(\epsilon)\vert = \sum_{k=1}^{n} x_k \log\vert \sigma_j (\epsilon_k) \vert= \big(h_{\calO} x\big)_j$$
   where $x=(x_1,...,x_{n})\in \Z^{n}$ and $(h_{\calO})_{i,k}=\log\vert \sigma_i (\epsilon_k) \vert$ is an $(n+1)\times n$ matrix. The regulator $\reg (\calO)$ is the absolute value of any $n \times n$ minor of $h_{\calO}$. 

   We 
   get 
    \begin{align}
       {(\sum_j (\log | \sigma_{j}( \epsilon) |)^2)^{1/2}}
& \leq  \Vert h_{\calO}x\Vert_2 .
   \end{align}
   where $\Vert \; \Vert_2$ denotes the Euclidean norm of $\R^{n+1}$.
{By Minkowski's theorem applied to the lattice $h_{\calO}\Z^n\subseteq\mathfrak a$ (of covolume $\sqrt{n+1}\reg(\calO)$),
there is $x\in\Z^n\setminus\{0\}$ with $h_\calO x$ in a suitable ball around $0$ and hence we can claim}
   \begin{align}
      {(\sum_j (\log | \sigma_{j}( \epsilon) |)^2)^{1/2}}
& {\leq  c'_n  \; \reg(\calO)^{\frac{1}{n}}}, \text{ for some } c'_{n} > 0.
   \end{align}
 We deduce that 
 $$ \left( \max\left\{ \log \Big(2^{-n(n+1)} \disc (\calO) \Big) , 0    \right\} \right)^{n}
 \leq (c'_{n} )^{n}\left( n   \sum_{j=1}^{n} 4j^2 \right)^{n/2}  \; \reg(\calO) .$$
 Note that $\sum_{j=1}^{n} 4j^2\leq 4 \frac{ n^3}{3} +O(n^{2})$.
 {This gives the inequality with} $\gamma_n = 2^{-n(n+1)}$ and $c_n >0$ that only depends on $n$.
\end{proof}

\subsection{Avoiding subfields}

The higher dimensional Cusick-Silverman inequality requires that we work with fields that have no {non-trivial} subfields.
In this section we will establish a {criterion} to show when non-trivial subfields do not arise for $K = \mathbb{Q}[X]/\langle f(X) \rangle$ for 
$f(X) = (X-a_0)(X-a_1) \cdots (X-a_{n}) - 1.$
As remarked before, without loss of generality, we can assume that $a_0=0 \leq a_1 \leq \ldots \leq a_n$ by ordering the integer coefficients and doing a change of variable $X \mapsto X-a_0$.

First, we prove {a lemma of the following type}.

\begin{lemma}
	\label{le:partition_galois}
	Let $f(X) \in \mathbb{Z}[X]$ be an irreducible polynomial of degree { $n+1$} and let $K$ be the number field $K=\mathbb{Q}[X]/\langle f(X) \rangle$. 
Let $\{\alpha_{0},\dots,\alpha_{n}\} \subseteq \overline{\mathbb{Q}}$ be the roots of $f$ and let $\tilde{K}=\mathbb{Q}(\alpha_{0},\dots,\alpha_{n})$ be the splitting field of $f$. Then, the following two statements are equivalent.
\begin{enumerate}

    \item There is a non-trivial subfield $\mathbb{Q} \subsetneq L \subsetneq K$ contained in $K$ with $\deg L = m \geq 2$, $m \mid n+1$.
    \item There exists {a reindexing of the roots} of $f$ in an $m\times \frac{n+1}{m}$ matrix such that the partition
		    $$\{\alpha_{i}\}_{i=0}^{n} = \bigsqcup_{l=0}^{m-1} \{\alpha_{l,0},\dots,\alpha_{l,{(\frac{n+1}{m}-1})}\}$$
		    is preserved by the Galois group $\Gal(\tilde{K}/\mathbb{Q})$ while acting on $\{\alpha_i\}_{i=0}^{n}$.
\end{enumerate}
Here, {s permutation group $G \subseteq S_{n+1}$ preserves} the partition when any element of the group $G$ shuffles rows of the partition i.e. for every $g\in G$ and $ 0\leq l \leq m-1$ there exists $0\leq k\leq m-1$ such that $g\{\alpha_{l,0},\dots,\alpha_{l,({\frac{n+1}{m}-1}})\}=\{\alpha_{k,0},\dots,\alpha_{k,{\frac{n+1}{m}-1}}\}$.
\end{lemma}
 
\begin{proof}
By Galois theory, the group $G=\Gal(\tilde{K}/\mathbb{Q})$ acts transitively on $\{\alpha_{i}\}_{i=0}^{n}$ since $f$ is irreducible. Suppose $K=\Q(\alpha_0)$.
Let $G_{f}$ be the (not necessarily normal) subgroup of automorphisms in $G$ that fix $K$ pointwise, which is also {the stabiliser in $G$ of $\alpha_0$}. 
We also know from Galois theory that ${[G:G_{f}]=n+1}$ and that the subfields $\mathbb{Q} \subsetneq L \subsetneq K$ are in bijection with intermediate subgroups $G_{f} \subsetneq H \subsetneq G$ via sending $H$ to the fixed subfield of $H$ within $K$ namely $L \mapsto H=\Gal(\tilde{K},L)$
and $H\mapsto L=\lbrace x \in K \; : \; hx=x \; , \forall h\in H \rbrace$.

Let $L \subsetneq K$ be a non-trivial subfield of degree $m\mid n+1$ and $H \supsetneq G_{f}$ be the corresponding subgroup of $G$ whose fixed field in $K$ is $L$. Then, the orbits of $\{\alpha_{0}\}$ 
under the action of $H$ give a subset $P$ of ${\{\alpha_0,\cdots,\alpha_{n}\}}$ with $\#P=[H:G_f]$. Due to $H\neq G_f,G$, we know that $P$ is a nontrivial subset of roots. We claim that $\{gP, g\in G/H \}$ is a nontrivial partition of the roots, where $G/H$ is the space of the left cosets of $H$ in $G$. From the transitivity of the action of $G$ on the roots, we know that $\cup_{g\in G/H}gP=\{\alpha_0,\cdots,\alpha_{n} \}$. But from $n+1=[G:G_f]=[G:H][H:G_f]=\#(G/H)\#P$, we obtain that $\{gP, g\in G/H \}$ is indeed a partition. The $G$-invariance of partition follows from the definition.  

On the other hand, if there is a partition of $\{\alpha_{i}\}_{i=0}^{n}$ 
that is invariant under the action $G$, one can construct $H$ as the set of elements of $G$ that permute the roots of the cell of partition containing $\alpha_0$.
This $H$ must be a strict subgroup of $G$ since $G$ acts transitively on all the roots. Whereas, we have $G_{f} \subsetneq H$ because the cardinality of a cell is $(n+1)/m$, and $H$ acts transitively on the cell containing $\alpha_0$.
\end{proof}

We will now state a sufficient condition to ensure that the field generated by a root of the split minus one $f(X) = \prod_{i=0}^{n}(X-a_{i}) - 1$ 
contains no non-trivial subfield for most parameters $a_{0},\dots,a_{n} \in \mathbb{Z}$.

\begin{lemma}
  \label{le:suff_no_subfield}
  There exists a large enough constant $C>(\sqrt{n+1})^{n+5}$ such that {for every family} of integers 
  $ a_{0}=0 \leq a_1 \leq  \dots \leq a_{n} \in \mathbb{Z}$ satisfying the following conditions:
    \begin{enumerate}
      \item for each $i \in \{0,1,\dots,n\}$, 
  \begin{equation}
	  \label{eq:inequality_hypothesis}
	  \prod_{\substack{j=0 \\ j \neq i}}^n (|a_{i} - a_{j}|-1)  >  C (\sqrt{a_n})^{n-1};
  \end{equation}
  \item for every {non-trivial} partition $\{a_{i}\}_{i=0}^{n} = \bigsqcup_{l=0}^{m-1} \{a_{l,0}, \ldots , a_{l,d-1}\}$ with $n+1=md$, {there exist distinct} $l_{1},l_{2} \in \{0,\dots,m-1\}$ such that
	  \begin{equation}
	      \sum_{j=0}^{d-1} a_{l_{1},j} \neq
	      \sum_{j=0}^{d-1} a_{l_{2},j}.
	      \label{eq:unequal}
	  \end{equation}
  \end{enumerate}
 Then the field $K=\mathbb{Q}[X]/ \langle f(X) \rangle$ has no non-trivial subfields, where we denote by $f(X)=\prod_{i=0}^n(X-a_i) -1$.
\end{lemma}

  \begin{proof}
  We argue {by contradiction}.
  Let us assume \cref{eq:inequality_hypothesis}, \cref{eq:unequal} and that $K$ admits a non-trivial subfield.
  
  By \cref{le:partition_galois} the Galois group $G = \Gal(\tilde{K}/\mathbb{Q})$ of $\tilde{K}=\mathbb{Q}(\alpha_{0},\dots,\alpha_{n})$ preserves a {non-trivial} partition of the roots of $f(X)$ 
  {i.e. there exist} $d,m\geq 2$ such that $dm=n+1$ and
    $\{\alpha_{i}\}_{i=0}^{n} = \bigsqcup_{l=0}^{m-1} \{\alpha_{l,0},\dots,\alpha_{l,{d-1}}\}$. 

Let us now show using \cref{eq:inequality_hypothesis} that
for all $i \in \{0,\dots,m-1\}$,
	  \begin{equation}
		 \sum_{j=0}^{d-1} \alpha_{ij}  
		 = 
		 \sum_{j=0}^{d-1} a_{ij} \in \mathbb{Z}.
		 \label{eq:they_are_int}
	  \end{equation}
It suffices to prove that for all $k\in \lbrace 0,1,\dots,m\rbrace$,  
  \begin{equation}
      \sum_{i=0}^{m-1} \left(\sum_{j=0}^{d-1} \alpha_{ij}\right)^{k} =
      \sum_{i=0}^{m-1} \left(\sum_{j=0}^{d-1} a_{ij}\right)^{k} .
      \label{eq:to_show}
  \end{equation}
  This holds for $k=0$. 
  The case of $k=1$ can be verified by relating the first non-leading coefficient of $f(X)$ and $f(X)+1$. 

  Let us prove \cref{eq:to_show} for every $m\geq k\geq 2$.
  The quantity on the left is an algebraic integer invariant under $G$ and hence is a rational integer in $\mathbb{Z}$.
  Since both sides are integers, it is sufficient to show that the difference is less than one.
First note that since $d \geq 2$, we have 
that $\sup \Big\lbrace \sum_{j=0}^{d-1} a_{ij}, \big\vert\sum_{j=0}^{d-1} \alpha_{ij} \big\vert \Big\rbrace \leq  \sum_{j=0}^{d-1} (a_{ij} +1) \leq 2\sum_{j=0}^{d-1} a_{ij}$.
Then we also use that $m, d \leq (n+1)/2$.
  
  \begin{align}
      \sum_{i=0}^{m-1} 
        \Bigg\vert  &\left(\sum_{j=0}^{d-1} \alpha_{ij}
                \right)^{k}
         - \left(\sum_{j=0}^{d-1} a_{ij}\right)^{k} 
        \Bigg\vert   \\
        &\leq   
        \sum_{i=0}^{m-1} k2^{k-1}
        \left[  
            \left(\sum_{j=0}^{d-1} a_{ij}
            \right)^{k-1}
            \sum_{j=0}^{d-1} |\alpha_{ij} - a_{ij}| 
        \right]\\
        &\leq  k 2^{k-1} (d a_n)^{k-1}
          \sum_{i=0}^{m-1} \sum_{j=0}^{d-1}  |\alpha_{ij} - a_{ij}|
         \\
        &\leq k 2^{k-1} (d a_n)^{k-1} (n+1)\sup_{0 \leq i \leq n} 
            \Big\lbrace|\alpha_{i} - a_{i}|\Big\rbrace \\
        & \leq \frac{ \big(\sqrt{n+1}\big)^{n+5} }{2} \big(\sqrt{a_n}\big)^{n-1} \sup_{0 \leq i \leq n} 
            \Big\lbrace|\alpha_{i} - a_{i}|\Big\rbrace \leq \frac{ \big(\sqrt{n+1}\big)^{n+5} }{2} C^{-1}.
  \end{align}
    Here in the last step, we use the error term in \cref{eq:inequality_of_roots} of \cref{lem:irred} and assumption \cref{eq:inequality_hypothesis} to remove the $\big(\sqrt{a_n}\big)^{n-1}$ term and get the upper bound.
  By hypothesis $C^{-1}\big(\sqrt{n+1}\big)^{n+5}<1$ hence \cref{eq:to_show} is true for all integers $k \in \{0,1,\cdots, m\}$. 
	  
We finish the proof since the transitive action of $G$ will map 
$\sum_{j=0}^{d-1} \alpha_{i_{1}j}$ to $\sum_{j=0}^{d-1} \alpha_{i_{2}j}$
for any $i_{1} \neq i_{2}$. 
On the other hand, each $\sum_{j=0}^{d-1} \alpha_{ij}$ is an integer 
by~\cref{eq:they_are_int} so it is fixed under the action of $G$. The only way this can happen is if for each $i_1 \neq i_2$ one gets
\begin{equation}
\sum_{j=0}^{d-1} a_{i_{1}j} = \sum_{j=0}^{d-1} a_{i_{2}j},
\end{equation}
which contradicts \cref{eq:unequal}.
\end{proof}

\subsection{Using the Cusick-Silverman inequality}

Let $n\geq 3$ be an integer, fix some integers $a_0=0<a_1 <\ldots < a_{n} \in \ZZ_+$, consider the polynomial $f(X)=(X-a_0)\cdots(X-a_n) -1 $, the number field $K=\Q[X]/\langle f(X)\rangle$, ordered embeddings $\sigma_0, \ldots , \sigma_n:K \rightarrow \R $ and the order $\calO=\ZZ[X]/\langle f(X)\rangle$.
Denote by $b_i=\alpha_0 - a_i$ for all $0 \leq i \leq n$ where $\alpha_0\in \overline{\Q}$ is the smallest root of $f$.

Let $\beta_1,\ldots, \beta_n$ be a generating set of $\calO^\times$. 
Denote by $y_{i,j} = \log \vert a_i - a_j \vert$ for all $ i\neq j \in \lbrace 0, \ldots , n \rbrace$.

\begin{proposition}\label{prop:separation-condition}
    Let $n\geq 3$ be an integer, let $C$ be a large enough constant such that $C>n+1$ and ${\log \gamma_{n} + n(n+1)\log C \geq 0}$.
    
    Consider $a_0=0<a_1 <\ldots < a_{n}$ such that 
    \begin{enumerate}
        \item[1.)] $ a_{j+1}-a_j - 1 > C \sqrt{a_n} $  for all $j=0, \ldots, (n-1)$ ;
        \item[2.)] for every {non-trivial} partition $\{a_{i}\}_{i=0}^{n} = \bigsqcup_{l=0}^{m-1} \{a_{l,0}, \ldots , a_{l,d-1}\}$ with $n+1=md$, {there exist distinct} $l_{1},l_{2} \in \{0,\dots,m-1\}$ such that
	  $$    \sum_{j=0}^{d-1} a_{l_{1},j} \neq
	      \sum_{j=0}^{d-1} a_{l_{2},j} .$$
    \end{enumerate}
    Let $f$ and $\calO$ be as above.
    Then the following statements {hold}.
    \begin{enumerate}
        \item[(a)] The shape of $\calO^\times$ lies in a compact set of ${\PGL_n(\R)/\PGL_n(\Z)}$.
        \item[(b)] There exists a constant $C_n>0$ and a finite set $\widetilde{\mathcal{M}}\subset \mathrm{M}_n(\Z)$ such that for every $\lbrace b_1, \ldots, b_n\rbrace \subset \lbrace \alpha_0 - a_0, \ldots, \alpha_0 - a_n \rbrace$ of rank $n$, there exists a shape change matrix $M\in \widetilde{\mathcal{M}}$ and a generating set of $\calO^\times$ such that for every $1 \leq i \leq n$,
    $$ b'_i = \pm \prod_{j=1}^n \beta_j^{M_{ij}}.$$
    \end{enumerate}

\end{proposition}

\begin{proof}

    By \cref{lem:irred} the polynomial $f$ is irreducible and by \cref{le:suff_no_subfield} the field $K=\Q[X]/\langle f(X)\rangle$ has no non-trivial subfields.
    We first use Cusick's inequality of \cref{prop:cusick} to deduce an upper bound of the index of $\Theta = \langle b_1,\ldots , b_n\rangle$ in ${\calO^\times}$.
    Then we prove that the shape of $\Theta$ is in a compact set.
    
    Since $K$ has no non-trivial subfield, by Cusick's inequality,
    $$\reg (\calO)  \geq {c_{n}}\big( {\log\gamma_{n}} + \log  \disc (\calO) \big)^n.$$
    
    On one hand, $\calO=\mathbb{Z}[\alpha]$ hence $\disc (\calO) = \prod_{i>j} (\alpha_i - \alpha_j)^2.$
    By going to the logarithm and using the approximations in \cref{lem:irred} and the hypothesis on the $a_i$ we deduce that
    $$\log \disc (\calO) = 2\sum_{i>j} y_{i,j} + O_n(\varepsilon) \geq {n(n+1)\log \big(C\sqrt{a_n}\big)} $$
    where
$$
\varepsilon=\sup_i\vert a_i-\alpha_i\vert,
\qquad y_{i,j}=\log|a_i-a_j|.
$$
Assuming $C$ is large enough, we have
$$
\log\gamma_n+n(n+1)\log C\geq0.
$$
Hence $\reg \calO \geq {c_{n} \Big( n(n+1)\log \sqrt{a_n} \Big)^n.}$
    
    On the other hand, $\Theta$ is generated by $b_1,\ldots, b_n$. 
    The image of each $b_j$ in $\R^{n+1}_0$ is given by
    $$ v_j = (\log \vert \sigma_0 (b_j) \vert, \ldots, \log \vert \sigma_n (b_j) \vert).$$
    Now by extracting directly the minor 
    $$\reg \Theta = \vert v_1 \wedge \ldots \wedge v_n \vert =
    \big\vert\det
    \begin{pmatrix}
        \log \vert\sigma_1(b_1)\vert & \cdots & \log \vert\sigma_n(b_1)\vert \\
        \vdots &  & \vdots \\
        \log \vert\sigma_1(b_n)\vert & \cdots & \log \vert\sigma_n(b_n)\vert
    \end{pmatrix}
    \big\vert .
    $$
   
    We first prove that the above determinant is nonzero and $\Theta$ has rank $n$. Here we fix the choice $b_i=\alpha_0-a_i$, due to $(\alpha_0-a_0)\cdots(\alpha_0-a_n)=1$, other choices also have rank $n$. Since each $b_i$ is a unit, we have 
    $$\sum_{j=1}^n \log|\sigma_j(b_i)|=\log|\prod_{j=1}^n\sigma_j(b_i)|=-\log|\sigma_0(b_i)|=-\log|\alpha_0-a_i|<0 $$
    The diagonal terms are negative and the off-diagonal terms are positive; since the row sums are negative, the matrix is strictly diagonally dominant, hence its determinant is nonzero. 

    As above, $\log \vert \sigma_j (b_i)  \vert = y_{i,j} + O(\varepsilon)$ for all $i\neq j$.
    Now, $\sigma_j (b_0\ldots b_n)=1$, so
    $\log(\vert \sigma_j (b_j)  \vert)=-\sum_{\overset{0\leq i \leq n}{i\neq j}} y_{i,j}+ O(\varepsilon) .$ To lighten notation, we denote by $y_{j,j}=-\sum_{\overset{0\leq i \leq n}{i\neq j}}y_{i,j}.$
    Hence
    $$\reg \Theta = 
    \big\vert \det
    \begin{pmatrix}
         y_{1,1} & \cdots & y_{1,n} \\
        \vdots &  & \vdots \\
        y_{n,1} & \cdots & y_{n,n}
 \end{pmatrix}\big\vert 
 + O_n(\varepsilon(\log a_n)^{n-1}) \leq n! \sup_{1\leq i,j\leq n} \vert y_{i,j}\vert^n + O_n(\varepsilon(\log a_n)^{n-1}).$$

We deduce the upper bound 
$$D:=[\calO^\times : \{\pm 1\}\Theta]= \frac{\reg (\Theta)}{\reg(\calO)} \leq \frac{ n! \sup_{1\leq i,j\leq n} \vert y_{i,j}\vert^n +O_n(\varepsilon(\log a_n)^{n-1})}{{c_{n} \Big(n(n+1)\log \sqrt{a_n}\Big)^n}}  . $$
Now using the hypothesis on the $a_i$, 
$$\sup_{1\leq i,j\leq n} \vert y_{i,j} \vert \leq n \log a_n.$$
Hence 
\begin{equation}
    \label{eq:bound_on_index}
        D \ll {\frac{n! 2^n}{c_{n}\,(n+1)^n}}
\end{equation}

We show that the shape of $\Theta$ lies in a compact set.
By {Mahler's compactness criterion}, it suffices to find an upper bound for the height of shape i.e. find a lower bound for the shortest vector of the $n$-lattice.
Note that the length of each generator in $\R^{n+1}_0$ is $O_n(\log a_n)$.
Let $v\in \oplus_{1 \leq j \leq n} \Z v_j$ be a non zero element.
By property of the wedge product, for every $1 \leq k \leq n$, 
$$ \Vert v\Vert \geq \sup_{1 \leq k \leq n} \frac{\Vert v \wedge_{j\neq k} v_j \Vert}{\Vert \wedge_{j\neq k} v_j \Vert}.$$
Since $v$ is non zero, we can find $1 \leq k \leq n$ such that $v$ and $\lbrace v_j \; , \; j\neq k \rbrace$ is a basis of $\R_0^{n+1}$. 
Then $\Vert v \wedge_{j\neq k} v_j \Vert \geq \reg \Theta\geq \reg\calO \gg  (\log a_n)^n$, by Cusick inequality.
Since $\Vert \wedge_{j\neq k} v_j \Vert =O\big((\log a_n)^{n-1}\big)$, we deduce that 
$$\Vert v\Vert \gg \log a_n .$$

To finish the proof, denote by $\beta_1,\cdots,\beta_n $ a generating set of $\calO^\times$.
Then there exists an integer matrix $M$ with $|\det M|=D$ such that for every $0 \leq i \leq n$
\begin{equation}
 \begin{pmatrix}
	\log |\sigma_i b_1| \\ \vdots \\ \log | \sigma_i b_n|
\end{pmatrix}=M\begin{pmatrix}
	\log |\sigma_i \beta_1| \\ \vdots \\ \log | \sigma_i \beta_n|
\end{pmatrix}.   
\label{eq:defi_of_M}
\end{equation}
Note that $M$ is defined up to right multiplication by $\PGL_n(\Z)$ since we are allowed to readjust our generating set $\beta_1,\cdots,\beta_n$.
Then one observes that $M \in \calM$ where 
\begin{equation}
\mathcal{M}=
\{ A \in \M_n(\mathbb{Z}) \mid 0<|\det(A)| \leq C_n \}/\sim,
\label{eq:matrix_coset}
\end{equation}
where $C_n$ is chosen large enough that $D\leq C_n$ by \cref{eq:bound_on_index} and $\sim$ is equivalence under right-action of $\PGL_n(\mathbb{Z})$. One needs to check that $\calM$ is a finite set.

Indeed, up to right multiplication by elements of $\PGL_n(\mathbb{Z})$, we can reduce each matrix in $\calM$ to lower triangular matrices with increasing diagonal entries and where {in each row}, the non-diagonal entries are strictly smaller in absolute values than the diagonal entry. 
Hence $\card \calM$ is a finite number and depends only on $n$.
\end{proof}
\begin{remark}
    If $K$ admits non-trivial subfields, our method is not sufficient to conclude $[\calO^{\times} : \langle b_1, \cdots, b_{n} \rangle ] < \infty$, namely because Cusick's inequality has a smaller exponent.
\end{remark}

Inspired by the previous proposition, we state the following definition.

\begin{defin}\label{def-canonical-M-basis}
Let $a_0=0 < a_1 < \cdots < a_n$ that satisfy the conditions of the above proposition.

Note that the choice of $M$ and the generating set of $\calO^\times$ depend on the choice of subset of cardinality $n$ of $\lbrace X-a_0, \ldots , X-a_n \rbrace$.
So for the rest of the paper, we will fix a \emph{canonical} subset of cardinality  $n$, that we denote for every $1 \leq i \leq n$ by 
$$b_i = X-a_i.$$
Then $M \in \widetilde{\mathcal{M}}$ and $\beta_1, \ldots, \beta_n$ will be defined canonically.  
\end{defin}

\subsection{Selmer polynomials and inert primes}
In this section, we setup some auxilary results that will help us in the shapes machine of \cref{se:machine}.

We first need the following lemma to show that there exists an endless supply of desirable primes.
Let us fix the notation
\begin{equation}
    f^{\sel}(X) = X^{n+1} + X^{n} - 1.
\end{equation}
This polynomial is a variant of the Selmer polynomial $X^{n+1} - X -1$. The Selmer polynomial is known to be irreducible and therefore $f^{\sel}$ is also irreducible.

\begin{lemma}\label{lem:irr-selmer}
For each $n\geq 2$, there exist infinitely many primes $p$ such that the polynomial $f^{\sel}(X)=X^{n+1}+X^{n}-1$ is irreducible over $\F_p$.
\end{lemma}
\begin{proof}
    Let $g(X)=X^{n+1}-X-1$ be the Selmer polynomial. Notice that $-X^{n+1}f^{\sel}(1/X)=g(X)$. Hence the Galois group of $f^{\sel}$ is isomorphic to the Galois group of $g$. It is classically known that the Selmer polynomial $g$ has Galois group isomorphic to the symmetric group $S_{n+1} $\cite{nartEquationsTypeXn+aX+b1979, osadaGaloisGroupsPolynomials1987}.

    Hence the Galois group of $f^{\sel}$ is also $S_{n+1}$ and contains $(n+1)$-cycles.
    If $p$ is not a factor of $\disc(f^{\sel})$, then the polynomial $f^{\sel}(X)$ is irreducible modulo $p$ if and only if the  Frobenius morphism $\mathrm{Fr}_p: a \mapsto a^p$ on $\overline{\mathbb{F}}_p$ lifts to an $(n+1)$-cycle in $S_{n+1}$. 
    {By the Chebotarev density theorem}, there exist infinitely many primes $p$ such that the Frobenius $\mathrm{Fr}_p$ is an $(n+1)$-cycle and the density of such primes is $1/(n+1)$ among all primes.    
\end{proof}
\label{ss:ipc}

Such primes $p \in \mathbb{Z}$ are called inert primes of the polynomial $X^{n+1}+X^n-1$. Because of the last lemma, we know that there are infinitely many such primes.

\subsection{About roots of unity}
\label{sec:roots-of-unity}

Let $\calO^{\sel}$ be the monogenic order $\mathbb{Z}[\alpha] \subseteq \mathbb{Q}[\alpha]$ where $\alpha \in \overline{\mathbb{Q}}$ is a solution of $f^{\sel}(X)=0$. 
We will try to analyse to roots of unity of these order $\mathcal{O}^{\sel}$.
In \cref{se:machine}, these roots of unity appear as an {artifact} of the difficulties posed {by} the non-tightness of the Cusick inequality of \cref{ss:cusick}. See \cref{re:cusick-tight}.

We take an inert prime $p$ for $f^{\sel}(X)$ which is guaranteed to exist because of \cref{lem:irr-selmer}.
We get that $\calO^{\sel}/p\calO^{\sel}$ is a field. 
We can define the $p$-adic norm on $\calO^{\sel} \setminus \{0\}$ by $$|x|_p=p^{-\max\{n\in\N,\ x\in p^n\calO \}},$$
and $|0|_{p}=0$. We {may then} look at the metric completion of $\mathcal{O}^{\sel}$ under $|\cdot|_{p}$ and one can convince {oneself} that the metric completion is $\mathcal{O}^{\sel} \otimes _{\mathbb{Z}} \mathbb{Z}_{p}$.

By Hensel's lemma, we know a bit about the roots of unities of $\calO\otimes \Z_p$.
\begin{lemma}
  \label{le:roots_of_unity}
  For an inert prime $p$ and for the Selmer polynomial of degree $n+1$, one has that the torsional units of $\calO^{\sel} \otimes \mathbb{Z}_p$ form a cyclic subgroup of size $p^{n+1}-1$.
\end{lemma}
\begin{proof}
Let $K^{\sel}_{p}$ be the field of fractions of $\mathcal{O}^{\sel} \otimes \mathbb{Z}_{p}$, (equivalently, $K^{\sel}_{p} 
= \mathbb{Q}_{p}[X]/\langle f^{\sel}(X) \rangle$.) 
We know that $K^{\sel}_p$ is unramified extension of $\mathbb{Q}_p$ since $K_p^{\sel} = \mathbb{Q}_p[X] / \langle f^{\sel}(X) \rangle$.

Let $\calO_{K_p^{\sel}}$ be the maximal order in $K^{\sel}_{p}$.
Then we have that 
\begin{equation}
    \mathcal{O}^{\sel}\otimes \mathbb{Z}_{p} \subseteq 
    \calO_{K_p^{\sel}}
    \subseteq K_{p}^{\sel}.
\end{equation}

We know by classical results of Hensel that the {group of torsion elements} of $(\calO_{K_p^{\sel}})^{\times}$ {is} isomorphic to $\mathbb{F}_{p^{n+1}}^{\times}$ since $K_p^{\sel}$ is unramified.

We further claim that in fact $\mathcal{O}^{\sel}\otimes \mathbb{Z}_{p}  =\calO_{K_p^{\sel}}$. We know that $\calO^{\sel}/p \calO^{\sel} \simeq \mathbb{F}_{p^{n+1}}$ because of \cref{lem:irr-selmer}. Hence $\calO^{\sel} \otimes \mathbb{Z}_p$ has a unique maximal ideal $p (\calO^{\sel} \otimes \mathbb{Z}_p)$ which is principal. Hence, $\calO^{\sel} \otimes \mathbb{Z}_p$ is a DVR (discrete valuation ring) and is therefore integrally closed. So it must be equal to the maximal order $\calO_{K_p^{\sel}}$.
\end{proof}

Because of
\cref{le:roots_of_unity}, we can {fix} a $(p^{n+1}-1)$-th root of unity 
$\zeta_{\sel} \in \calO^{\sel}  \otimes \Z_p $. Then for any $k \in \N$ due to $$(\calO^{\sel}  \otimes \Z_p)/p^k(\calO ^{\sel} \otimes\Z_p)\simeq \calO^{\sel} /p^k\calO ^{\sel}  ,$$ 
the element $\zeta_{\sel}$ descends to {a $(p^{n+1}-1)$-th root of unity} in the quotient ring $\calO^{\sel} /p^k\calO^{\sel}$.
This $\zeta_{\sel}$ will be a part of our input congruence parameters in \cref{ss:input}.

Here is an important technical lemma which deals with roots of unity in $\calO \otimes \Z_p$ after projection modulo $\pm p^k \calO^{\sel} = p^{k} \mathcal{O}^{\sel}$.

\begin{lemma}
 \label{le:hensel}
  Let $\mathcal{O}^{\sel}$ be as before and let $p$ be an inert prime. Then the following is an exact sequence for every $k \geq 1$.
  \begin{equation}
  1 \rightarrow  \frac{1+ p \mathcal{O}^{\sel}}{ 1 + p^{k} \mathcal{O}^{\sel}} \rightarrow  \left(\frac{\mathcal{O}^{\sel}}{p^{k} \mathcal{O}^{\sel} } \right)^{\times} \rightarrow  
  \left( \frac{\mathcal{O}^{\sel}}{p \mathcal{O}^{\sel}}\right)^{\times} \rightarrow  1
  \end{equation}
Hence, for an integer $N$ coprime to $p$, one has that the solutions of $X^{N}-1=0$  in $\mathcal{O}^{\sel}/p^{k} \mathcal{O}^{\sel}$ are lifts of roots of unities of $(\mathcal{O}^{\sel}/p \mathcal{O}^{\sel})^{\times}$.
\end{lemma}

\begin{proof}
    Since every unit in $\mathcal{O}^{\sel}/p \mathcal{O}^{\sel}$ lifts to a unit in $\mathcal{O}^{\sel} \otimes \mathbb{Z}_{p}$,
    we know that the reduction map
    $$
    (\mathcal{O}^{\sel}/p^{k}\mathcal{O}^{\sel})^{\times}
    \longrightarrow(\mathcal{O}^{\sel}/p\mathcal{O}^{\sel})^{\times}
    $$
    is surjective. Its kernel is clearly the preimage of
    $1+p\mathcal{O}^{\sel}$, therefore the exact sequence holds.

The second statement is a trivial consequence of the first statement.
\end{proof}

\subsection{Exploiting depth}

Let us write a key algebraic and very general lemma about local behaviour of powers of a fixed non-torsional unit in an order.

\begin{lemma}
    \label{le:powers_of_non-torsion_units}
    
Let $\calO \subseteq K$ be an order in a number field $K$.  Let $\beta \in \calO^{\times}$ be a non-torsional unit. 
Let $p$ be an odd prime such that $\calO \otimes \Z_p$ is a DVR in which $p$ is a uniformizer. 
Let $\zeta \in \calO \otimes \Z_p$ be a generator of the cyclic group of torsional units of $(\calO \otimes \Z_p)^{\times}$.

Then, there exists $l_0$ depending on $p,\calO,\beta$ such that for any $k \geq l_0$ 
  $$\{n\in\Z,\ \beta^n\in \langle\zeta\rangle(1+p^k(\calO\otimes\Z_p)) \}= p^{k-l_0} \cdot \Z,  $$  
    where $\langle\zeta\rangle $ is the finite group generated by $\zeta$ in $\calO\otimes\Z_p$.
\end{lemma}
\begin{proof}
Let us temporarily denote $\calO_p = \calO \otimes \Z_p$ for simplicity.

Since $\beta$ is a non-torsional unit, $\beta^j \notin {\langle \zeta \rangle} \subseteq (\calO_p)^{\times}$ for every $j \geq 1$.
The image of $\beta$ under the projection to $\calO/p\calO$ is invertible. 
Note that $(\mathbb{Z}/p\mathbb{Z})^\times \subset \langle \zeta \rangle.$
Hence one can define $$l_0 = \max \lbrace l \; , \;  \beta \in \langle \zeta\rangle(1 + p^l \calO_p) \rbrace < \infty.$$
Observe that $\langle \zeta \rangle (1 + p^{l_0} \calO_p ) \subseteq \calO_p^{\times}$ is a subgroup under multiplication. So we have for all $j\in \Z$
$$\beta^{j} \in \langle \zeta\rangle ( 1 + p \calO_p) \Leftrightarrow \beta^{j} \in \langle \zeta\rangle ( 1 + p^{l_0} \calO_p) .$$

Then, by definition of $l_0$ one has for every $k \leq l_0$,
$$\beta^l \in \langle \zeta\rangle(1 + p^{k} \calO_p) \Leftrightarrow l \in \Z.$$
So this proves the statement for $k = l_0$.

Assume by induction for some $k\geq l_0$ that
$$
\beta^l \in \langle \zeta\rangle ( 1 + p^k \calO_p ) \Leftrightarrow l \in p^{k-l_0} \Z.
$$
We want to show the following statement for $k+1$: 
\begin{equation}
    \label{eq:induct}
    \beta^l \in \langle \zeta\rangle ( 1 + p^{k+1} \calO_p ) \Leftrightarrow l \in p^{k-l_0+1} \Z.
\end{equation}

Let $\beta^n \in \langle \zeta\rangle ( 1 + p^{k+1} \calO_p)$. Then $\beta^n \in \langle \zeta\rangle( 1 + p^{k} \calO_p)$.
So $n=n_0 n_1$ where $n_0 = p^{k-l_0}$. We need to prove that $p\mid n_1$ for the $\Rightarrow$ implication of \cref{eq:induct}.
By induction assumption, $\beta^{n_0}=  \zeta' (1 + p^k y)$ with $y \in \calO_p \setminus p\calO_p$ and $\zeta' \in \langle \zeta\rangle$.
Note that $$\beta^n =  (\beta^{n_0})^{n_1} = \zeta'^{n_1} (1 + p^k y)^{n_1} \in \langle \zeta\rangle( 1 + n_1 p^k y + p^{2k}\calO_p) .$$
Hence $\beta^n \in \langle \zeta\rangle( 1 + p^{k+1}\calO_p )\Leftrightarrow p\mid n_1,$ since ${y \notin p\calO_p}$.
By induction, for every $k\geq l_0$, we have 
 $\lbrace l\in \Z \; , \; \beta^l \in  \langle \zeta\rangle (1 + p^k \calO_p) \rbrace  = p^{k-l_0} \Z.$
\end{proof}

\subsubsection{Going back to Selmer}

In the next Lemma, we now consider the monogenic order of the Selmer polynomial and take an inert prime $p$ for Selmer polynomial.
We give a condition on non-torsional units of different orders (of different number field) for which the conclusion of the previous Lemma is the same.

\begin{lemma}\label{lem:root-unity}
Let $p$ be an odd inert prime for Selmer polynomial, i.e. such that $\calO^\sel \otimes \Z_p$ is an integral domain,
let $b\in (\calO^\sel)^\times$ be a non-torsional unit, let $\zeta_\sel$ be a generator of the cyclic group of torsional units of $(\calO^\sel \otimes \Z_p)^\times$.  
Then there exists an integer $l_0\in \N$ (which depend only on $b,p,\zeta_\sel$) such that the following holds.

For every $k\geq l_0$,
for every order $\calO$ in a totally real field of degree $n+1$ such that 
$$ 
\frac{\calO^{\sel}}{p^k\calO^{\sel}} \simeq 
 \frac{\calO}{p^k\calO},
 $$
 for any torsional unit $\zeta\in(\calO\otimes\Z_p)^\times$ such that $\zeta \equiv \zeta_\sel \pmod{p^k \calO}$,
 for any unit $\beta \in \calO^\times$ such that $\pi_k(\beta)^d\in \langle \zeta \rangle b$ with $ \gcd(p,d)=1$, we have

$$     \left\lbrace n\in \Z , \;  \beta^n \in \langle \zeta\rangle( 1+p^k \calO )\right\rbrace = p^{k-l_0}\Z $$
where $\langle\zeta\rangle$ is the finite group generated by $\zeta$ in $(\calO\otimes\Z_p)^\times$ and $\pi_k:\calO \rightarrow \calO/p^k\calO$ is the canonical projection.

Furthermore, there exists $ 0\leq \epsilon< p^{n+1}-1$ (depending on $\beta$ and independent of $k$) such that for every $k\ge l_0$
$$     {\beta^{p^{k-l_0}} \equiv \zeta^{\epsilon p^k} \pmod{p^k \calO}.}$$
\end{lemma}

\begin{proof}

Let $b$ and $p$ be as in the statement. Then by \cref{le:powers_of_non-torsion_units}, there exists an integer $l_0\in \N$ such that for every $k \geq l_0$ 
  $$\{n\in\Z,\ b^n\in \langle\zeta_\sel\rangle(1+p^k(\calO^\sel \otimes\Z_p)) \}= p^{k-l_0} \cdot \Z.  $$ 
In particular, there exists $\zeta_0\in \langle \zeta_\sel\rangle$ and $y_0 \in \calO^\sel \setminus p\calO^\sel$ such that 
$$b=\zeta_0(1+p^{l_0}y_0).$$

Fix now an integer $k\geq l_0$, an order $\calO$, a torsional unit $\zeta\in(\calO\otimes\Z_p)^\times$ and a non-torsional element $\beta \in \calO^\times$ as in the statement.
Now by \cref{le:powers_of_non-torsion_units} again, consider 
$$l_1 = \max\lbrace l , \; \beta \in \langle \zeta \rangle(1 + p^l(\calO\otimes\Z_p))\rbrace.$$
We shall prove that $l_1=l_0$.
There exists $\zeta_1 \in \langle \zeta\rangle$ and $y_1\in(\calO\otimes\Z_p)^\times$ such that $\beta = \zeta_1(1+p^{l_1} y_1)$.
Then $$\beta^d = \zeta_1^d(1+p^{l_1}y_1)^d \in  \langle\zeta\rangle (1 + d p^{l_1}y_1 + p^{2 l_1} (\calO\otimes\Z_p)).$$

Using that $\pi_k(\beta)^d \in \langle \zeta\rangle b$ and the identification between $\calO/p^k\calO$ and $\calO^\sel/p^k\calO^\sel$,
we deduce that $dp^{l_1}\mid p^{l_0}$ and $p^{l_0}\mid d p^{l_1}$.
Now by assumption $d$ and $p$ are coprime, hence $l_1=l_0$.
   
For the torsional parts, choose a residue $\epsilon$ modulo $p^{n+1}-1$ such that
$$\beta\equiv \zeta^{\epsilon p^{l_0}} \pmod{p^{l_0} \calO} .$$
Then we notice that for every $l_0 \leq l \leq k$
$$\beta^{p^{l-l_0}} \equiv (\zeta^{\epsilon p^{l_0}})^{p^{l-l_0}} = \zeta^{\epsilon p^{l}} \pmod{ p^{l} \calO}. $$
\end{proof}

\section{Shapes machine}
\label{se:machine}

In this section, we explain the working of the ``shapes machine'', by which we mean \cref{prop:structure-suborders}. 

\subsection{Input data}
\label{ss:input}

To choose a field and an order inside the number field, we define input congruence parameters.

These are given by 
\begin{enumerate}
\item An inert prime $p \in \Z$ for $X^{n+1}+X^{n}-1$ 
    \item A tuple $c \in \lbrace \sum_{i=1}^n e_i \rbrace \cup_{j=1}^n \lbrace -e_j \rbrace$.
    \item A generator $\zeta_{\sel}$ of the group of {roots of unity} of $\calO^{\sel}  \otimes\Z_p$.
    \item An integer $N \in \N$.
\end{enumerate}

With any input parameters $(p,N,c,\zeta_{\sel})$, we define a set of tuples $(a_1, \ldots , a_n) \in \Z^n$ such that
$$a_i  \equiv c_i \pmod{p^N} \text{ for each } i =1,\cdots,n.$$
We then say that $(a_1, \cdots , a_{n}) \in \Z^{n}$ satisfies the input congruence parameters $(p,N,c,\zeta)$.

\begin{remark}
Note that if $a\in \Z^n$ satisfies input congruence parameters $(p,N,c,\zeta)$ for an inert prime $p$ for the  Selmer polynomial, then \cref{eq:unequal} of \cref{le:suff_no_subfield} is satisfied (easy to check).
We still require separation of the coefficients $a_0, a_1, \ldots, a_n$ to deduce that $\mathbb{Q}[X]/\langle (X-a_0)\cdots(X-a_n)-1 \rangle$ has no non-trivial subfields.
\end{remark}

\subsection{Output congruence parameters}
\label{ss:output}

We define output congruence parameters to be
a set containing the following 
\begin{enumerate}
\item A prime $p \in \Z$ and an integer $N\in \N$
    \item Tuples $\phi, \psi \in \Z^{n}$, which we consider as linear forms on $\Z^{n}$ via $x \mapsto \langle \phi, x\rangle$ and $x \mapsto \langle \psi, x\rangle$ respectively,
    \item Integers $l, \epsilon \in \N$ with $0\leq\epsilon<p^{n+1}-1$, and put $m_p=(p^{n+1}-1)/(p-1)$.
\end{enumerate}

We say that a lattice $\Lambda \subseteq \Z^{n}$ is defined by the output congruence parameters $(p,N,l, \epsilon,\phi,\psi)$ {{up to depth $k$}} if $\Lambda$ is the set of integers $(z_1,\cdots,z_{n})\in\Z^{n}$ such that 
	\begin{align}
		\psi(z_1,\cdots,z_{n})&=0\mod  p^{k-l}\\ 
		\epsilon p^l \psi (z_1,\cdots,z_{n})+\phi (z_1,\cdots,z_{n})&=0\mod m_p,
		\label{eq:depthk}
	\end{align}

\subsection{Working of the ``shapes machine''}

Choose an input congruence parameter $(p,N,c,\zeta_{\sel})$ 
as in~\cref{ss:input}. 
Let $a \in \Z^{n}$ be a tuple satisfying the input congruence conditions given by $(p,N,c,\zeta_{\sel})$ and consider the ``split minus one'' polynomial $f(X)$ with parameters given by $a$. That is,
\begin{equation}
    f(X) = X(X-a_{1})(X - a_{2}) \dots (X-a_{n}) - 1,
    \label{eq:splitmone}
\end{equation}
where $(a_{1},\dots,a_{n})=a$. 
Let $\calO$ be the monogenic order $\mathbb{Z}[\alpha]$ where $\alpha \in \overline{\mathbb{Q}}$ is a solution of $f(X)=0$.

The main observation to exploit here is that because of the choice of input of \cref{ss:input}, one gets that for $k \leq N$ we have 
\begin{equation}
\frac{\calO^{\sel}}{p^k\calO^{\sel}} \simeq 
 \frac{\calO}{p^k\calO}.
 \label{eq:order_iso}
\end{equation}
Indeed, this is because by our assumption on $c$ in \cref{ss:input}. We get when $c=(1,1,1,\dots)$
\begin{equation}
f(X) \equiv
 X(X-1)^n-1  = f^{\sel}(X-1) \pmod{p^N}.
\end{equation}
and otherwise
\begin{equation}
f(X) \equiv
 (X+1)X^n-1 = f^{\sel}(X)  \pmod{p^N}.
\end{equation}

Because of \cref{eq:order_iso}, one gets that there is a $\zeta \in \calO$ such that
\begin{equation}
    \label{eq:zeta_iso}
    \zeta \equiv \zeta_{\sel} \pmod{p^N \calO}.
\end{equation}

\subsection{Congruence parameters for congruence suborders}

Recall the definition of the finite set $\calM$ from \cref{eq:matrix_coset}. Note that even though $\calM$ is a set of equivalence classes, we will take a system of {representatives} of each equivalence class and choose a fix a finite set of matrices $\widetilde{\calM} \subseteq \M_n(\mathbb{Z})$ with $\card \widetilde{\calM} = \card \calM$

For the proposition below, we set up the following notation.
Let $a=(a_1,\cdots,a_n) \in \mathbb{N}^n$.
Denote by $f(X)=(X-a_0)\cdots(X-a_n)-1$ and by $\calO_a = \Z[X]/\langle f(X) \rangle$.
Consider a root of unity $\zeta$ satisfying \cref{eq:zeta_iso}, the canonical choice of shape change matrix $M = M(a)\in \tilde{\mathcal{M}}$ and generating set $\beta_1, \ldots, \beta_n$ of $\calO_a^\times$ as in \cref{def-canonical-M-basis}.
We denote
the subgroup
\begin{equation}
 \calO_{p^k,a}^\times =\calO_a^{\times} \cap (\Z + p^k \calO_a ).
\end{equation}

\begin{proposition}\label{prop:structure-suborders}
Consider the constant $C_n$ defined in \cref{eq:matrix_coset}.
Consider an input parameters $(p,N,c,\zeta_\sel)$
as defined in \cref{ss:input}. Assume that $p>\max(2,n+1,C_n)$ and $N \geq l=l(p,c,\zeta_\sel)$ where $l(p,c,\zeta_\sel)$ is the $l$ coming from \cref{lem:root-unity} after choosing the $b$ used there as $X$ or $X-1$ depending on $c$
Then there exists a finite set of output parameters $$\Big\lbrace (p,N,l, \epsilon, \psi_M,\phi_M) ; \; M\in \widetilde{\mathcal{M}} \text{ and } {0\le\epsilon<p^{n+1}-1}   \Big\rbrace$$ such that the following holds.

Whenever $a = (a_1, \cdots, a_{n}) \in \N^n$ is such that 
\begin{itemize}
    \item condition 1.) of \cref{prop:separation-condition} is satisfied for a large enough constant $C>0$; i.e $ a_{j+1}-a_j > C \sqrt{a_n}$ for all $0 \leq j \leq n-1$ with $a_0=0$,
    \item the input congruence condition of parameters $(p,N,c,\zeta_\sel)$ is satisfied,
\end{itemize}
then the following conditions are imposed on $\calO_a$.

Take $M=M(a)$ as in the notations above. Then, for all $l \leq k \leq N$, we have
$$\calO_{p^k,a}^\times=\{\pm\beta_1^{x_1}\cdots\beta_{n}^{x_{n}}, \; (x_1,\cdots,x_{n})\in\Lambda_M \} $$
where $\Lambda_M$ is the lattice defined up to depth $k$ by the output congruence condition of parameters $(p,N,l,\epsilon, \psi_M,\phi_M)$. Meaning, by \cref{eq:depthk}, that $\Lambda_M$ is the set of integer solutions of the system
$$ \left\lbrace 
\begin{array}{c l}	 
\psi_M(x) \equiv 0 &\pmod{p^{k-l}}    \\
\epsilon p^l \psi_M(x) +\phi_M(x) \equiv 0 &\pmod{m_p}.
\end{array} \right.$$

	\begin{remark}\label{rem:structure-suborders}
	\begin{itemize}
	\item The final $\psi_M$ would really depend on $M$ and also on $c$. 
	\item The parameters $\epsilon$ and $\phi_M$ are one of finitely choices. For each $M=M(a) \in \widetilde{M}$ one of these finitely many can be selected.
	    \item The canonical shape change matrix and generating set of $\calO_a^{\times}$ exist because condition 2.) of \cref{prop:separation-condition} holds modulo $p^N$ so it has to be true in $\Z$.
	    \item {For $n\ge3$, the Cusick inequality is not as sharp as in the totally real cubic case.}
	    So our strategy to calculate units of suborder requires to choose a root of unity $\zeta_\sel$.

	    \item Later, we will choose $2n+1$ different primes to impose input congruence conditions and parametrize the parameters $c$ and $\zeta_\sel$ on the primes.
	    Note that for the same tuple of $a_i$s, there is a canonical choice of shape change matrix because of the way $\widetilde{\calM}$ is defined and the way $\beta_1,\cdots,\beta_n$ are chosen.
	\end{itemize}
	
	\end{remark}

\end{proposition}

\begin{proof}
	We apply \cref{prop:separation-condition} to $a_0, \ldots, a_n$.
	Its hypothesis 2.) holds modulo $p^N$ by choice of $c \in \cup_{j=1}^n \lbrace -e_j \rbrace \cup \lbrace \sum_{j=1}^n e_j\rbrace$ hence hypothesis 2.) is true in $\Z$. 

    Recall that by canonical choice given in \cref{def-canonical-M-basis},
    for every $1 \leq i \leq n$,
$$b_i = \varepsilon_i \prod_{j=1}^n \beta_j^{m_{i,j}},$$
where $b_i=X-a_i$ for every $1 \leq i \leq n$.

Consider the adjoint matrix of $M$, that we denote by $M^{\adj}=(\tilde{m}_{ij})_{1\leq i,j\leq n}$.
We have $M^{\adj}M=(\det M) I_{n}$.
Denote by $d= \det M$.

    Define the following map for every $c \in \cup_{i=1}^n \lbrace  -e_i\rbrace \cup \lbrace \sum_{i=1}^n e_i \rbrace$,
    $$ v(c)= \left\lbrace \begin{array}{ll}
        c & \text{ if } c=\sum_{i=1}^n e_i \\
        (n+1)c + \sum_{j=1}^n e_j & \text{ otherwise.}
    \end{array} \right. $$

Let $1 \leq k \leq N$.
We work on the quotient finite ring $\calO/p^k\calO$ and up to multiplication by $\pm 1$. 
Let $\pi_k$ be the quotient map to ${\calO/p^k\calO}$.
There are two cases depending on $c$.
\begin{enumerate}
    \item If $c=\sum_{i=1}^n e_i$, then we denote by $b=\pm (X-1)$ and $\tilde{b}=b_1$ so that $\pi_k(\tilde{b})=b$.
    \item Otherwise, there exists $1 \leq i \leq n$ such that $c=-e_i$. Then denote by $b=\pm X$ and $\tilde{b}=b_0$ so that $\pi_k(\tilde{b})=b$.
\end{enumerate}
Now, in either case, by input congruence condition or using that $b_0\cdots b_n = 1$, 
we deduce that for every $1 \leq j \leq n$,  
    $$ \pi_k(b_j)=b^{v_j(c)}. $$

	\textbf{Step 1: Finding common generator of $\pi_k(\beta_j)$:}

    We consider
	\begin{equation}\label{equ:lj-cj}
	   (h_1, \ldots, h_n)^t:= M^{\adj}v(c).
	\end{equation}
	In the rest of the proof, abusing notations, we will write $v_j$ instead of $v_j(c)$ for the coordinates of vector $v(c)$.
	 Then from \cref{eq:defi_of_M}, we have 
	\begin{equation}\label{equ:beta-j}
	\begin{split}
	  b^{h_j}&=\prod_{i=1}^n \big( b^{v_i} \big)^{\tilde{m}_{ji}}
	  =\prod_{i=1}^n \pi_k(b_i)^{\tilde{m}_{ji}} \\
	  &=\prod_{i=1}^n \prod_{l=1}^n \pi_k(\beta_l^{m_{il}})^{\tilde{m}_{ji}}
	  =\prod_{l=1}^n\pi_k \Big(\beta_l^{\sum_{i=1}^n m_{il}\tilde{m}_{ji}} \Big)\\
	 b^{h_j} &= \pi_k(\beta_j)^{d}.
	\end{split}\end{equation}
	Now choose an integer solution $(f_j)$ of
	$$ \sum_{j=1}^n f_j h_j=h=\gcd(h_j). $$
	Consider $ \beta=\prod_{j=1}^n \beta_j^{f_j}$. From \cref{equ:beta-j} we have
	$$
	 \pi_k(\beta)^{d}=\prod_{j=1}^n\pi_k(\beta_j)^{d \; f_j}=\prod_{j=1}^n b^{h_jf_j} =b^{h}.
	 $$
	Since $M(h_j)^t=d v$ {and $\gcd(v(c))=1$}, we know that ${h\mid d}$, so an $h$-th root of unity appears, since
	  \begin{equation}\label{equ:pi-k-beta}
	      \pi_k(\beta^{d/h} \cdot \tilde{b}^{-1} )^{h}=1.
	 \end{equation}
	 Using the assumption that $p > C_n$ and since $ C_n \geq d$ by equation \cref{eq:matrix_coset} we know that $\gcd(d,p) = 1$ and therefore by \cref{le:hensel}, we deduce that 
\begin{equation}\label{equ:pi-k-beta-principal}
	      \pi_k(\beta)^{d/h} \in \langle \zeta \rangle \; b.
	 \end{equation}
	 
Moreover for every $1 \leq j \leq n$
$$ ( \pi_k(\beta)^{d})^{h_j/h}=(b^{h})^{h_j/h}=b^{h_j}= \pi_k(\beta_j)^{d}. $$
Since the image of $\zeta$ generates the torsional units modulo scalars, in $(\calO/p^k\calO)^\times/(\Z/p^k\Z)^\times$ there is a unique $w_j\in \{0,\dots,m_p-1\}$ such that
\begin{equation}\label{equ:betaj}
\pi_k(\beta_j)=\pi_k(\beta)^{h_j/h}\zeta^{w_j}\quad\text{in }(\calO/p^k\calO)^\times/(\Z/p^k\Z)^\times.
\end{equation}
Let us now define the linear forms with integer coefficients
\begin{equation}\label{equ:psi}
\psi_M(x_1,\cdots, x_{n})=\sum_{j=1}^n \frac{h_j}{h}x_j.
\end{equation}
\begin{equation}\label{equ:phi}
\phi_M (x_1, \cdots , x_n) = \sum_{j=1}^n w_j x_j. 
\end{equation}
Note that the parameters $(h_j)$ and $(w_j)$ depend on $M$ and $c$, so does $\psi_M\in \Z^n$.

	\textbf{Step 2: Solving equations}
	
To sum up the previous step, we have
	\begin{align}
		\calO_{p^k, a}^\times&= \calO_{a}^\times \cap ( \Z + p^k \calO ) \\
		&=\{\pm\prod_j\beta_j^{x_j},\ \prod_j\beta_j^{x_j}\in \Z+p^k\calO  \} \\
		&=\{\ \beta^{\psi_M (x)}\zeta^{\phi_M(x)}\in \Z+p^k\calO  \},
	\end{align}
	where the linear forms $\psi_M, \phi_M$ only depend on the matrix $M$ and {$\zeta$ is a fixed generator of the cyclic group of roots of unity in $\calO/p^k\calO$, namely the image of $\zeta_{\sel}$}.
	
	It is equivalent {to solving} the equation 
	\begin{equation}\label{equ:beta-zeta}
		\beta^{\psi_M(x_1,\cdots,x_{n})} \zeta^{\phi_M(x_1,\cdots,x_{n})}\in \Z+p^k\calO.
	\end{equation}
	The strategy is to solve {the non-torsional part} first using \Cref{lem:root-unity} and then cancel out the root of unity.
	
	First, the solution of \cref{equ:beta-zeta}, must satisfy
	$$ 	\beta^{\psi_M(x_1,\cdots,x_{n})}\in{\langle \zeta\rangle(1+p^k\calO)}.  $$
	By \Cref{lem:root-unity} and since $\beta^{d/h} \in \langle \zeta \rangle b $, there exist $l\in \N$, depending only on $p,b,\zeta_\sel$ such that
	$$ 	\beta^{\psi_M(x_1,\cdots,x_{n})}\in{\langle \zeta\rangle(1+p^k\calO)}\Longleftrightarrow p^{\max\{k-l,0 \}}|\psi_M(x_1,\cdots,x_{n}),  $$
	
	Furthermore, there exists ${0\le\epsilon<p^{n+1}-1}$ depending on $M$ such that
   $$ {\beta^{p^{k-l}} \equiv \zeta^{ \epsilon p^k}} \pmod{p^k \calO} .$$
       The torsion factor in \Cref{equ:beta-zeta} must have residue in $\mathbb F_p^\times$, hence the tuple $(x_1,\cdots, x_n)$ is a solution to the system
   $$\left\lbrace \begin{array}{c c l}
    \psi_M(x_1,\cdots,x_n) &\equiv 0 &    \pmod{p^{k-l}}  \\
    \epsilon p^l  \psi_M(x_1,\ldots,x_n)  + \phi_M(x_1,\ldots,x_n)   &\equiv 0 & \pmod{m_p},\qquad m_p=\frac{p^{n+1}-1}{p-1}.
   \end{array}\right.$$
\end{proof}

\section{Density of \texorpdfstring{$S$}{S}-arithmetic diagonal orbits}
\label{se:lal}

 This section contains the construction from congruence equations with $2n+1$ primes to S-arithmetic setting. Then we prove the lattice approximation lemma.

\subsection{Projection of an \texorpdfstring{$S$}{S}-arithmetic diagonal orbit}

We are interested in some special elements $x \in \mathbb{Q}_S^{n-1}$ which we describe below.

\begin{defin}	\label{de:rfr}
	Let $S$ be a finite set of primes. 
	We say that $x \in \mathbb{Q}_S^{n-1}$ is \textbf{rational} and \textbf{strongly full rank}
	if there exists a subset $S' \subseteq S$
	of $2n+1$ primes that admits a partition
	$ S' = S_0 \sqcup S_1 \sqcup S_2$ with $S_0 = \{p_0\}$ and $S_1=\{p_1, \ldots, p_n\}$ and $S_2=\{q_1, \ldots, q_n\}$  such that the following holds:
	\begin{enumerate}
		
		\item $x_{p}\in \mathbb{Q}^{n-1}$, for every $p \in S'$.
		For every $1 \leq i \leq n-1$, we denote by $x^{(i)}_p\in \mathbb{Q}$ the $i$-th coordinate of $x_p$;
		\item For every $1 \leq j \leq n$, we have $x_{p_j}=x_{q_j}$;
		\item 
		Every $n \times n$ minor of
		the matrix $A_x \in \M_{ (S_0\cup S_1) \times n}(\mathbb{Q})$ 
		is non-zero:
		\begin{equation}
			A_x = 
			\begin{pmatrix}
				1 & x_{p_{0}}^{(1)} & x_{p_{0}}^{(2)} & \cdots & x_{p_{0}}^{(n-1)} \\
				1 & x_{p_1}^{(1)} & x_{p_1}^{(2)} & \cdots & x_{p_1}^{(n-1)} \\
				\vdots&\vdots & \vdots & &\vdots &\\
				1 & x_{p_n}^{(1)} & x_{p_n}^{(2)} & \cdots & x_{p_n}^{(n-1)} \\
			\end{pmatrix}.
		\end{equation}

	\end{enumerate}

\end{defin}

\begin{proposition}
	\label{pr:lal_high_rank}
	Let $S$ be a set of primes such that $\card S \geq 2n+1 $,
	let $x = (x^{(i)})_{1 \leq i \leq n-1} \in \mathbb{Q}_S^{n-1}$
	be \textbf{rational} and \textbf{strongly full rank} in the sense of above \cref{de:rfr}.
	Consider a point $\eta(x)\in \Ni(\mathbb{Q}_S) \subseteq \G(\mathbb{Q}_{S})$. 
	Then the orbit
	\begin{equation}
		\Hi(\mathbb{Q}_{S}) \eta( x)  \G(\mathcal{O}_{S})  \subseteq  \G(\mathbb{Q}_{S^{*}}) / \G(\mathcal{O}_{S}) 
		\label{eq:injection_of_orbit}
	\end{equation}
	maps to a dense subset under the projection map 
	\begin{equation}
		\pi: \G(\mathbb{Q}_{S^{*}}) / \G(\mathcal{O}_{S}) \rightarrow  \G(\mathbb{R})/\G(\mathbb{Z}).
	\end{equation}
\end{proposition}

\begin{remark}
	Some important things to remark about the left side of~\cref{eq:injection_of_orbit}. 
	\begin{enumerate}
		
		\item 
		We are implicitly using the inclusion defined in \cref{eq:implicit_inclusion} to consider the action of $\Hi(\mathbb{Q}_{S})$ and $\eta(x)$ on $\G(\mathbb{Q}_{S^{*}})/\G(\mathcal{O}_{S})$ in~\cref{eq:injection_of_orbit}. 
		\item 
		One could replace the action of $\Hi(\mathbb{Q}_{S})$ with just the matrices $a_{S}(k)$ defined in \cref{eq:defi_of_ask} for $k \in \mathbb{Z}^{S}$, since the projection map $\pi: \G(\mathbb{Q}_{S^{*}})/\G(\mathcal{O}_{S})$ quotients on the left with $\G(\mathbb{Z}_{S})$.
	\end{enumerate}
\end{remark}

For our application, it is sufficient to work with the following stronger version of \cref{pr:lal_high_rank}. We will prove this proposition in \cref{ss:banana}.

\begin{proposition}	\label{pr:strong_lal}
	Let $S$ and $x$ be as in~\cref{pr:lal_high_rank}. Consider the multiplicatively closed subset $\mathcal{H}_{\leq 0} \subseteq \Hi(\mathbb{Q}_{S})$ defined as
	\begin{equation}
		\mathcal{H}_{\leq 0} = \{a_{S}(k) \mid k \in (\mathbb{Z}_{\leq 0})^{S}\},
	\end{equation}
	where $a_{S}(k)$ is as in~\cref{eq:defi_of_ask}. Then the image
	\begin{equation}
		\pi \big(\mathcal{H}_{\leq 0} \cdot \eta(x) \cdot \G(\mathbb{Q}_{S})\big) \subseteq     \G(\mathbb{R})/\G(\mathbb{Z}),
	\end{equation}
	is dense in $\G(\mathbb{R})/\G(\mathbb{Z})$.

\end{proposition}

Clearly Proposition~\ref{pr:lal_high_rank} is implied by~\cref{pr:strong_lal}.

\subsection{Equivalent representation of the diagonal orbit}

We note the following lemma about dual lattices.

\begin{lemma}
	\label{le:duality}
	Let $\Lambda = g\cdot \mathbb{Z}^{n}$ be a lattice inside $\mathbb{R}^{n}$ for some $g \in \GL_{n}(\mathbb{R})$. Then, 
	the lattice 
	$$\Lambda^{\sharp} =\{x \in \mathbb{R}^{n}  \mid x^{t} y \in \mathbb{Z} \text{ for all } y \in \Lambda \}$$ 
	is given by 
	\begin{equation}
		\Lambda^{\sharp} = (g^{-1})^{t} \mathbb{Z}^{n}.
	\end{equation}
	\begin{proof}
    It follows from the fact that $x^t g z \in \Z$ for all $z \in \Z^n$ is equivalent to $( g^t x )^t \in \Z^n$, hence the Lemma.
	\end{proof}

\end{lemma}

We will use the following equivalent representation of the $\pi(\Hi(\mathbb{Q}_{S}) \eta(x) \G(\mathcal{O}_{S}))$.
Recall the notation $[x]_{k}$ for $x \in \mathbb{Q}_{S}^{n-1}$ defined in \cref{ss:defi_of_approx} and the notation $a(\cdot)$ from \cref{ss:diagonal}.  

\begin{lemma}
	\label{le:equivalent_representation}
	Consider for some
	$x\in \mathbb{Q}_{S}^{n-1}$
	the matrix $\eta(x)\in \Ni(\mathbb{Q}_{S})$.
	Then
	for $k \in \mathbb{Z}^{S}$ (and, in part~\ref{it:equivalent-lattice}, assume $k_p\geq0$ and $x_p\in\mathbb Z_p^{n-1}$ for all $p\in S$),
	define
	\begin{equation}
		y = 
		\eta(- [ x]_{k}) a\Big(\prod_{p \in S} p^{k_{p} } \Big) .
	\end{equation}
	
	Then 
	\begin{enumerate}
	    \item 
	    	The coset 
	$y\cdot \G(\mathbb{Z}) \in \G(\mathbb{R})/\G(\mathbb{Z})$ satisfies 
	\begin{equation}
		y\cdot \G(\mathbb{Z})= \pi\big( a_{S}(-k) \eta(x) \G(\mathcal{O}_{S}) \big) .
	\end{equation}
	\item 
		The point $y \G(\mathbb{Z}) \in \G(\mathbb{R}) / \G(\mathbb{Z})$ does not depend on the choice of the approximation $[x]_k \in \mathbb{Q}^{n-1}$.
	\item\label{it:equivalent-lattice}
	The lattice $\Lambda = y  \mathbb{Z}^{n} $ is given upto scaling by
	\begin{equation}
		\Lambda =\{  (z_{1},\dots,z_{n})\in \mathbb{Z}^{n} \mid  z_{1} + x^{(1)}_{p} z_{2} + \dots + x^{(n-1)}_{p} z_{n} = 0 \pmod{p^{k_{p}}} \text{ for each } p \in S  \}.
	\end{equation}
	\end{enumerate}

\end{lemma}
\begin{proof}
\begin{enumerate}
    \item 	We will now use the $\iota$ defined in \cref{eq:implicit_inclusion} explicitly.
	
	For each $p \in S$ and $i \in \{1,\dots,n-1\}$, one has that $x_{p}^{(i)} - [x]^{(i)}$ must be in $p^{k_p} \mathbb{Z}_{p}$. 
	One has for each $p \in S$
	the elements of $\G(\mathbb{Q}_{p})$ satisfy
	\begin{align}
		& \Big(\iota(~a_{S}(-k)\;\eta(x)~) \cdot \eta( - [x]_{k} )   \; a\left(\prod_{q \in S}^{} q^{k_{q}} \right) \Big)_{p}=   \\
		& 
		\begin{pmatrix}
			\prod_{q \in S \setminus \{p\}}^{} q^{k_{q}}  &  p^{-k_{p}} (x^{(1)} - [x^{(1)}]_{k} ) & \dots & p^{-k_{p}}(x^{(n-1)} - [x^{(n-1)}]_{k} )\\
			& &  & \\ 
			& 1 & & \\
			& &  & \\ 
			& & \ddots & \\ 
			& &  & \\ 
			& & & 1 
		\end{pmatrix} .
	\end{align}
	Since, $\prod_{q \in S \setminus \{p\}} q^{k_{q}} \in\mathbb{Z}_{p}^{\times}$ for every $p \in S$, the matrix therefore lies in $\G(\mathbb{Z}_{p})$.
	
	We get by our construction that
	$$\pi ( \iota (a_{S}(-k) \eta(x)) \G(\mathcal{O}_{S})) = 
	\pi (\iota(a_{S}(-k) \eta(x)) \eta(-[x]_{k}) a\left(\prod_{p \in S}{p^{k_{p}}} \right) \G(\mathcal{O}_{S})).$$
	Hence we obtain from the place at infinity that the image is
	\begin{equation}
		y \cdot \G(\mathbb{Z})=   \eta(- [x]_{k}) a\left(\prod_{p \in S}{p^{k_{p}}}\right)  \G(\mathbb{Z}).
	\end{equation}
	
	\item 
  Follows from the proof of the previous section. We don't need to specify which approximation $[x]_k$ we use.

	\item 
    	For this part, we use \cref{le:duality} to conclude that if $\Lambda = y \cdot \mathbb{Z}^{n}$ then 
	\begin{equation}
		\Lambda^{\sharp}  =  \eta([x]_{k})^{t} a\left(\prod_{p \in S}p^{-k_{p}}\right) \mathbb{Z}^{n}.
	\end{equation} 
	Since $(\Lambda^{\sharp})^{\sharp}  = \Lambda$, one knows that
	\begin{equation}
		\Lambda = \{ z \in \mathbb{R}^{n} \mid z^{T} w \in \mathbb{Z} \text{ for all } w \in \Lambda^{\sharp}\}.
	\end{equation}
	This implies that $z \in \Lambda$ if and only if 
	\begin{align}
		z_{2},z_{3},\dots,z_{n} \in \mathbb{Z}\text{ and } 
		z_{1} + z_{2} [x^{(1)} ]_{k} + \cdots + z_{n} [x^{(n-1)}]_{k} \in \Bigl(\prod_{p\in S}p^{k_{p}}\Bigr) \mathbb{Z}.
	\end{align}

\end{enumerate}

\end{proof}

\subsection{Banana lemma}
\label{ss:banana}
The main tool to show \cref{pr:lal_high_rank} is a Theorem due to Dani--Raghavan \cite[Proposition 1.5]{DaniRaghavan1980}.
We use a generalised version of their result in \cite[Theorem 1.4]{shah1996limit}.
We adapt it to our choice of Lie group, our choice of lattice, our choice of horospherical group and omit the part of about measures. 
 \begin{theorem}{\bf{ Banana Lemma } \cite[Theorem 1.4]{shah1996limit}} \label{le:banana}
    Let $\lbrace a(i)\rbrace_{i\in \N}$ be a sequence in (the closed positive Weyl chamber) 
    $\lbrace \diag(b_1, \ldots , b_n) \; : \;  b_1  \geq \ldots \geq  b_n >0 \rbrace,$ 
    such that 
    $$ \left\lbrace
    \begin{array}{ll} 
    \log a_1(i)- \log a_{2}(i) \longrightarrow + \infty &\text{ as } i \rightarrow + \infty \\
    \sup_{i\in \N} \log a_j(i) - \log a_{j+1}(i) < + \infty &\text{ for every } 2 \leq j \leq n-1 .
    \end{array}\right.
    $$
    Denote by $U^+=\lbrace g \in \G(\R) \; : \; a(i)^{-1} g \; a(i) \rightarrow I_n \text{ as } i \rightarrow +\infty \rbrace$ and by $\pi : \G(\R) \rightarrow \G(\R) / \G (\Z)$ the natural projection.
    Let $\rho : \G(\R) \rightarrow \G(\R)$ be a Lie group automorphism.
    
    \noindent Then for any Borel set $\Omega \subset U^+$ having strictly positive Haar measure,
    \begin{equation}\label{eq:banana-basepoint}\pi \left( \overline{\; \bigcup_{i\in \N}   \rho\big( a(i) \; \Omega\big) \, \G(\Z) \;  } \right) = \G(\R)/\G(\Z).\end{equation}
 \end{theorem}

\begin{proof}[Proof of \cref{pr:strong_lal} (preparation)]

	We sketch the proof of \cref{pr:strong_lal}. 
	Without loss of generality, we assume that $S=S_0 \sqcup S_1 \sqcup S_2$ with $S_0 = \{p_0\}$ and $S_1 = \{p_{1},\ldots,p_{n}\}$ and $S_2=\{q_1, \ldots, q_n \}$ where $S$ and $x\in \Q_S^{n-1}$ satisfy the conditions of \cref{de:rfr}. 
	This can be done by working with a smaller subset of the group $\mathcal{H}_{\leq 0}$.
	
	Now let 
	\begin{equation}
	\label{eq:omega_1}
		\Omega_1 = \pi(\mathcal{H}_{\leq 0} \cdot \eta(x) \G(\mathcal{O}_S))    
	\end{equation}
	as given in the proposition.
	We construct $\hat{ A}_{x} \in \G(\mathbb{Q})$ such that the closed set $\Omega = \hat{A}_{x} \cdot \overline{\Omega_1}$ satisfies 
	the hypothesis of \cref{le:banana}. 
	
	First pick $\rho \in \mathbb{Z}_{\geq 1}$ so that the matrix $A_{x}$ given in \cref{de:rfr} 
	satisfies $\rho \cdot A_{x} \in \M_{(n+1) \times n}(\mathbb{Z})$.
	Then the matrix $\hat{A}_{x}$ is given by
	\begin{equation}\label{Ahat-x}
		\hat{A}_{x}=  \rho \cdot 
		\begin{pmatrix}
			1 & x_{p_{1}}^{(1)}& \dots  & x_{p_{1}}^{(n-1)}\\ 
			1 & x_{p_{2}}^{(1)} &\dots & x_{p_{2}}^{(n-1)} \\ 
			\vdots & & \ddots  & \\
			1 & x_{p_{n}}^{(1)} & & x_{p_{n}}^{(n-1)}
		\end{pmatrix},
	\end{equation}
	which is the matrix created from the lower $n$ rows of $\rho A_{x}$.
	We denote $\det( \hat{ A}_{x} ) = \gamma \in \mathbb{Z}$. By our assumption in \cref{de:rfr}, we know that $\gamma \neq 0$.

    For all $k\in \N^S$, consider an approximation $[x]_k\in \Z^{n-1}$ of $x$.
    We compute $\hat{A}_{x} \eta(-[x]_k)a(\prod_{p\in S} p^{k_p})$.
    We denote by $k=(k_0, k_1, \ldots, k_n, l_1, \ldots, l_n)$, by $\overline{k}= (k_1, \ldots, k_n, l_1, \ldots, l_n)$ and by $x_{p_i}=x_{q_i}=x_i \in \Q^{n-1}$ for all $1 \leq i \leq n$.
    
    \begin{align}
		\hat{A}_x \cdot \eta(-[x]_k)a\left(\prod_{p\in S} p^{k_p}\right)  
		& =  
		\rho  \cdot  
		\begin{pmatrix}
			1 & x_{1}^{(1)} - [x]_{k}^{(1)}& \dots  & x_{1}^{(n-1)} - [x]_{k}^{(n-1)}\\ 
			1 & x_{2}^{(1)} - [x]_{k}^{(1)} &\dots & x_{2}^{(n-1)} - [ x]_{k}^{(n-1)}\\ 
			\vdots & & \ddots  & \\
			1 & x_{n}^{(1)} - [x]_{k}^{(1)}& & x_{n}^{(n-1)} - [x]_{k}^{(n-1)}
		\end{pmatrix}
		a\left(  \prod_{p \in S} p^{k_{p}}\right) \\ 
		& =  
		\begin{pmatrix}
			p_{1}^{k_{1}}q_1^{l_1} & &   & \\ 
			& p_{2}^{k_{2}}q_{2}^{l_{2}} & & \\
			& & \ddots & \\
			& & & p_{n}^{k_{n}}q_{n}^{l_{n}}
		\end{pmatrix} 
		\cdot B(\overline{k},[x]_k) \cdot 
		\begin{pmatrix}
			p_0^{k_{0}} & &   & \\ 
			& 1 & & \\
			& & \ddots & \\
			& & & 1  
		\end{pmatrix} 
		, \\
		&= \diag(p_1^{k_1}q_1^{l_1}, \ldots,p_n^{k_n}q_n^{l_n} ) \cdot B(\overline{k},[x]_k) \cdot a(p_0^{k_0}),
		\label{eq:magic_multiplication}
	\end{align}
	where  $B(\overline{k},[x]_k)$ is the matrix given as 
	\begin{align}
		B(\overline{k},[x]_k):=\rho \cdot 
		\begin{pmatrix}
			\frac{\prod_{p\in S_1\cup S_2} p^{k_{p}}}{p_1^{k_1}q_1^{l_1} } 
			& \frac{x_{1}^{(1)} - [x]_{k}^{(1)}}{p_{1}^{k_{{1}}}q_1^{l_1} }
			& \dots  
			& \frac{x_{1}^{(n-1)} - [x]_{k}^{(n-1)} }{p_{1}^{k_{{1}}}q_1^{l_1} }\\ 
			&&&\\
			\frac{\prod_{p\in S_1\cup S_2} p^{k_{p}}}{p_2^{k_2}q_2^{l_2} }  
			& \frac{x_{2}^{(1)} - [x]_{k}^{(1)}}{p_2^{k_2}q_2^{l_2} } &\dots 
			& \frac{x_{2}^{(n-1)} - [x]_{k}^{(n-1)} }{p_2^{k_2}q_2^{l_2}}\\ 
			\vdots & & \ddots  & \\
			\frac{\prod_{p\in S_1 \cup S_2} p^{k_{p}}}{p_n^{k_n}q_n^{l_n} } 
			& \frac{x_{n}^{(1)} - [x]_{k}^{(1)} }{p_n^{k_n}q_n^{l_n}} &\dots 
			& \frac{x_{n}^{(n-1)} - [x]_{k}^{(n-1)} }{p_n^{k_n}q_n^{l_n}}\\ 
		\end{pmatrix}.
		\label{eq:garbage_matrix}
	\end{align}

Then, the diagonal invariance in \cref{prop:diag-inv} and the horospherical pieces of $\Omega$ are the contents of the rest of this section.
\end{proof}

	\subsubsection{Diagonal invariance}

\begin{proposition}\label{prop:diag-inv}
    Let $S$ and $x$ be as in \cref{pr:lal_high_rank} and $\hat{A}_x$ be as in \cref{Ahat-x}.
    Then the closed set $\Omega = \hat{A}_x \cdot \overline{ \pi(\mathcal{H}_{\leq 0} \cdot \eta(x) \G(\mathcal{O}_S)) }$
    is invariant under $\mathrm{A}(\R)$: namely
    $$ \mathrm{A}(\R) \cdot \Omega = \Omega,$$
    where $\mathrm{A}(\R) = \lbrace \diag(t_1, \ldots, t_n) \; : \; t \in (\R_{>0})^{n}\rbrace \subset \G(\R)$.
\end{proposition}

\begin{lemma}\label{le:diagonal_primes}
Let $p_{1},\ldots,p_{n}, q_1, \ldots, q_{n-1} \in \mathbb{P}$ be some distinct primes.
Then for every $M\in \Z_{\geq 1}$, the set 
$  \lbrace  \diag(p_1^{k_1}q_1^{l_1}, \ldots, p_{n-1}^{k_{n-1}}q_{n-1}^{l_{n-1}}, p_n^{k_n}) \; : \; (k,l)\in M\N^{2n-1} \rbrace $ is dense in $\mathrm{A}(\R)\subset \G(\R)$.
\end{lemma}
\begin{proof}
    Let $\diag(e^{t_1}, \ldots, e^{t_{n-1}}, 1 ) \in \mathrm{A}(\R)$ with $t\in \R^{n-1}$, let $\varepsilon>0$.
    We prove that there exists $(k,l)\in M \N^{2n-1}$ such that
    $$\Big\vert (k_1 \log p_1 + l_1 \log q_1 - k_n \log p_n, \ldots,k_{n-1} \log p_{n-1} + l_{n-1} \log q_{n-1} - k_n \log p_n) - t \Big\vert \leq \varepsilon.$$  
    For all $1 \leq i \leq n-1$, since $p_i$ and $q_i$ are distinct prime numbers, the additive subgroup generated by $\lbrace \log p_i, \log q_i \rbrace$ is dense in $\R$, so there exists $T_i>0$ such that $\lbrace k_i \log p_i + l_i \log q_i \; : \; k_i, l_i \in M\N \rbrace$ is an $\varepsilon$-net in $[T_i, + \infty)$.
    Denote by $T= \max_{i=1}^{n-1} T_i$.
    Then the set below is an $\varepsilon$-net in $[T , +\infty )^{n-1} $. 
    $$\lbrace (k_i \log p_i + l_i \log q_i)_{1 \leq i \leq n-1} \; : \; (k_i,l_i) \in M\N^2 \text{ for all } 1 \leq i \leq n-1 \rbrace.$$ 
    
    Since $[T , +\infty )^{n-1}=T \sum_{i=1}^{n-1} e_i + \sum_{i=1}^{n-1} \R_{\geq 0} e_i$ is a convex, non-degenerated cone, and since the vector $M\log p_n \sum_{i=1}^{n-1} e_i$ is in the interior of the cone $\sum_{i=1}^{n-1} \R_{\geq 0} e_i$, we can pull the $\varepsilon$-net so that it covers $t$: we choose $k_n\in M\N$ such that $\inf_{i=1}^{n-1} t_i  - 1 > T - k_n \log p_n$. \\
    Then $t \in (T-k_n \log p_n) \sum_{i=1}^{n-1} e_i + \sum_{i=1}^{n-1} \R_{\geq 0} e_i$ and 
    $$\lbrace (k_i \log p_i + l_i \log q_i - k_n \log p_n)_{1 \leq i \leq n-1} \; : \; (k_i,l_i) \in M\N^2 \text{ for all } 1 \leq i \leq n-1 \rbrace$$ is an $\varepsilon$-net in $(T-k_n \log p_n) \sum_{i=1}^{n-1} e_i + \sum_{i=1}^{n-1} \R_{\geq 0} e_i$. \\
    Therefore, there exists $(k_i,l_i)_{1 \leq i \leq n-1} \in (M\N^2)^{n-1}$ at distance less than $\varepsilon$ to $t$.
\end{proof}

We add another preparatory Lemma on properties of the integer matrices $B(\overline{k},[x]_k)$ introduced in \cref{eq:garbage_matrix}.

 \begin{lemma}
   \label{le:diagonal_equivalence_of_cosets}
   Consider $B(\overline{k},[x]_{k})$ the matrix defined in~\cref{eq:garbage_matrix}. 
   Recall that $\gamma = \det \hat{A}_x \neq 0$ which is the $n\times n$ minor extracted from $A_x$ in \cref{de:rfr} by removing the first line. 
   \begin{enumerate}

	\item $\det B (\overline{k} ,[x]_{k})= \gamma$ for all $k \in \mathbb{Z}_{\geq 0}^{S}$ and $[x]_k$.
	   
	\item \label{it:horosphere_calculation}
	For all $k\in (\Z_{\geq 0})^S$ and all $v,v'\in \R^{n-1}$
	$$B( \overline{k} ,v)=B( \overline{k} ,v')\eta\left(\frac{v'-v}{\prod_{p\in S_1 \cup S_2 }p^{k_p}} \right) $$
	\item For all $k\in (\Z_{\geq 0})^S$ and all $v,v'\in \R^{n-1}$
	$$ B(\overline{k},v) a(p_0^{k_0}) = B(\overline{k},v') a(p_0^{k_0}) \eta \left( \frac{v'-v}{\prod_{p\in S} p^{k_p} } \right).$$
	\item For all $v,v'\in \R^{n-1}$ such that $v-v'\equiv 0 \pmod{p^{k_p}}$ for all $p\in S$, 
	$$ B(\overline{k},v) a(p_0^{k_0}) \G(\Z) = B(\overline{k},v') a(p_0^{k_0}) \G(\Z).$$
	In other words, for every $k\in \N^S$ the map $v \mapsto B(\overline{k},v)$ is defined in $\big(\R/(\prod_{p\in S} p^{k_p} )\Z \big)^{n-1}$.  
	\item 
	\label{it:diagaonal_part}
	For all $k, k' \in \mathbb{Z}_{\geq 0}^{S}$ with $k_p=k'_p$ for $p\in S_0$ and $v_p(\gamma)\leq k'_p \leq k_p$ for all $p\in S_1 \cup S_2$, and all approximations $[x]_k$,
		\begin{align}
		    \overline{k}' \equiv \overline{k} &\pmod{\varphi(\gamma p_{0}^{k_{0}})} \implies \\
		    & \quad \quad B(\overline{k},[x]_{k}) a(p_{0}^{k_{0}})
		    \G(\mathbb{Z}) = B(\overline{k}',[x]_{k}) a(p_{0}^{k_{0}})\G(\mathbb{Z})
		\end{align}

   \end{enumerate}
   
 \end{lemma}
\begin{proof}
    
	\begin{enumerate}

    \item 
 	It follows from~\cref{eq:magic_multiplication} that
	$\det \hat{A}_x = \det B(\overline{k},[x]_k)$ for every choice of $(k_{p})_{p\in S}$.
    \item By matrix multiplication, we have
	    \begin{align}
	    &B(\overline{k},v')\eta\left(\frac{v'-v}{\prod_{p\in S_1 \cup S_2} p^{k_p} } \right)\\
	        =&\rho \cdot 
		\begin{pmatrix}
			\frac{\prod_{p\in S_1 \cup S_2} p^{k_{p}}}{p_1^{k_1}q_1^{l_1} } 
			& \frac{x_{1}^{(1)} - v'^{(1)}}{p_{1}^{k_{{1}}} q_1^{l_1} }& \dots  & \frac{x_{1}^{(n-1)} - v'^{(n-1)} }{p_{1}^{k_{{1}}} q_1^{l_1} }\\ 
			\frac{\prod_{p\in S_1 \cup S_2} p^{k_{p}}}{p_2^{k_2}q_2^{l_2} }
			& \frac{x_{2}^{(1)} - v'^{(1)}}{p_{2}^{k_{{2}}}q_2^{l_2} } &\dots & \frac{x_{2}^{(n-1)} - v'^{(n-1)} }{p_{2}^{k_{{2}}} q_2^{l_2} }\\ 
			\vdots & & \ddots  & \\
			\frac{\prod_{p\in S_1 \cup S_2} p^{k_{p}}}{p_n^{k_n}q_n^{l_n} } 
			& \frac{x_{n}^{(1)} - v'^{(1)} }{p_{n}^{k_{{n}}} q_n^{l_n} } &\dots & \frac{x_{n}^{(n-1)} - v'^{(n-1)} }{p_{n}^{k_{{n}}} q_n^{l_n}}\\ 
		\end{pmatrix}\eta\left(\frac{v'-v}{\Pi_{p\in S_1 \cup S_2 }p^{k_p}} \right)\\
		=&\rho \cdot 
		\begin{pmatrix}
			\frac{\prod_{p\in S_1 \cup S_2} p^{k_{p}}}{p_1^{k_1}q_1^{l_1} } 
			& \frac{x_{1}^{(1)} - v^{(1)}}{p_{1}^{k_{{1}}}q_1^{l_1}}& \dots  & \frac{x_{1}^{(n-1)} - v^{(n-1)} }{p_{1}^{k_{{1}}}q_1^{l_1}}\\ 
			\frac{\prod_{p\in S_1 \cup S_2} p^{k_{p}}}{p_2^{k_2}q_2^{l_2} }
			& \frac{x_{2}^{(1)} - v^{(1)}}{p_{2}^{k_{{2}}}q_2^{l_2}} &\dots & \frac{x_{2}^{(n-1)} - v^{(n-1)} }{p_{2}^{k_{{2}}}q_2^{l_2}}\\ 
			\vdots & & \ddots  & \\
			\frac{\prod_{p\in S_1 \cup S_2} p^{k_{p}}}{p_n^{k_n}q_n^{l_n} }
			 & \frac{x_{p_{n}}^{(1)} - v^{(1)} }{p_{n}^{k_{{n}}}q_n^{l_n}} &\dots & \frac{x_{p_{n}}^{(n-1)} - v^{(n-1)} }{p_{n}^{k_{{n}}}q_n^{l_n}}\\ 
		\end{pmatrix}
		= B(\overline{k},v).
	    \end{align}
	\item This identity follows from the fact that $a(p_0^{-k_0}) \eta(w) a(p_0^{k_0}) = \eta\Big(\frac{w}{p_0^{k_0}}\Big).$
	\item If $v-v' \equiv 0 \pmod{p^{k_p}}$ for all $p\in S$, then $\eta \left( \frac{v'-v}{\prod_{p\in S} p^{k_p} } \right) \in \G(\Z)$, hence the equality.
	
    \item 
    
 To show this, we want to use \cref{le:easy_hecke}. 
 
 We first study the matrix $B(\overline{k},[x]_k) a(p_{0}^{k_{0}}) \pmod{ \gamma p_{0}^{k_{0}}}$. 
 By the above result, changing the choice of $[x]_{k} \in \mathbb{Z}^{n-1}$ to $[x]_{k} + v \cdot \prod_{p \in S}{p^{k_{p}}}$ for some $v \in \mathbb{Z}^{n-1}$ will still preserve the $\G(\Z)-$coset.
Let us therefore replace $[x]_{k}$ by an even better approximation to establish our claim.
	
	Now consider the prime factorization of $\gamma p_{0}^{k_{0}}= \pm \prod_{p' \in T} p'^{r_{p'}}$ for some finite set of primes $T \subseteq \mathbb{P}$ and some exponents $r \in (\mathbb{Z}_{\geq 1})^{T}$. 
	We will enlarge the set of primes from $S$ to $S'=S\cup T$ and upgrade $x \in \mathbb{Q}^{S}$ to $x' \in \mathbb{Q}^{S'}$ 
	by simply setting for $p' \in T \setminus S \subseteq S',\; x'_{p'} = 0$ and $x_{p}' = x_{p}$ for $p \in S$. 
	Then any approximation $[x']_{k,r} \in \mathbb{Z}^{n-1}$ is also a valid replacement for the approximation $[x]_{k} \in \mathbb{Z}^{n-1}$.
	
	We now study the matrix $B(\overline{k},[x]_k) a(p_{0}^{k_{0}})$ modulo $\prod_{p' \in T} p'^{r_{p'}}$ by going through each prime in $T$.
	
	If $p' \in T \setminus S$, then
	for each $i \in \{0,\ldots,n\}$, the value of $p_{i}^{k_{i}}\pmod{p'^{r_{p'}}}$ depends only on $k_{i}$ modulo $\varphi(p'^{r_{p'}})$ since $p_{i}$ generates a finite cyclic group in $(\mathbb{Z}/p'^{r_{p'}} \mathbb{Z})^{\times}$.
	Furthermore, by replacing $[x]_{k}$ with $[x']_{k+r}$, we ensure that $[x']_{k+r} \equiv 0\pmod{p'^{r_{p'}}}$.
	Then by analyzing the expression in~\cref{eq:garbage_matrix}, fixing $(k_{1},\dots,k_{n})$ modulo $\varphi(p'^{r_{p'}})$ preserve the coset $B(\overline{k},[x]_{k}) a(p_{0}^{k_{0}}) \G(\mathbb{Z})$.
	
	Now suppose $p=p' \in S \cap T$.
	Then one has $[x']_{k+r} \equiv x_{p} \pmod{p'^{r_{p'}}}$. 
	The factor $\rho$ ensures that the matrix $\rho x_{q}$ is an integer for each $q \in S \setminus \{p\}$. 
	Furthermore, for each $q \in S \setminus \{p\}$, all the negative powers $q^{-k_{q}}$ exist modulo $p'^{r_{p'}}$ and depend only on the value of $k_{q}$ modulo $\varphi(p'^{r_{p'}})$. 
	On the other hand, for each $i \in \{1,\dots,n-1\}$the term $p^{-k_{p}}(x_{p}^{(i)}-[x']_{k+r}^{(i)}) \equiv 0\pmod{p^{r_{p}}}$ since $x_{p}-[x']_{k+r}$ is divisible by $p^{k_{p}+r_{p}}$ by construction.
	\end{enumerate}
\end{proof}

\begin{proof}[Proof of \cref{prop:diag-inv} (diagonal invariance)]

Consider $\hat{A}_x$ and $\gamma$ defined in \cref{ss:banana}.
We want to prove that $d y\in\Omega$ for all $d \in \mathrm{A}(\R)$ and $y \in \Omega$.
By \cref{le:equivalent_representation} and a diagonal argument, it is sufficient to assume that $y = \hat{A}_x \eta(-[x]_k) a(\prod_{p\in S}p^{k_p}) \G(\Z)$ for some $k\in \N^S$ with $k_p\geq v_p(\gamma)$ for $p\in S_1\cup S_2$, and $[x]_k \in \Z^{n-1}$ an approximation of $x$.
We find sequences $\lbrace z_m \rbrace \subset \Omega$ and $\lbrace d(m) \rbrace \subset \mathrm{A}(\R)$ such that $z_m=d(m) y$ where $d(m)\rightarrow d$.

Recall by \cref{eq:magic_multiplication} that
	\begin{equation}
		y = 
		\diag(p_{1}^{k_{1}}q_1^{l_1} , \dots ,p_{n}^{k_{n}}q_n^{l_n}) \cdot B(\overline{k},[x]_k)  \cdot a(p_{0}^{k_{0}}) \G(\mathbb{Z}).
	\end{equation}
Now for all $k' \in \lbrace 0\rbrace \times (\varphi(\gamma p_0^{k_0}) \Z_{\geq 0} )^{S_1 \cup S_2}$, by \cref{le:diagonal_equivalence_of_cosets}.\ref{it:diagaonal_part}, 
\begin{align}
B(\overline{k},[x]_k)  \cdot a(p_{0}^{k_{0}}) \G(\mathbb{Z}) &= B(\overline{k},[x]_{k+k'})  \cdot a(p_{0}^{k_{0}}) \G(\mathbb{Z}) \\
&= B(\overline{k}+\overline{k}',[x]_{k+k'}) \cdot a(p_{0}^{k_{0}}) \G(\mathbb{Z}) .
\end{align}
Hence, 
$$ \begin{array}{rl}
    \hat{A}_x \eta(-[x]_{k+k'}) \cdot a\Big(\prod_{p\in S}p^{(k+k')_p}\Big) \G(\Z)   & =  \diag(p_{1}^{k_1+k'_1}q_1^{l_1+l'_1} , \dots ,p_{n}^{k_n+k'_{n}}q_n^{l_n+l'_n}) \\
    & \quad \quad \quad  B(\overline{k}+\overline{k}',[x]_{k+k'})  \cdot a(p_{0}^{k_{0}}) \G(\mathbb{Z})  \\
     & =  \diag(p_{1}^{k'_{1}}q_1^{l'_1} , \dots ,p_{n}^{k'_{n}}q_n^{l'_n})  \\ 
    & \quad \quad \diag(p_{1}^{k_{1}}q_1^{l_1} , \dots ,p_{n}^{k_{n}}q_n^{l_n})
     B(\overline{k},[x]_{k})  \cdot a(p_{0}^{k_{0}}) \G(\mathbb{Z}) \\
     &= \diag(p_{1}^{k'_{1}}q_1^{l'_1} , \dots ,p_{n}^{k'_{n}}q_n^{l'_n}) \cdot y.
\end{array} $$
Using \cref{le:diagonal_primes}, we choose a sequence
$$
\lbrace k(m)\rbrace_{m\geq0}\subset
\lbrace0\rbrace\times(\varphi(\gamma p_0^{k_0})\Z_{\geq0})^{S_1\cup S_2}
$$
such that
$$
\diag\big(p_1^{k_1(m)}q_1^{l_1(m)},\ldots,p_n^{k_n(m)}q_n^{l_n(m)}\big)
\longrightarrow d
\qquad(m\to+\infty).
$$
We deduce diagonal invariance for
\begin{align}
    z_m &= \hat{A}_x \eta(-[x]_{k+k(m)}) \cdot a\left(\prod_{p\in S}p^{(k+k(m))_p}\right) \G(\Z) \\
    &= \diag\big(p_1^{k_1(m)}q_1^{l_1(m)},\ldots,p_n^{k_n(m)}q_n^{l_n(m)}\big) \cdot y \longrightarrow d \cdot y \text{ as } m \text{ goes to } + \infty.
\end{align}
\end{proof}

\subsubsection{Horospherical piece}

We prove the following statement.
\begin{proposition}\label{prop-horos-piece}
Let $\Omega$ be as in \cref{prop:diag-inv}.
Then there exists $F_0 \in \G(\R)$ such that
$$\Omega \supset \eta\big( \R^{n-1} \big) F_0 \G(\Z).$$
\end{proposition}

\begin{proof}    
We give an overview of the proof.
Recall that 
$$ \hat{A}_x \eta(-[x]_k)a\left(\prod_{p\in S}p^{k_p} \right) \G(\Z)= \diag(p_1^{k_1}q_1^{l_1},\ldots,p_n^{k_n}q_n^{l_n}) B(\overline{k},[x]_k) a(p_0^{k_0}) \G(\Z).$$
Denote by $\gamma = \gamma' \prod_{p\in S} p^{r_p}$ with $\gamma'$ coprime with any prime in $S$ and $(r_p)_{p\in S} \in (\Z_{\geq 0})^S$.

\begin{itemize}
    \item[(i)] In \cref{lem-solving-horosp-eq}, we construct solutions $(v, \overline{k})\in \Z^{n-1} \times  \N^{S_1\cup S_2}$ of the following equation
\begin{equation}
    B(\overline{k},[x]_k) a(p_0^{k_0})\G(\Z)
= \eta (v) B(\overline{r},[x]_k ) a(p_0^{k_0}) \G(\Z).
\end{equation}
\end{itemize}
Denote by $V(k_0)\subset \Z^{n-1}$ the projection of the set of such solutions e.g. such that for every $v\in V(k_0)$, there exists $k=k_0:\overline{k}$ such that 
\begin{align}
    \hat{A}_x \eta(-[x]_k)a\left(\prod_{p\in S}p^{k_p} \right) \G(\Z)
    &= \diag(p_1^{k_1}q_1^{l_1},\ldots,p_n^{k_n}q_n^{l_n}) \cdot
    a(p_0^{k_0})  \\ 
    & \quad \quad \quad
      \eta\bigg(\frac{v}{p_0^{k_0}} \bigg) \cdot a(p_0^{-k_0}) B(\overline{r},[x]_k) a(p_0^{k_0}) \cdot \G(\Z).
\end{align}
In the same \cref{lem-solving-horosp-eq}, we prove that $V(k_0)$ is additively invariant by $p_1^{r_1}q_1^{s_1}p_0^{k_0+r_0} \Z^{n-1}$.
By diagonal invariance and by \cref{le:diagonal_equivalence_of_cosets}, 4.), the problem reduces to proving that: 
\begin{itemize}
    \item[(ii)] $\frac{V(k_0)}{ p_0^{k_0}}$ equidistributes in $\R^{n-1}$ as $k_0$ goes to infinity in \cref{lem-using-bourgain};
    \item[(iii)] the following set is bounded in $\G(\R)/\G(\Z)$ in \cref{lem-non-escape2}.
$$\left\lbrace  a(p_0^{-k_0}) B(\overline{r},[x]_{k_0:\overline{r}}) a(p_0^{k_0}) \G(\Z) \right\rbrace_{k_0 \in \N} \subset  \G(\R)/\G(\Z).$$
\end{itemize}

Write
$$
z_{k_0}=a(p_0^{-k_0})B(\overline{r},[x]_{k_0:\overline{r}})
a(p_0^{k_0})\G(\Z).
$$
By (iii), there is a sequence $k_m\to\infty$ such that
$z_{k_m}\to F_0\G(\Z)$ for some $F_0\in\G(\R)$.

Let $\xi\in\R^{n-1}$. Using the density supplied by
\cref{lem-using-bourgain}, choose $v_m\in V(k_m)$ such that
$$
\frac{v_m}{p_0^{k_m}}\longrightarrow\xi.
$$
For each $m$, choose $\overline{k}^{(m)}$ such that
$(v_m,\overline{k}^{(m)})\in W(k_m)$. The point given by the
left-hand side of the identity above belongs to $\Omega$.
Multiplying the right-hand side on the left by the inverse of its
diagonal prefactor and using the diagonal invariance of $\Omega$, we obtain
$$
\eta\bigg(\frac{v_m}{p_0^{k_m}}\bigg)z_{k_m}\in\Omega.
$$
Taking the limit as $m\to\infty$ and using that $\Omega$ is closed gives
$$
\eta(\xi)F_0\G(\Z)
=\lim_{m\to\infty}
\eta\bigg(\frac{v_m}{p_0^{k_m}}\bigg)z_{k_m}
\in\Omega.
$$
Since $\xi\in\R^{n-1}$ was arbitrary, this proves the proposition.
\end{proof}

\paragraph{Solving the unipotent equation}
We use the criteria given in Corollary~\ref{co:difficult_hecke} and Corollary~\ref{co:less_difficult_hecke}.
Recall that the $(i,j)-$th minor of the matrix $P\in \mathrm{M}_n(\R)$ is denoted by $\Delta_{i,j}(P)$. 
Denote by $B(x_0)=\hat{A}_x \eta(-x_0)$ and by $\nu_p = val_{p}(\Delta_{1,1}(B(x_0)))$ for all prime number $p$.
Abusing notations, denote by $\nu_0=\nu_{p_0}$.

Take $\gamma = \det B(x_0)$.
Denote by $\hat{S}$ the subset of prime divisor of $\gamma$ not in $S$.
Recall that we took ${r} \in \N^{S  \cup \hat{S}}$ such that
$\gamma =  \prod_{p\in S\cup \hat{S}} p^{r_p}$ where $r_p$ are the coordinates of ${r}$.
Define $\overline{r} \in \N^{S_1 \cup S_2}$ to be the projection of $r$.
For every $2 \leq i \leq n$ denote by $\hat{\gamma}_i = p_1^{-r_1}q_1^{-s_1}p_i^{-r_i}q_i^{-s_i}  \gamma.$

\begin{lemma}\label{lem-solving-horosp-eq}
{

}
     For every $k_0 \geq  \nu_{0}$, consider $W(k_0)$ the set of tuples $(v,\overline{k}) \in \Z^{n-1} \times \N^{S_1 \cup S_2}$
such that
$$\left\lbrace 
\begin{array}{ccl}
    \overline{k} &\equiv& \overline{r} \pmod{\varphi(\gamma \prod_{p\in S\cup \hat{S}} p^{\nu_p} )}\quad\text{and}\quad \overline{k}\ge\overline{r}\quad \text{in }\N^{S_1 \cup S_2}\\
v^{(i)} &\equiv& (-1)^{i+1} \frac{\Delta_{i,1}( B(x_0)) }{ \Delta_{1,1} (B(x_0)) }  \bigg( \frac{p_i^{r_i}q_i^{s_i}}{p_1^{r_1}q_1^{s_1}}  - \frac{p_i^{k_i} q_i^{l_i}}{p_1^{k_1}q_1^{l_1}} \bigg) \pmod{p_0^{k_0} \hat{\gamma}_i } \text{ for all } 2 \leq i \leq n  \\
    v^{(i)} &\equiv& 0 \pmod{p_1^{r_1}q_1^{s_1} p_i^{r_i}q_i^{s_i} } \text{ for all } 2 \leq i \leq n 
\end{array}
\right\rbrace,$$
where we index $v = (v^{(2)},\cdots , v^{(n)})$.
Then for every $(v,\overline{k})\in W(k_0)$, 
$$B(\overline{k}, [x]_k) a(p_0^{k_0}) \G(\Z)= \eta(v) B(\overline{r},[x]_k) a(p_0^{k_0})\G(\Z). $$
\end{lemma}
By definition of $W(k_0)$, its projection to $\Z^{n-1}$ is well defined in $ \big( \Z/\gamma p_0^{k_0}\Z \big)^{n-1} $.
\begin{proof}
Consider the prime number decomposition $\gamma =  \prod_{p\in S\cup \hat{S}} p^{r_p} $.
By multiplicative properties of the totient function, $\varphi(\gamma) = \prod_{p\in S \cup \hat{S}} \varphi(p^{r_p}).$

First, note that $\overline{k}\equiv \overline{r} \pmod{\varphi( p^{r_p+\nu_{p}})}$ for every $p\in S \cup \hat{S}$.  
Hence for every $2 \leq i \leq n$, every $p\in (S \cup \hat{S}) \setminus \lbrace p_1,q_1,p_i,q_i\rbrace$  
$$ \frac{p_i^{r_i}q_i^{s_i}}{p_1^{r_1}q_1^{s_1}}  - \frac{p_i^{k_i} q_i^{l_i}}{p_1^{k_1}q_1^{l_1}} \equiv 0 \pmod{p^{r_p + \nu_{p} }} .$$
There are no primes in $p_0^{k_0} \hat{\gamma}_i$ outside of $(S \cup \hat{S}) \setminus \lbrace p_1,q_1,p_i,q_i\rbrace$.
Therefore, for every $2 \leq i \leq n$, the entry $v^{(i)}$ is an integer.
By Chinese remainder Lemma, adding the last equation about $v^{(i)}$, such an integer is well defined modulo $p_0^{k_0} \gamma$.

To simplify notations, denote by $B(\overline{k},[x]_k) = B(\overline{k})$ and $B(\overline{r},[x]_k) = B(\overline{r})$.
As in Corollary~\ref{co:difficult_hecke}, consider the line vector in $(\Q^{n-1})^t$ extracted from the matrix $w = \Big( \big( B(\overline{r})\big)_{1,j} \Big)_{2 \leq j \leq n} $.
We prove that elements of $W(k_0)$ satisfy the sufficient conditions given in Corollary~\ref{co:less_difficult_hecke}, namely that
\begin{itemize}
    \item[1.)] ${B}(\overline{k})-{B}(\overline{r}) \equiv 0 \pmod{\gamma}$,
    \item[2.)] $v \equiv 0 \pmod{\gamma}$,
    \item[3.)] 
$
    \Delta_{1,1}(B(\overline{k})) \big(v \; \submat B(\overline{r}) + w) \big)  $ \\  
\hspace*{1cm} $ \; + \; \big( (-1)^{2+1} \Delta_{2,1}(B(\overline{k})) ,\ldots, (-1)^{n+1}  \Delta_{n,1}(B(\overline{k})) \big) \; \submat B(\overline{r}) \equiv 0 \pmod{ p_0^{k_0+r_0} }$

\end{itemize}

\textbf{Verifying condition 1.):}
By analysing matrix coefficients of $B(\overline{k})$ and $B(\overline{r})$, condition 1.) it is sufficient to have 
$ p_i^{k_i}q_i^{l_i} \equiv p_i^{r_i} q_i^{s_i} \pmod{\gamma}$
for all $1 \leq i \leq n$ and $x_i - [x]_k \equiv 0 \pmod{p_i^{k_i+r_i}q_i^{l_i + s_i}}$ for all $1 \leq i \leq n$.

We check that for all $1 \leq i \leq n$, we have
 $p_i^{k_i}q_i^{l_i} \equiv p_i^{r_i} q_i^{s_i}$ modulo $p^{r_p}$ for all $p\in S \cup \hat{S}$.
Now
$\overline{k}\equiv \overline{r} \pmod{ \varphi(p^{r_p+\nu_p})}$ for all $p\in S \cup \hat{S}$.
Consequently, for every $1 \leq i \leq n$, 
we deduce that $p_i^{k_i} q_i^{l_i} \equiv p_i^{r_i}q_i^{s_i} \pmod{p^{r_p+\nu_p}}$ for every $p\in \lbrace p_0 \rbrace \cup \hat{S}$. 
Let us now deduce the same for every prime in $S_1 \cup S_2$.
Fix $1 \leq j \leq n$ and check the equations modulo $p_j^{r_j}$, then:
\begin{itemize}
	\item For every $1 \leq i \leq n$ with $i\neq j$, the prime $p_i,q_i$ are coprime with $p_j$, hence $p_i^{k_i}q_i^{l_i} \equiv p_i^{r_i}q_i^{s_i} \pmod{p_j^{r_j}}$ since $(k_i, l_i) \equiv (r_i,s_i) \pmod{\varphi(p_j^{r_j})}$. 
    \item When $i=j$, it is sufficient to have $\overline{k}_i \geq \overline{r}_i$ and $x_i - [x]_k \equiv 0 \pmod{p_i^{k_i+r_i}}.$
\end{itemize}
The above statements hold also for $q_j^{s_j}$, up to replacing  $p_j$ with $q_j$ in the first point and $p_i^{k_i+r_i}$ with $q_i^{l_i+s_i}$ in the third point.

\textbf{Verifying condition 2.):} 
Using that $\overline{k}\equiv \overline{r} \pmod{\varphi( \gamma \prod_{p\in S\cup \hat{S}}p^{\nu_p})}$, implies by properties of the totient function that $\overline{k}\equiv \overline{r}$ modulo $\varphi(p^{\nu_p+r_p} )$ for every $p\in S \cup \hat{S}$.

Since every prime in $S_1 \cup S_2 $ is coprime with $p_0$ as well as every prime in $\hat{S}$, we deduce that $\frac{p_i^{r_i}q_i^{s_i}}{p_1^{r_1}q_1^{s_1}}  - \frac{p_i^{k_i} q_i^{l_i}}{p_1^{k_1}q_1^{l_1}} \equiv 0 $ modulo $p_0^{\nu_{0} + r_0} \prod_{p\in \hat{S}}p^{r_p+\nu_p} $
for every $2 \leq i \leq n$.

Fix $2 \leq i \leq n$. 
It remains to check that $v^{(i)}$ is divisible by $p^{r_p+\nu_p}$ for every $p\in S_1 \cup S_2$.
For every $2 \leq j \leq n$ with $j\neq i$, since $p_j, q_j$ is coprime with $p_1,q_1,p_i,q_i$, we have that $p_j^{r_j+\nu_{p_j}} q_j^{s_j+\nu_{q_j}}$ divide $\frac{p_i^{r_i}q_i^{s_i}}{p_1^{r_1}q_1^{s_1}}  - \frac{p_i^{k_i} q_i^{l_i}}{p_1^{k_1}q_1^{l_1}}$. 
Consequently, $v^{(i)} \equiv 0$ modulo $p^{r_p+\nu_p} $ for all $p\in ( S \cup \hat{S}) \setminus \lbrace p_1, q_1, p_i, q_i \rbrace$. 
For $p_1,q_1, p_i,q_i$, we use the assumption that $p_1^{r_1}q_1^{s_1} p_i^{r_i}q_i^{s_i}$ divide $v^{(i)}$ and deduce condition 2.).

By \cref{le:diagonal_equivalence_of_cosets} 4.), we just need to choose an approximation $[x]_k \in \Z^{n-1}$ such that for every $p\in S_1 \cup S_2$ 
$$ x_p - [x]_k \equiv 0 \pmod{p^{k_p+r_p+\nu_p}}.$$
This will not change the projection point in $\G(\R)/\G(\Z)$.

\textbf{Verifying condition 3.):}
By definition of the adjoint matrix (using minors) and by recalling that $w$ is the non-diagonal part of the first line of the matrix $B(\overline{r})$, we deduce that
\begin{equation}\label{vector-inversion}
\begin{split}
   w \cdot \Big( \submat B(\overline{r}) \Big)^{-1} 
    &=  \Big(  B(\overline{r})_{1,2}, \ldots , B(\overline{r})_{1,n}  \Big)  \cdot \frac{ \Big( \submat B(\overline{r}) \Big)^{adj} }{ \Delta_{1,1}( B(\overline{r}) ) } \\
        &= -\left( (-1)^{2+1} \frac{ \Delta_{2,1}(B(\overline{r}))}{\Delta_{1,1}(B(\overline{r}))}, \ldots, (-1)^{n+1} \frac{ \Delta_{n,1}(B(\overline{r}))}{\Delta_{1,1}(B(\overline{r}))}  \right).
\end{split}
\end{equation}
Replacing the inverted vector, we deduce that condition 3.) is satisfied if 

\begin{equation}
    \begin{split}
         \Delta_{1,1}(B (\overline{k})) \; v
         & +  \Delta_{1,1}(B (\overline{k})) \left( (-1)^{2} \frac{ \Delta_{2,1}(B(\overline{r}))}{\Delta_{1,1}(B(\overline{r}))}, \ldots, (-1)^{n} \frac{ \Delta_{n,1}(B(\overline{r}))}{\Delta_{1,1}(B(\overline{r}))}  \right)    \\
        & - \Delta_{1,1}(B (\overline{k})) \left( (-1)^2 \frac{\Delta_{2,1}(B(\overline{k}))}{\Delta_{1,1}(B(\overline{k})) } ,\ldots, (-1)^n  \frac{\Delta_{n,1}(B(\overline{k}))}{ \Delta_{1,1}(B (\overline{k})) } \right) \;  \equiv 0 \pmod{p_0^{k_0+r_0}} .
    \end{split}
\end{equation}
By factorizing each line of the minor $\Delta_{i,1}(B(\overline{r}))$ where $1 \leq i \leq n$, we have
$$ \Delta_{i,1}(B(\overline{r})) = \frac{ p_i^{r_i}q_i^{s_i} }{\prod_{p\in S_1 \cup S_2}p^{r_p} } \Delta_{i,1}(B(0,[x]_k))  \text{    and    } \frac{\Delta_{i,1}(B(\overline{r})) }{ \Delta_{1,1}(B(\overline{r})) } = \frac{p_i^{r_i}q_i^{s_i}}{p_1^{r_1}q_1^{s_1}} \frac{\Delta_{i,1}( B(0, [x]_k)) }{ \Delta_{1,1} (B(0, [x]_k))  } .  $$
Similarly, for every $1 \leq i \leq n$,
$$\Delta_{i,1}(B(\overline{k})) = \frac{ p_i^{k_i}q_i^{l_i} }{\prod_{p\in S_1 \cup S_2}p^{k_p} } \Delta_{i,1}(B(0,[x]_k))
\text{    and    }
 \frac{\Delta_{i,1}(B(\overline{k})) }{ \Delta_{1,1}(B(\overline{k})) } = \frac{p_i^{k_i}q_i^{l_i}}{p_1^{k_1}q_1^{l_1}} \frac{\Delta_{i,1}( B(0, [x]_k)) }{ \Delta_{1,1} (B(0, [x]_k))  } .$$
Furthermore $x_0 \equiv [x]_k \pmod{p_0^{k_0+r_0}} $ implies that
$ B(0,[x]_k) \equiv B(x_0) \pmod{p_0^{k_0+r_0}} .$
So we can replace the minors of $B(0,[x]_k)$ with minors of $B(x_0)$ in the congruence equation.
Noting now that every $p\in S_1 \cup S_2$ is coprime with $p_0$ we further simplify.
Hence condition 3.) is satisfied if, for every $2 \leq i \leq n$,
$$ \Delta_{1,1}( B(x_0)) \left( v^{(i)} + (-1)^i \frac{\Delta_{i,1}( B(x_0)) }{ \Delta_{1,1} (B(x_0)) }  \bigg( \frac{p_i^{r_i}q_i^{s_i}}{p_1^{r_1}q_1^{s_1}}  - \frac{p_i^{k_i} q_i^{l_i}}{p_1^{k_1}q_1^{l_1}} \bigg) \right) \equiv 0 \pmod{p_0^{k_0+r_0}} .$$
Consequently, our assumptions on $v^{(i)}$ are sufficient to obtain condition 3.).
\end{proof}

\paragraph{Equidistribution of the unipotent solutions}
The main result of this paragraph is \cref{lem-using-bourgain}.
We adopt the same notations as \cref{lem-solving-horosp-eq}. 
Our proof is based on \cref{corol-bourgain}.
Here is a preparatory lemma.
\begin{lemma}

\label{lemma:index_of_subgroup}
Let $\hat{\gamma}, \varphi \in \N$ be non-zero integers. 
Let $p_0, p_1, q_1$ be distinct prime numbers such that $p_1, q_1 \nmid \hat{\gamma}$.
    Consider for every $k\geq 0$ the subgroup of $\Z/(\hat{\gamma} p_0^k\Z)$ given by
    $$ H_k=	\{   p_1^{k_1} q_1^{l_1}\in \Z/(\hat{\gamma} p_0^{k})\Z:\  k_1 = l_1=0\mod \varphi \; \}. $$
    Then there exists $c>0$ depending on $\hat{\gamma},p_0,p_1,q_1, \varphi$ such that $|H_k|\geq c \cdot p_0^k$ for all $k \in \Z_{\geq 0 }$.
    
\end{lemma}
\begin{proof}
Denote by $j_0 = \min \lbrace j \; : \; \vert H_j\vert>1 \rbrace.$
Then either $p_1^{\varphi}$ or $q_1^{\varphi}$ is non-trivial in $\Z/\hat{\gamma}p_0^{j_0}\Z$.
Denote by $x$ the non trivial element and by $d_0=d_0(\hat{\gamma},p_1,q_1,\varphi)$ the order of $x$ in $(\Z/\hat{\gamma} p_0^{j_0}\Z)^\times$.
Note that the order of $x$ in $(\Z/\hat{\gamma} p_0^{j_0+k}\Z)^\times$ is $d_0p_0^{k}$.
Indeed, there exists $m_0 \in \Z$, coprime with $\hat{\gamma}p_0^{j_0}$ such that $x^{d_0}  = 1 + \hat{\gamma}p_0^{j_0}m_0$.
Furthermore, for every $k\geq 1$ there exists $m_k\in \Z$ such that 
$x^{d_0 p_0^k} = 1 + \hat{\gamma}p_0^{j_0+k} m_0 + \hat{\gamma}^2 p_0^{2j_0}m_k.$
Consequently, $\vert H_{k}\vert \geq d_0 p_0^{k-j_0}$ for every $k\geq j_0$.
Choose $c= \min( d_0 p_0^{-j_0}, 1)$ so that $ \vert H_k\vert \geq c p_0^k$ for every $k\geq 0$.
\end{proof}
\begin{lemma}\label{lem-using-bourgain}
       For every $k_0 \geq  \nu_{0}$, denote by $V(k_0)$ the projection into $\Z^{n-1}$ of $W(k_0)$.
       
Then $V(k_0)/p_0^{k_0}$ is dense in $\R^{n-1}$ as $k_0 \rightarrow + \infty$.
\end{lemma}
\begin{proof}[Proof of \cref{lem-using-bourgain}.]
To lighten the notations, denote by $\varphi = \varphi( \gamma \prod_{p\in S\cup \hat{S}} p^{\nu_p} ) .$
For every $2 \leq i \leq n$, denote by $m_i=p_1^{r_1}q_1^{s_1}p_i^{r_i}q_i^{s_i}$, by
$$H_i(k_0) = \lbrace p_i^{k_i}q_i^{l_i}  \in \Z/\hat{\gamma}_i p_0^{k_0}\Z \; : \; k_i \equiv l_i \equiv 0 \pmod{ \varphi ) } \rbrace.$$
By \cref{lemma:index_of_subgroup}, since $p_i, q_i \nmid \hat{\gamma}_i$ there exists a constant $c_i $ depending only on $\hat{\gamma}_i, p_0, p_i, q_i, \varphi $ such that $\vert H_i(k_0)\vert\geq c_i \; p_0^{k_0} $.
Hence for $k_0$ large enough, $\vert H_i(k_0) \vert \geq \big(\hat{\gamma}_i p_0^{k_0}\big)^{1/2}$.

Consider the factorisation $\Delta_{1,1} B(x_0) = \Delta(i) \; \widehat{\Delta}(i) $
where $\Delta(i)$ is coprime with $p_0\hat{\gamma}_i$ and $\widehat{\Delta}(i)$ is in the multiplicative span of $p_0$ and $S\cup\hat{S}\setminus\{p_1,q_1,p_i,q_i\}$, e.g. the prime factors of $p_0\hat{\gamma}_i$.

For every $k_1, l_1 \in  \varphi\N$ and all $k_0\in \N$, consider
$$ X^{(i)}(k_0 ; k_1, l_1) =  \frac{(-1)^{i+1}\Delta_{i,1}( B(x_0)) \; p_i^{r_i}q_i^{s_i}}{ m_i \; \Delta(i)   \; p_1^{r_1}q_1^{s_1} }   \bigg(1  - \frac{H_i(k_0)}{ p_1^{k_1}q_1^{l_1}} \bigg) \subset \Z/\hat{\gamma}_i p_0^{k_0}\Z.$$
Denote by $\widetilde{X}^{(i)}(k_0; k_1, l_1)$ its preimage in $\Z$.
Note that $$ V(k_0) = \bigcup_{k_1, l_1 \in \N} \;  \prod_{i=2}^n \frac{ m_i }{\widehat{\Delta}(i) } \widetilde{X}^{(i)}(k_0; \varphi k_1, \varphi l_1) .$$

By \cref{corol-bourgain}, since $m_i, p_1, q_1, \Delta(i), p_i, q_i \in (\Z/ \hat{\gamma}_i p_0^k \Z)^\times$, the sets $X^{(i)}(k_0; k_1, l_1)$ equidistribute in $\R/\Z$ as $k_0$ goes to infinity and the speed of equidistribution is the same for every $k_1, l_1$.
Hence gaps between consecutive elements in $\widetilde{X}^{(i)}(k_0; \varphi k_1, \varphi l_1) $
are $o(\hat{\gamma}_i p_0^{k_0}) = o(p_0^{k_0})$ as $k_0$ goes to infinity.
Consequently, $V(k_0)/p_0^{k_0}$ is dense in $\R^{n-1}$.
\end{proof}

\paragraph{Non-escape of the starting point}

\begin{lemma} \label{lem-non-escape1}
Let $\overline{B}, F \in \G(\Q)$ be invertible matrices. 
Denote by $\overline{B}=\begin{pmatrix}
		b & v \\ t & B
	\end{pmatrix}$ where $t \in \Q^{n-1}$ and $v \in (\Q^{n-1})^t$.
Assume $\det B\neq 0$.

Then for every prime number $p$ the following sequence in $\G(\R)/\G(\Z)$ is bounded.
	$$  \Big( a(p^{-k}) \; \overline{B} \; a(p^{k}) F \G(\Z)\Big)_{k\in\N}  . $$
More precisely, all limit points are included in the following bounded set of $\G(\R)/\G(\Z)$:
$$ \begin{pmatrix}
b - v B^{-1}t & 0 \\
B [0,1]^{n-1} & B
\end{pmatrix} \G(\Z) F \G(\Z). $$
\end{lemma}

\begin{proof}
Note first that since $\overline{B}$ and $B$ are invertible, $\det \overline{B} = \big(b-vB^{-1}t\big) \det B \neq 0.$ 
Furthermore, $\G(\Z) F \G(\Z)$ is a Hecke orbit of the shape point $F \G(\Z)$, which is finite.
Therefore the set we claim contains all limit points is bounded.

It is sufficient to prove that limit points of $\big( a(p^{-k}) \; \overline{B} \; a(p^{k}) \G(\Z)\big)_{k\in\N}$ are in 
$$ \begin{pmatrix}
b - v B^{-1}t & 0 \\
B [0,1]^{n-1} & B
\end{pmatrix} \G(\Z). $$
First note that
\begin{equation}\label{eq:lem-non-esc}
    a(p^{-k}) \; \overline{B} \; a(p^{k}) = 
    \begin{pmatrix}
    b & p^{-k} v \\ p^k t & B
    \end{pmatrix}.
\end{equation}
Now consider the integer part $\lfloor B^{-1}p^k t \rfloor$ of the vector $B^{-1}p^k t$ and denote by $\lbrace B^{-1}p^k t \rbrace$ the fractional part so that 
$B^{-1}p^k t = \lfloor B^{-1}p^k t \rfloor + \lbrace B^{-1}p^k t \rbrace.$
Then note that 
\begin{equation}
    \begin{pmatrix}
    b & p^{-k} v \\ p^k t & B
    \end{pmatrix} 
    \begin{pmatrix}
    1 & 0 \\ - \lfloor B^{-1} p^{k} t \rfloor & I_{n-1}
    \end{pmatrix}
    = 
    \begin{pmatrix}
    b - v B^{-1}t + p^{-k} v \lbrace B^{-1}p^k t\rbrace & p^{-k} v \\ B \lbrace B^{-1}p^k t\rbrace & B
    \end{pmatrix}.
\end{equation}
Therefore 
$$  a(p^{-k}) \; \overline{B} \; a(p^{k}) \G(\Z) = \begin{pmatrix}
    b - v B^{-1}t + p^{-k} v \lbrace B^{-1}p^k t\rbrace & p^{-k} v \\ B \lbrace B^{-1}p^k t\rbrace & B
    \end{pmatrix} \G(\Z).  $$
We conclude by noting that the terms $ p^{-k} v \lbrace B^{-1}p^k t\rbrace $ and $p^{-k}v$ converge to $0$ as $k\rightarrow + \infty$, since $\lbrace B^{-1}p^k t\rbrace \in [0,1]^{n-1}$.    
\end{proof}

\begin{lemma}\label{lem-non-escape2}
For every $k\in \N^S$, recall that we denote by $\overline{k}=(k_p)_{p\in S_1 \cup S_2}$ and $k=k_0:\overline{k}$.
Independently from any choice of approximation $[x]_{k_0:\overline{r}}$, the following sequence is bounded.
    $$\left\lbrace  a(p_0^{-k_0}) B(\overline{r},[x]_{k_0:\overline{r}}) a(p_0^{k_0}) \G(\Z) \right\rbrace_{k_0 \in \N} \subset  \G(\R)/\G(\Z).$$  
\end{lemma}

\begin{proof}
By \cref{le:diagonal_equivalence_of_cosets} item 3, for every $k_0$ and every approximation $[x]_k$, 
$$B(\overline{r},[x]_k) a(p_0^{k_0}) =B(\overline{r}, x_0) a(p_0^{k_0}) \eta\bigg( \frac{x_0 - [x]_k}{ p_0^{k_0}\prod_{p\in S_1 \cup S_2} p^{r_p}} \bigg) .$$
Denote by $D>0$ a $p_0$-unit such that $D x_0\in\Z^{n-1}$, by $\hat{\gamma}_0=\prod_{p\in S_1 \cup S_2} p^{r_p}$, and by $d$ the multiplicative order of $p_0$ in $\Z/\hat{\gamma}_0\Z$.
We argue that $ \eta\big( \frac{x_0 - [x]_k}{p_0^{k_0}\hat{\gamma}_0} \big) \G(\Z)$ only depends on the residue modulo $d$ of $k_0$.
Moreover $D(x_0-[x]_k)/p_0^{k_0}\in\Z^{n-1}$, so all these unipotent parameters lie in $(D\hat{\gamma}_0)^{-1}\Z^{n-1}$, a finite set modulo $\Z^{n-1}$.
Denote by $j(k_0)\in [0,d-1]\cap \Z$ such that $k_0 \equiv j(k_0) \pmod{d}$. Indeed using that $[x]_{j(k_0):\overline{r}} \equiv [x]_{k} \equiv x_p \pmod{p^{r_p}}$ for every $p\in S_1 \cup S_2$ and $[x]_{j(k_0):\overline{r}} \equiv [x]_{k} \equiv x_0 \pmod{p_0^{j(k_0)}}$ since $k_0 \geq j(k_0)$,
\begin{equation}
 \begin{split}
\eta\bigg( \frac{x_0 - [x]_k}{p_0^{k_0}\hat{\gamma}_0} \bigg) \G(\Z)
& = 
\eta\bigg( \frac{x_0 - [x]_{j(k_0):\overline{r}}}{ p_0^{j(k_0)}\hat{\gamma}_0} \bigg) \\
&\quad\quad \quad \eta\bigg( \frac{[x]_{j(k_0):\overline{r}} - [x]_{k}  }{p_0^{j(k_0)}\hat{\gamma}_0} \bigg) 
\eta\bigg( \frac{ (x_0 - [x]_{k} )/p_0^{k_0} -  (x_0 - [x]_{k}  )/p_0^{j(k_0)}}{\hat{\gamma}_0} \bigg) 
\G(\Z) \\
&= \eta\bigg( \frac{x_0 - [x]_{j(k_0):\overline{r}}}{ p_0^{j(k_0)}\hat{\gamma}_0} \bigg) \G(\Z).
\end{split}
\end{equation}
Hence, for every $k_0\in \N$
$$a(p_0^{-k_0}) B(\overline{r},[x]_{k_0:\overline{r}}) a(p_0^{k_0}) \G(\Z) = a(p_0^{-k_0}) B(\overline{r}, x_0) a(p_0^{k_0}) \; \eta\bigg( \frac{x_0 - [x]_{j(k_0):\overline{r}}}{ p_0^{j(k_0)} \hat{\gamma}_0 } \bigg)\G(\Z). $$
We conclude using the previous \cref{lem-non-escape1} on $B(\overline{r},x_0),\eta\big( \frac{x_0 - [x]_{j:\overline{r}}}{ p_0^{j} \hat{\gamma}_0 } \big) \in \G(\Q)$ for every fixed integer $0 \leq  j \leq d-1$. 
\end{proof}

\begin{proof}[Proof of \cref{pr:strong_lal}]
Let
$$
\Omega_1=\pi(\mathcal{H}_{\leq 0}\cdot\eta(x)\G(\mathcal{O}_S)),
\qquad
\Omega=\hat{A}_x\cdot\overline{\Omega_1}.
$$
By \cref{prop:diag-inv}, $\Omega$ is closed and invariant under $\mathrm{A}(\R)$.
By \cref{prop-horos-piece}, there exists $F_0\in\G(\R)$ such that
$$
\eta(\R^{n-1})F_0\G(\Z)\subseteq\Omega.
$$

For every $i\in\N$, set $a(i)=a(e^i)$. For this sequence we have
$U^+=\eta(\R^{n-1})$. Fix
$$
E=\eta([0,1]^{n-1}),
$$
which has positive Haar measure in $U^+$. Applying
\cref{le:banana} with $\rho(g)=F_0^{-1}gF_0$, after conjugating $a(i)$, $E$ and $\G(\Z)$ by $F_0^{-1}$ gives
$$
\pi\left(\overline{\bigcup_{i\in\N}a(i)E F_0\G(\Z)}\right)
=\G(\R)/\G(\Z).
$$
For every $i$, we have $E F_0\G(\Z)\subseteq\Omega$, and hence
$a(i)E F_0\G(\Z)\subseteq a(i)\Omega=\Omega$.
Since $\Omega$ is closed, it follows that
$\Omega=\G(\R)/\G(\Z)$. Therefore
$$
\hat{A}_x\cdot\overline{\Omega_1}=\G(\R)/\G(\Z),
\qquad\text{so}\qquad
\overline{\Omega_1}=\G(\R)/\G(\Z).
$$
Thus $\Omega_1$ is dense in $\G(\R)/\G(\Z)$, as claimed.
\end{proof}

\section{Proof of the main theorem}
\label{se:main-proof}

First we fix the family of input parameters:
\begin{enumerate}
    \item fix $S$ a set of $2n+1$ prime numbers that are inert for the Selmer polynomial with $p>C_n$ for all $p\in S$;

    \item for all $1 \leq i \leq n$, denote by 
$c(p_i)=c(p_{n+i})= -e_i$ and $c(p_0)= \sum_{j=1}^n e_j$;

    \item for all $p\in S$, choose a generator $\zeta_{\sel, p}$ of the group of unities of $\calO^{\sel} \otimes \Z_p$.
\end{enumerate}

{We show there is an output family $\lbrace (p,l, \epsilon_p,\phi_p,\psi_p), \; p\in S \rbrace$, a sequence $(a_1(N),\ldots,a_n(N))_{N\in \N}$ such that
\begin{itemize}
    \item the shape of $\log \calO^\times$ converges to a point in $\PGL_n(\R) / \PGL_n(\Z)$;
    \item for each depth $N$, the tuple $a(N)$ satisfies input family $\lbrace (p,c(p), \zeta_{\sel,p}), \; p\in S \rbrace$ up to depth $N$; 
    \item for all $k_S \in \lbrace 1, \ldots, N \rbrace^S$, the units $\calO_{k_S}^\times=\cap_{p\in S} ( \Z + p^{k_p} \calO ) $ are extracted from $\calO^\times$ using the output family at depth $k_S$;
    \item the output congruence conditions $\lbrace (p, l_p, \epsilon_p, \psi_p, \phi_p, k_p) : p\in S, \text{ and } 0 \leq  k_p\leq N \rbrace_{p\in S}$ only allow depth to change.
\end{itemize}
}

Fix a large enough constant $C$ so that condition 1.) of \cref{prop:separation-condition} holds.
Using Chinese remainder Lemma, choose $(a_0,\cdots,a_{n})\in \Z^{n+1}$ that satisfies every input congruence condition $(p,c_{(p)},\zeta_{\sel,p}, N)$ for every $p\in S$ and such that for every $1 \leq j \leq n$
$$ a_{j}-a_{j-1} > C \sqrt{a_n}.$$
Now as in \cref{def-canonical-M-basis}, consider $f(X)= \prod_{i=0}^n (X-a_i) - 1$, denote by $\calO=\Z[X]/\langle f\rangle$ the monogenic order of the number field $K=\Q[X]/\langle f\rangle$ and by 
$b_i = X-a_i$ for all $1 \leq i \leq n$.
Choose a canonical shape change $M \in \tilde{\calM}$ and generators $\beta_1,\ldots,\beta_n$ of $\calO^\times$ such that for all $1 \leq i \leq n$
$$ b_i = \prod_{j=1}^n \beta_j^{m_{ij}}.$$

Now for each prime $p\in S$, apply \cref{prop:structure-suborders} for this choice of integer tuple $(a_0,\ldots, a_n)$, for input parameters $(p, c(p), \zeta_{\sel,p}, N)$.
There exists output parameters $(p,N,l_p,\epsilon_p,\psi_p,\phi_p)$ where $0\leq\epsilon_p<p^{n+1}-1$ and $\psi_p$, $\phi_p$ are linear forms on $\Z^n$
such that for every $0 \leq k_p \leq N$, 
$$\calO_{p^{k_p}}^\times=\bigg\{\pm \prod_{l=1}^n\beta_l^{x_l},\ x\in L_p(k_p-l_p) \bigg\} $$
where for every $k_p \in \Z$, 
\begin{equation}\label{equ:lik}
	L_p(k_p)=\left\lbrace x \in \Z^n \; \bigg\vert
	\begin{array}{cl}   
	\psi_p(x)&=0\mod p^{ \max (k_p, 0) }  \\
	\epsilon_p p^{l_p}\psi_p(x)+\phi_p(x)&=0\mod m_p,\qquad m_p=\frac{p^{n+1}-1}{p-1}
	\end{array} 
	\right\rbrace.
\end{equation}

As in the proof of \cref{prop:structure-suborders}, consider the map 
$$p \mapsto v(p) = \left\lbrace
\begin{array}{ll}
    c(p) & \quad \text{if  } p=p_0\\
    (n+1)c(p) + \sum_{i=1}^n e_i & \quad \text{otherwise}
\end{array}\right.   $$ 
Recall that the entries of $\phi_p$ are positive integers that divide $\gcd(\det M,  p^{n+1}-1)$ and that 
\begin{equation}\label{equ:psi-i}
\psi_{p}= \frac{1}{\gcd(M^{\adj}v(p))} M^{\adj}v(p).
\end{equation}

Since all the primes in $S$ are distinct, then for every $(k_p)_{p\in S} \in \N^S$ such that $l_p \leq k_p \leq N$ for all $p\in S$, we have
\begin{equation}\label{equ:suborders}
\calO_{\prod_{p\in S} p^{k_p}}^\times=\bigcap_{p\in S}\calO_{ p^{k_p}}^\times=\bigg\{\pm \prod_{i=1}^n \beta_i^{x_i},\quad  x\in \bigcap_{ p\in S}L_p(k_p-l_p) \bigg\}.  
\end{equation}

For a tuple $k_S\in\N^S$, we define 
\begin{equation}\label{equ:lk}
L(k_S)=\bigcap_{p\in S}L_p(k_p).
\end{equation}

\begin{proposition}\label{prop:density-lk}
The set of all lattices $L(k_S)$ where $k_S\in\N^S$ is dense in $\PGL_n(\R)/\PGL_n(\Z)$. 
\end{proposition}

\begin{proof}
	
For the purpose of density it suffices to consider $k_S\in\N^S$ with $k_p\geq l_p$ for every $p\in S$.

	We separate the congruence equations into two parts, one part that involves $k_S$, the system 
	$$\lbrace \psi_p(x)=0\pmod{ p^{k_p} } \text{ for all }  p\in S \rbrace, $$ and the other part which is independent of $k_S$: the system
	$$\lbrace \epsilon_p p^{l_p} \psi_p(x)+ \phi_p(x)=0\pmod{ m_p } \text{ for all }  p\in S \rbrace .$$

Applying \cref{lem:density-constant-matrix} to the $n\times S$
integer matrix $\Phi=(\epsilon_p p^{l_p}\psi_p+\phi_p)_{p\in S}$,
there exists an integer matrix $\gamma\in\mathrm{M}_n(\Z)$ such that
$\det\gamma\neq0$ and such that $\gamma\Z^n$ is the solution to
$$
(\Phi)^t x\in\diag(m_p)_{p\in S}\Z^S.
$$
	
	Then, note that $L(k_S)=\gamma J(k_S)$ where
	$$ J(k_S)=\{ x\in\Z^n,\  \psi_p(\gamma (x))=0\mod p^{k_p}, \text{ for all } p \in S \} .$$
    By \cref{lem:using-S-arithmetics}, the family $\lbrace J(k_S) \rbrace_{k_S \in \N^S}$ is dense, hence the density of $\lbrace L(k_S) \rbrace_{k_S \in \N^S}$.
\end{proof}
\begin{lemma}\label{lem:density-constant-matrix}
For all integer matrix $\Phi \in \mathrm{M}_{n,S}(\Z)$, there exists a square matrix of integer entries $\gamma \in \mathrm{M}_n(\Z)$ such that $\det \gamma \neq 0$ and 
$$\lbrace x \in \Z^n, \; (\Phi)^t x \in \diag(m_p)_{p\in S} \Z^S \rbrace = \gamma \Z^n.$$
\end{lemma}
\begin{proof}
We solve the following linear diophantine equation where $x \in \Z^n$ is unknown and denote by $(e_{p})_{p\in S}$ the canonical basis of $\Z^{S}$,
$$(\Phi)^t x \in  \diag(m_p)_{p\in S} \Z^S.$$
We argue by induction on the cardinality of $S$. The case where $\Phi=0$ is solved by $\gamma=I_n$.
Assume $S$ contains only one prime $p$ and $\Phi\neq 0$.
Apply Smith normal form, there exists $V \in \SL(n,\Z)$ and $d_p\in \Z_{\geq 1}$ such that $(\Phi)^t = d_p e_p^t \;  V$ and denote by $\gamma = V^{-1}a\big( \frac{ \mathrm{lcm}(d_p, m_p)}{d_p} \big)$.
It follows that $\det \gamma\neq 0$ and $\gamma\Z^n$ is the solution to the above equation.

If $S$ contains $n+1$ distinct primes for some integer $n\geq 1$, fix $p\in S$.
By induction, we solve the equation for $(\Phi_q)_{q \in S \setminus \{ p \} } \in \mathrm{M}_{n, S \setminus \{ p \} }(\Z)$ and we denote by $\gamma_n\Z^n$ its solution.
Note that $$ \lbrace x \in \Z^n , \; (\Phi)^t x \in \diag(m_q)_{q\in S} \Z^S \rbrace 
= \gamma_n \lbrace y \in \Z^n , \; (\Phi_p)^t\; (\gamma_n y) \in \Z m_p e_p \rbrace.$$
Again, by induction, we solve the system of equation (in $y$) for $ (\gamma_n)^t \; \Phi_p$ and denote by $\gamma_0\Z^n$ the solution.
Consequently, $\gamma_n \gamma_0 \Z^n$ is the solution to the equation for $n+1$ primes where both $\gamma_n$ and $\gamma_0$ have non-zero determinant.
\end{proof}

\begin{lemma}\label{lem:using-S-arithmetics}
      
	For any integer matrix $\gamma \in \mathrm{M}_n(\Z)$ with non-zero determinant, 
	for all $k_S\in \N^S$, denote by
	$$J(k_S) = \lbrace x \in \Z^n , \; \psi_p(\gamma x) = 0 \pmod{p^{k_p}} , \; \text{for all } p\in S \rbrace.$$
	Then the system of $S$ equations satisfy the highly full rank condition and the set $ \lbrace J(k_S) \rbrace_{k_S \in \N^S}$ is dense in $\PGL(n,\R) / \PGL(n,\Z)$.
\end{lemma}	

\begin{proof}
	We partition $S=S_0 \sqcup S_1 \sqcup S_2$ where $S_0 = \lbrace p_0\rbrace$ and $S_1 = \lbrace p_1 , \ldots, p_n \rbrace$ and $S_2 = \lbrace p_{n+1}, \ldots, p_{2n} \rbrace$.
	We relabel the elements in the following way $S_2= \lbrace q_1, \ldots, q_n \rbrace$.
	From \cref{equ:psi-i}, we know the vectors  $(\psi_p)_{p\in S_0 \cup S_1}$ arranged in a $(S_0 \cup S_1)\times n$ matrix form the following full rank matrix
	$$A_*=
	\begin{pmatrix}
	h_0^{-1} &    & \\
	 & \ddots &  \\
	 &  & h_n^{-1}  
	\end{pmatrix}
	\begin{pmatrix}
	1 &  \cdots  & 1\\
	-n &  & (1) \\
	 & \ddots  &   \\ 
		 (1) &  & -n  
	\end{pmatrix} (M^{adj})^t \gamma,$$
	where $h_i = \gcd(M^{adj} v(p_i) )$ for all $0 \leq i \leq n$.
	We obtain the same matrix $A_*$ by arranging the vectors in $(\psi_p)_{p\in S_0 \cup S_2}$ in a $(S_0 \cup S_2) \times n$ matrix.
	Indeed, $M^{adj}$ and $\gamma$ are square integer matrix of non-zero determinant and the non-square factor is full rank.
	
	Note that using Chinese remainder and relabeling of $S_2$, for all $k_0 \in \N$ and $k_{S_1} \in \N^{S_1}$ and $l_{S_2} \in \N^{S_2}$, 
$$	\begin{array}{rl}
    J(k_0, k_{S_1}, l_{S_2})
    &=\{ x\in\Z^n ,\; \psi_p(\gamma x)=0\mod p^{k_p}, \text{ for all }p\in S\} \\
    &=\{ x\in\Z^n ,\; A_* x \in \diag(p_0^{k_0}, p_1^{k_1}q_1^{l_1}, \ldots, p_n^{k_n}q_n^{l_n} ) \Z^{S_0 \cup S_1}  \}.
	\end{array}$$

In order to apply \cref{le:equivalent_representation}, we need to find $g\in \SL(n,\Z)$ and a vector $r\in \Q^{n+1}$ of non-zero, rational entries such that there exists $x\in \Q_{S_0 \cup S_1}^{n-1} \cap \Q_{S_0 \cup S_2}^{n-1}$ rational, strongly full rank such that
$$ \begin{pmatrix}
r_0 &  & \\
& \ddots & \\
& & r_n
\end{pmatrix} A_* g = 
\begin{pmatrix}
1 & x_{0}^{(1)} & \cdots & x_{0}^{(n-1)} \\
\vdots & \vdots & & \vdots \\
1 & x_{n}^{(1)} & \cdots & x_{n}^{(n-1)}
\end{pmatrix} .$$
There exists $v\in \Z^n$ with $\gcd(v)=1$ such that $A_*v$ has only non-zero entries.
Otherwise, we would have $A_*\Z^n \subset \cup_{i=0}^n \lbrace x_i=0 \rbrace$. Since $A_* \Z^n$ is a $\Z$ module of rank $n$, it would imply in $A_* \Z^n \subset \lbrace x_i=0 \rbrace$ for some $0 \leq i \leq n$, meaning that the $i-$th line of $A_*$ is zero. This is in contradiction with the fact that $A_*$ is full rank.

We choose $v$ a bit more carefully as follows. 
For each $p\in S$, write the corresponding row of $A_*$ as $d_pu_p$, with $d_p>0$ and $u_p$ a primitive linear form.
By CRT, we may arrange $u_p(v)\not\equiv0\pmod p$ for every $p\in S$.
Since $\gcd(v)=1$, we can find a generating basis of the $\Z$-module $\Z^n$ containing $v$ and denote by $g\in \SL(n,\Z)$ the corresponding matrix whose first column is $v$.
Set $r_i = (A_* v)_i^{-1}$ for all $0 \leq i \leq n$.
Apply \cref{le:equivalent_representation} on $A_x=\diag(r_0, \ldots,r_n) A_* g$ where $x\in \Q^{n-1}$ is rational and strongly full rank in the sense of Definition~\ref{de:rfr}.
For $t_S\in\N^S$, put $k_p=t_p+v_p(d_p)$. Then \cref{le:equivalent_representation} identifies $g^{-1}J(k_S)$ with the lattice at $\pi(a_S(-t_S)\eta(x)\G(\calO_S))$. Since $g$ is fixed, \cref{pr:strong_lal} proves the required density.
\end{proof}

\begin{proof}[Proof of \cref{th:main_orders}]
We now finish the proof of \cref{th:main_orders}. Fix a nonempty open set $U\subset\PGL_n(\R)/\PGL_n(\Z)$.

\begin{enumerate}
\item[\bf Step 1.]
For each $N$, choose $a_1,\ldots,a_n$ satisfying the hypotheses of \cref{prop:separation-condition} and the input congruence conditions up to depth $N$, and write $a(N)=(a_1(N),\ldots,a_n(N))$. With $\log$ as in \cref{eq:logmap}, define $V_N=(\log b_1,\ldots,\log b_n)$ and $H_N=(\log\beta_1,\ldots,\log\beta_n)$ in the same coordinates on $\mathfrak a$. Then $V_N=H_NM^t$.
\item[\bf Step 2.]
Apply \cref{prop:structure-suborders} at each prime $p\in S$. This produces output parameters $l_p,\epsilon_p,\psi_p,\phi_p$ such that the unit lattices of the resulting suborders are given by \cref{equ:lik,equ:suborders}. Only finitely many choices are needed for $M$ and for these parameters.
\item[\bf Step 3.]
Pass to a subsequence on which $M$ and all the output parameters are constant. By the proof of \cref{prop:separation-condition}, the entries of $V_N/\log a_n(N)$ are bounded and $|\det V_N|\gg_n(\log a_n(N))^n$. After passing to a further subsequence, $V_N/\log a_n(N)\to V\in\mathrm{GL}_n(\R)$. Since $M$ is constant, $H_N/\log a_n(N)\to H=VM^{-t}$ as matrices in the same coordinates.
\item[\bf Step 4.]
By \cref{prop:density-lk}, the set of all $L(t)$ for $t \in \N^{S}$ is dense in $\PGL_n(\R)/\PGL_n(\Z)$. Choose $t\in\N^S$ such that $HL(t) \cdot \PGL_{n}(\mathbb{Z})\in U$. 
Since $H_N/\log a_n(N)\to H$, we may advance the subsequence until $N\geq l_p+t_p$ for every $p\in S$ and $H_N/\log a_n(N)L(t)\cdot\PGL_n(\Z)\in U$.
By \cref{equ:suborders}, this point is the shape of $\mathcal O_{(l+t)}$ in \cref{eq:orders}. Therefore, we have found a shape point lying in $U$.
\end{enumerate}

\end{proof}

\section{Positive diagonal action and totally positive units}

In this section, we will consider a slight modification of the main theorem.

Consider the group of diagonal matrices $A_{+} \subseteq \mathrm{SL}_{n+1}(\mathbb{R})$ given as 
\begin{equation}
    A_{+} = \{ 
\begin{bmatrix}
\lambda_{0} &  & &  \\
&  \lambda_{1} & &  \\ 
& & \ddots & \\ 
& & & \lambda_{n}
\end{bmatrix} \mid \prod_{i=0}^{n} \lambda_{i} = 1, \lambda_{i} > 0
    \}.
\end{equation}
The Lie algebra of $A_{+}$ is still $\mathfrak{a} \simeq \mathbb{R}^{n}$ and in fact the logarithm map is now an isomorphism 
\begin{align}
\log A_{+} & \rightarrow   \mathfrak{a} \simeq \mathbb{R}^{n} \\
\diag(a_{0},\dots,a_{n})  & \mapsto   (\log a_{i})_{i=0}^{n}
\end{align}

For this setup, the compact orbits of $A_{+} \curvearrowright \mathrm{SL}_{n+1}(\mathbb{R})/\mathrm{SL}_{n+1}(\mathbb{Z})$ are slightly different from the $A$-action. Firstly, one can generalize the definition of shapes in the following form.
\begin{defin}
For a compact $A_{+}$-orbit $F \subseteq \mathrm{SL}_{n+1}(\mathbb{R})/\mathrm{SL}_{n+1}(\mathbb{Z})$, one defines $\mathcal{P}_{+}(F)$ to the set of periods given as 
\begin{equation}
    \mathcal{P}_{+}(F) = \{ a \in A_{+} : a g \mathrm{SL}_{n+1}(\mathbb{Z}) = g \mathrm{SL}_{n+1}(\mathbb{Z}), \text{ for all } g \mathrm{SL}_{n+{1}}(\mathbb{Z}) \in  F \}  .
\end{equation}
The shapes $\sh_{+}(F)$ of $F$ is denoted to be the lattice $\log \mathcal{P}_{+}(F) \subseteq \mathbb{R}^{n}$ identified as a point in $\PGL_{n}(\mathbb{R})/\PGL_{n}(\mathbb{Z})$.
\end{defin}

To get the number theoretic object equivalent to $\sh_{+}(F)$, one has the following definition.

We define the shape of an order as the following. Given a number field $K$ and an order $\mathcal{O} \subseteq K$, we denote 
\begin{equation}
\mathcal{O}^{\times}_{+} = \{ u \in  \mathcal{O}^{\times} \mid \sigma(u) > 0 \text{ for every embedding }\sigma: K \rightarrow  \mathbb{R}\}
\end{equation}
It is clear that if $u \in \mathcal{O}^{\times}, u^{2} \in \mathcal{O}^{\times}_{+}$. Hence $\mathcal{O}^{\times}_{+} \subseteq \mathcal{O}^{\times}$ is a finite index subgroup.
Then, we modify the definition of $\sh(F)$ to count only the units in $\mathcal{O}^{+}$.
\begin{defin}
Given a totally real number field $K$ of degree $n+1$,
an ordered set of embeddings $\sigma_{0},\dots,\sigma_{n}$ and an order $\mathcal{O} \subseteq K$, 
one denotes $\sh_+(\mathcal{O})$ to be the point $h \PGL_{n}(\mathbb{Z}) \in  \PGL_{n}(\mathbb{R})/\PGL_{n}(\mathbb{Z})$  
such that the lattice $\log (\mathcal{O}^{\times}_{+}) =  h \cdot \mathbb{Z}^{n}$, where $h \in \GL_n({R})$ is some matrix and $\log : \mathcal{O}^{\times}\rightarrow \mathbb{R}^{n}$ is the map in~\cref{eq:logmap} after identifying $\mathfrak{a} \simeq \mathbb{R}^{n}$ with respect to the bijection in~\cref{eq:lin_bijection}.
\end{defin}

Then, one observes the following equivalence. We skip the proof because it follows from the proof of~\cref{pr:equivalence_of_shapes} with minor modifications.
\begin{proposition}
  The following two subsets of $\PGL_n(\mathbb{R})/\PGL_n(\mathbb{Z})$ are equal.
  \begin{enumerate}
      \item $\{\sh_{+}(F)  \mid F \text{ compact $A_{+}$-orbits in }\SL_{n+1}(\mathbb{R})/\SL_{n+1}(\mathbb{Z}) \}$,
      \item
      $$\begin{aligned}
      \{\sh_{+}(\mathcal{O})\mid K\text{ totally real subfield with an ordered set of embeddings},\\
      \mathcal{O}\subseteq K\text{ some order},\ \deg K=n+1\}.
      \end{aligned}$$
  \end{enumerate}
\end{proposition}

We would like to argue that the following theorem is also true.
\begin{theorem}
For $n\geq3$, the following set is dense in $\PGL_{n}(\R)/\PGL_{n}(\Z)$
$$\{\sh_{+}(F)  \mid F \text{ compact $A_{+}$-orbits in }\SL_{n+1}(\mathbb{R})/\SL_{n+1}(\mathbb{Z}) \}.$$
\label{th:main_pos}
\end{theorem}
We will only sketch the details of the proof and outline what modifications must be made to get the proof of \cref{th:main_pos} from the proof of \cref{th:main}.

\begin{enumerate}
    \item 
We work with the same orders $\mathcal{O}$ constructed in \cref{eq:defi_of_O}.
Observe that the index of the subgroup $\log \mathcal{O}^{\times}_{+} \subseteq \log \mathcal{O}^{\times}$ is at most $2^{n}$. 
Hence if $\mathcal{O}$ varies along some sequence of orders such that $\sh(\mathcal{O})$ converges to a point, up to varying the orders $\mathcal{O}$ along subsequence of orders, we can assume that $\sh_{+}(\calO)$ also converges to some point. 
\item 
Therefore we are only left with showing that when we lift to the suborders $\mathcal{O}_{(k)} \subseteq \mathcal{O}$
defined in \cref{eq:orders},
we get that the points $\sh_{+}(\mathcal{O}_{(k)})$ approximate an $S$-arithmetic diagonal orbit. This requires some modifications in ``shape machine'' which we list below.
\begin{enumerate}
    \item Instead of considering a generating set of units $\{\beta_i\}_{i=1}^n$ of $\calO^{\times}$ in  \cref{prop:separation-condition}, we will work with a generating set $\{\beta_i^{\plus}\}_{i=1}^n$ of $\calO^{\times}_{\plus}$. The finite set of matrices $\widetilde{M}$ may have to be enlarged to ensure that for some $M \in \widetilde{\calM}$ one gets
   \begin{equation}
       b_i^{2} =  \prod_{j=1}^{n} (\beta_j^{\plus})^{M_{ij}}.
   \end{equation}
   where $b_i = X-a_i$ is a canonized set of units (cf. Definition~\ref{def-canonical-M-basis}). 
   
  \item We generalize  \cref{prop:structure-suborders} by starting with $\beta_1^{\plus} , \cdots , \beta_n^{\plus}$ instead of $\beta_1,\cdots, \beta_n$. With the updated definition of $\widetilde{\calM}$, the conclusion remains the same with the starting from a positive basis of units and exploiting depth. Hence, for varying $k \leq N$ one has 
  \begin{equation}
      \calO^{\times}_{\plus} \cap (\Z + p^{k}\calO_a) 
      = 
      \{ (\beta_1^{\plus})^{x_1} \cdots (\beta_n^{\plus})^{x_n} : (x_1,\cdots,x_n) \in \Lambda_M \},
  \end{equation}
  where $\Lambda_M$ is as given in the system of equations in  \ref{prop:structure-suborders}.
  
  To push through \cref{prop:structure-suborders}, we can work with a larger lower bound $C_n$ (with $p > C_n$) to accommodate the enlarged $\widetilde{\mathcal{M}}$. 
  Then $p \nmid 2 \det M$ and the same idea of taking the adjoint of $M$ pushes through.
  
\end{enumerate}
\end{enumerate}

\bibliographystyle{alpha}
	
\bibliography{references}

\end{document}